\documentclass{Seip-axSIH}

\newtheorem{theorem}{Theorem}
\newtheorem{proposition}[theorem]{Proposition}
\newtheorem{lemma}[theorem]{Lemma}
\newtheorem{corollary}[theorem]{Corollary}

\theoremstyle{definition}
\newtheorem{construction}[theorem]{\noindent{\font\=cmssi10\Construction}\bf}
\newtheorem{constructions}[theorem]{\noindent{\font\=cmssi10\Constructions}\bf}
\newtheorem{definition}[theorem]{\font\=cmssi10\Definition\bf}
\newtheorem{definitions}[theorem]{\font\=cmssi10\Definitions\bf}
\newtheorem{example}[theorem]{\font\=cmssi10\Example\bf}

\newtheorem{remark}[theorem]{\font\=cmssi10\Remark\bf}

\newtheorem{conventions}[theorem]{\noindent{\font\=cmssi10\Conventions}\bf}

\def\FrS{\roman{F\kern.2mmr\kern.3mmS}}
\def\tiDee^#1_#2{\hbox{\font\=cmssi10\D\kern.4mm}{}^{\,#1}_{\hbox{\font\=cmmi6\#2}}} 
\def\tidee^#1_#2{\hbox{\font\=cmssi10\d\kern.4mm}{}^{\,#1}_{\hbox{\font\=cmmi6\#2}}} 
\def\tdelta^#1_#2{\hbox{\font\=cmtex11\\char'012}{}^{\,#1}_{\hbox{\font\=cmmi6\#2}}} 
\def\ucategory#1{\kern.2mm\underline{\roman{\kern-.2mm#1\kern-.2mm}}\kern.6mm}
\newcommand\bosy{\boldsymbol}
\newcommand\mfrak{\mathfrak}
\def\yBGN#1{\kern0.37mm\lower1mm\hbox{$^{^[}$}\kern-0.33mm^{\hbox{\font\=cmr6\#1}}\kern-0.2mm\lower1mm\hbox{$^{^]}$}} 

\def\Tau{\hbox{\font\=cmmi10 scaled\magstep1\\char'034}\kern0.15mm}
\def\Tihprod{\prod{_{\kern-0.3mm_{\bold{ti}\,}}}}
\def\keltop{{k_{}}_{\hbox{\font\=cmr5\el\kern.2mmt}}} 
\def\kappatv{\kappa\kern.3mm\lower.7mm \hbox{\font\=cmr5\t}\lower.67mm\hbox{\font\=cmr6\v}} 
\def\taukv{\tau\ar{\fiveroman kv}} 
\def\skcap{\kern.15mm\sqcap\kern-2.6mm\raise.3mm\hbox{\font\=cmr5\k}\kern1.3mm} 
\def\cgvcap{\kern.8mm\raise.5mm\hbox{\font\=cmbsy5\\char'144}\kern-.3mm\overline{\raise.475mm\hbox{\font\=cmr5\cgv}}\kern-.3mm\raise.5mm\hbox{\font\=cmbsy5\\char'145}\kern.7mm} 

\def\cder_#1{\kern.15mm\hbox{\font\=cmss10\D}{_{\kern-0.1mm}}_{#1}\,} 
\def\DerSe{\hbox{\font\=cmssi10\D\kern.3mm}\lower.7mm\hbox{\font\=cmr5\Se}\kern.6mm} 
\def\varSe{\hbox{\font\=cmtex11\\char'012}\subtext{Se}\kern.4mm}
\def\varSei^#1Ï{\hbox{\font\=cmtex11\\char'012}_{\kern.2mm\hbox{\font\=cmr5\Se}}^{\kern.6mm{#1}}\kern.4mm} 
\def\vProd{\kern.4mm\text{-\,}\roman P\kern-.7mm\lower.7mm\hbox{\font\=cmr5\rod}} 
\def\vDifq^1_#1{\text{\kern.4mm-\,}\roman{Difq}^{\hbox{\font\=cmr6\kern.3mm\1}}_{\,\boldsymbol#1}} 
\def\vbarDelta_#1{\text{\kern.4mm-\,}\bar\Delta{_{\kern.3mm}}_{\boldsymbol#1}\kern.5mm}
\def\cgbarDelta{\raise1.6mm\hbox{\font\=cmr5\{cg}\!}\bar\Delta\kern.45mm}
\def\rcgbarDelta{\raise1.6mm\hbox{\font\=cmr5\{rcg}\!}\bar\Delta\kern.45mm}
\def\CalDBGN^#1{\Cal D_{\hbox{\font\=cmr5\BGN}}^{\kern.6mm{#1}}\kern.15mm}
\def\CalDSeip^#1{\Cal D_{\hbox{\font\=cmr5\Seip}}^{\kern.6mm{#1}}}
\def\CalCSe0{\Cal C\kern.2mm\lower.7mm\hbox{\font\=cmr5\Se\kern.3mm\font\=cmr6\0}{}}
\def\LCS{\roman{LCS}\kern.4mm}
\def\cgVS{\roman{cgVS}\kern.4mm}
\def\scgVS{\roman{scgVS}\kern.4mm}
\def\basTop{\roman{bT}\lower.65mm\hbox{\font\=cmr6\o}\lower.42mm\hbox{\font\=cmr5\p}\kern.6mm} 
\def\tFunso{\roman{tF}\kern-.2mm\lower.7mm\hbox{\font\=cmr6\so}\kern.7mm} 
\def\tConco{\roman tC\kern-.3mm\lower.75mm\hbox{\font\=cmr6\co}\kern.3mm} 
\def\gFtwoilter{\roman{gF}\kern.2mm\raise1.25mm\hbox{\font\=cmr5\2}\kern-1.65mm\lower.75mm\hbox{\font\=cmr5\ilter}\kern.7mm} %
\def\Concu{C\kern-.2mm_{\font\=cmr6\lower.15mm\hbox{\kern.1mm\cu}}} 
\def\Link{\Cal L\subtext{\fivemath k}\sn}
\def\Lincg{\Cal L\kern.2mm\lower.6mm\hbox{\font\=cmr6\c}\lower.45mm\hbox{\font\=cmr5\g}\kern.3mm} 

\def\Cal{\mathcal}
\def\Eps{\hbox{\font\=cmmi10 scaled\magstep1\\char'017}\kern0.15mm}
\def\Iota{\kern.15mm\hbox{\font\=cmmi10 scaled\magstep1\\char'023}\kern0.2mm}
\def\Nu{\hbox{\font\=cmmi10 scaled\magstep1\\char'027}\kern0.25mm}
\def\uvarPi{\kern.15mm\underline{\kern-.15mm\varPi\kern-.85mm}\kern.85mm}
\def\uOmega{\kern.3mm\underline{\kern-.3mm\Omega\kern-.3mm}\kern.3mm}

\def\sigmaa{\char'033}
\def\tauu{\char'034}
\def\fRe{\hbox{\font\=cmr9\f\kern.1mm}\roman{Re}\kern.75mm}
\def\fIm{\hbox{\font\=cmr9\f\kern.1mm}\roman{Im}\kern.65mm}
\def\vecc#1{\kern-.5mm\vec{\kern.5mm#1}}
\def\TVS{\roman{TVS}\kern0.37mm}%
\def\LCS{\roman{LCS}\kern0.37mm}%
\def\BaS{\roman{BaS}\kern0.37mm}%
\def\dimHa{{\rm dim_{_{\kern.2mm Ha}}}}
\def\rajou{{}^{}{\Cal B}_{s\,}}
\def\Lis{\Cal L\kern.3mmis\,} 
\def\Linb{\Cal L\lower.7mm\hbox{\kern.1mm\font\=cmmi6\b}}
\def\LL^#1{L\kern0.15mm\raise.4mm\hbox{$^{#1}$}\kern0.15mm}

\def\dualbeta{^{\kern0.4mm\prime}_{\kern-.2mm\raise.95mm\hbox{$_{_\beta}$}}} 
\def\Nbh{\Cal N_{\font\=cmmi6\lower.15mm\hbox{\kern.1mm\bh\kern.15mm}}}
\def\Topma{\roman{{Top_{}}_{\hbox{\font\=cmr6\ma}}\kern.15mm}}
\def\co{\hbox{\font\=cmmi12\c}\kern.15mm\lower.15mm\hbox{$_{\rm o}$}}
\def\prodc{\prod{_{_{\kern-.3mm\bold c\kern.15mm}}}}
\def\vsprod_#1_#2{\prod\kern-0.3mm{}_{_{\roman{#1}\sp{#2}\,}}} 
\def\vscoprod_#1_#2{\coprod\kern-0.3mm{}_{_{\roman{#1}\sp{#2}\,}}} 
\def\expnota^#1]_#2{\,^{#1\,]{_{}}_{\roman{#2}}}} 

\def\bold#1{{\bf#1}}
\def\roman#1{{\rm#1}}
\def\limu_#1{\lim\kern-5.5mm\lower1.5mm\hbox{$_{#1}\ $}}
\def\oseoy{\raise1.9mm\hbox{\kern.5mm\font\=cmr5\o}\kern-1.7mm y}
\def\O{{}^{}\Cal O}
\def\Univ{\hbox{\font\=cmssbx10\U}} 
\def\Pows{\Cal P\kern-.4mm_s\kern.3mm}
\def\lei{      {}_{ {}^{\,\downarrow\text{\hskip-2.1mm}       }  }  \cap       }
\def\lei{\hbox{\kern.45mm$_{^\downarrow}\kern-1.280mm\cap\kern.85mm$}}
\def\Zepp{{{{Z\!\!\!Z^{\phantom{l}}}^{{}_{{}^{\!}\!+}}}}}  

\def\inve{\lower.85mm\hbox{$^{^-}$}\kern-.5mm{}^\iota}

\def\fvalue{\hbox{\kern.2mm\font\=cmr10\\char'022\kern-.2mm}} 
\def\ffvalue{\hbox{\kern.2mm\font\=cmr7\\char'022\kern-.2mm}} 
\def\image{\hbox{\font\=cmr10\\char'022\kern-1mm\char'022}} 
\def\iimage{\hbox{\font\=cmr7\\kern.3mm\char'022\kern-.7mm\char'022\kern-.3mm}} 
\def\images{\hbox{\font\=cmr10\\char'022\kern-1mm\char'022\kern-1mm\char'022}} 
\def\weco{\kern.15mm\hbox{\font\=cmtt10\\char'054}\kern.4mm} 
\def\cdotn{\kern-.2mm\cdot\kern-.2mm} 
\def\timesn{\kern-.2mm\times\kern-.2mm} 
\def\ttimes{\hbox{\kern-.2mm${}\times\kern-2.5mm\lower.8mm\hbox{\font\=cmr5\t}\kern1.8mm$}} 
\def\ttimesn{\hbox{\kern-.2mm${}\times\kern-2.5mm\lower.8mm\hbox{\font\=cmr5\t}\kern1.4mm$}} 
\def\ktimes{\hbox{\kern-.2mm${}\times\kern-2.5mm\lower1mm\hbox{\font\=cmr5\k}\kern1.5mm$}} 
\def\vstimes{\kern.95mm\raise.45mm\hbox{\font\=cmbsy6\\char'002}\kern-2.3mm\lower.9mm\hbox{\font\=cmr5\vs}\kern1.05mm} 
\def\risti#1{{}^{}\,{\times}_{{\!}_{{}_{\!\!\!{#1}{}^{\!}\ }}}}
\def\Circ{\kern.9mm\hbox{\font\=cmbsy10\\char'016}\kern.9mm}
\def\cardplus{\hbox{$\kern.77mm+\kern-1.95mm\raise.23mm\hbox{$_{_{\roman c}}$}\kern1.33mm$}}%
\def\ordplus{\hbox{$\kern.78mm+\kern-1.97mm\raise.23mm\hbox{$_{_{\roman o}}$}\kern1.22mm$}}%
\def\Examplee{{\font\=cmssi10\E\kern.15mmx\kern.15mma\kern.15mmm\kern.14mmp\kern.17mml\kern.15mme}\kern.3mm. }
\def\Examples{{\font\=cmssi10\E\kern.15mmx\kern.15mma\kern.15mmm\kern.14mmp\kern.17mml\kern.15mme\kern.15mms}\kern.3mm. }
\def\Remarkk{{\font\=cmssi10\R\kern.15mme\kern.15mmm\kern.15mma\kern.15mmr\kern.15mmk\kern.15mm. }}
\def\Remarkss{{\font\=cmssi10\R\kern.15mme\kern.15mmm\kern.15mma\kern.15mmr\kern.15mmk\kern.15mms}\kern.3mm. }
\def\N{{I\!\!N}} 
\def\No{{I\!\!N\kern-.54mm\lower.15mm\hbox{$_{\rm o}$}}} 
\def\iNo{I\!\!{N_{}}_{\kern-.22mm{\rm o}}} 
\def\Nopot#1{I\!\!N\kern-.54mm\lower.15mm\hbox{$_{\rm o}$}\kern-.7mm{}^{#1}} 
\def\potNo{^{\kern.37mm I\!\!{N_{}}_{\kern-.22mm{\rm o}}}} 
\def\minus{\kern.2mm\lower1.05mm\hbox{$^-$}}
\def\pplus{\raise.22mm\hbox{\font\=cmr5\\char'053}}
\def\mminus{\raise.18mm\hbox{\font\=cmsy5\\char'000}}
\def\plusinftyy{\raise.18mm\hbox{\font\=cmr5\\char'053}\infty}
\def\minusinftyy{\raise.18mm\hbox{\font\=cmsy5\\char'000}\infty}
\def\plusinfty{\lower1.05mm\hbox{$^+$}\infty}
\def\minusinfty{\lower1.05mm\hbox{$^-$}\infty}
\def\Qe{\hbox{$Q\kern-2.6mm\raise.2mm\hbox{\font\=cmssqi8\I}\kern1.7mm$}}
\def\Re{I\!\!R}
\def\Rep{{{I\!\!R^{\phantom{l}}}^{{}_{{}^{\!}\!+}}}}

\def\Ce{{\hbox{$C\kern-2.5mm\raise.2mm\hbox{\font\=cmssqi8\I}\kern1.48mm$}}}
\def\imag{\kern.15mm\lower.6mm\hbox{$^{^*}$}\kern-1.8mm\imath\kern.1mm} 
\def\Ke{{{}^{}I\!\!K}}
\def\ebiF{\kern.1mm\hbox{\font\=cmmib8\F}\kern.5mm} 
\def\ebiT{\kern.1mm\hbox{\font\=cmmib8\T}\kern.6mm} 
\def\ebiU{\kern.1mm\hbox{\font\=cmmib8\U}\kern.5mm} 
\def\biit#1{\hbox{\font\=cmmib10\#1}} 

\def\bmii#1#2{\hbox{\font\=cmmib#1\#2}}

\def\fssi#1{\hbox{\font\=cmssi10\#1}\kern0.15mm} 
\def\smb#1{\hbox{\font\†=cmmi8\†#1\kern.3mm}} 
\def\bCal#1{\hbox{\font\=cmbsy10\#1}} 
\def\eCal#1{\kern.1mm\hbox{\font\†=cmbsy8\†#1\kern.4mm}} 
\def\ecal#1{\kern.1mm\hbox{\font\†=cmsy8\†#1\kern.3mm}} 
\def\ncal#1{\kern.1mm\hbox{\font\†=cmsy9\†#1\kern.3mm}} 
\def\vcal#1{\kern-.1mm\vec{\kern.2mm\hbox{\font\†=cmsy7\†#1}\kern.3mm}} 

\def\concc{{}^{}{}^{}\!{\bf{\hat{\phantom w}}}\!{}^{}}
\def\id{\kern.3mm\roman{id}\kern.7mm}
\def\idv{\hbox{\font\=cmr10\id}\kern.25mm\lower.8mm\hbox{\font\=cmr7\v}\kern.3mm} 
\def\idm{\hbox{\font\=cmr10\id}\kern.25mm\lower.8mm\hbox{\font\=cmr6\m}\kern.3mm} 
\def\seq#1{\langle#1\rangle}
\def\Seq#1{\big\langle#1\big\rangle}
\def\ymp{{}^{}{\Cal N}_o\,}
\def\SemiNor{\Cal S_{_N}\kern0.15mm}
\def\vecs{\upsilon\kern-0.3mm\lower.15mm\hbox{$_s$}\kern0.2mm} 
\def\vecss{\hbox{\font\=cmitt10\v}\kern-0.1mm\lower.15mm\hbox{$_s$}\kern0.2mm} 
\def\bnull#1{\hbox{\font\=cmssbx10\0}{}_{\font\=cmmi6\lower.15mm\hbox{\kern-.1mm\#1\kern.15mm}}} 
\def\bnulla#1_#2{\hbox{\font\=cmssbx10\0}{}_{\font\=cmmi6\lower.15mm\hbox{\#1\kern-.1mm}}\lower.3mm\hbox{$_{_{#2}}$}} 
\def\bzero#1{\hbox{\font\=cmbx10\0}{}_{\font\=cmmi6\lower.15mm\hbox{\kern-.1mm\#1\kern.15mm}}} 
\def\dom{{{}^{}{\rm dom}\,{}_{{}^{}}}}
\def\domm{\kern0.15mm{\rm dom}{^{\kern.3mm\hbox{\font\=cmr6\2}}}\,}
\def\domr#1{\roman{dom}^{\font\=cmr6\raise.0mm\hbox{\kern.3mm\#1}}}

\def\rng{{}^{}{\rm rng}\,{}_{{}^{}}}
\def\CPi#1{C\kern-.2mm\lower.05mm\hbox{$_{_\Pi}$}\kern-1.52mm{}^{#1}}
\def\CinftyPi{C\kern.4mm\raise.3mm\hbox{$^\infty$}\kern-3.35mm_{_\Pi}\kern1.45mm}
\def\CinftyS{\Cinfty\kern-3.9mm_{_{\Cal S}}\kern1.45mm}
\def\Cinfty{C\kern.4mm\raise.3mm\hbox{$^\infty$}\kern.15mm}
\def\Cinftyzero{\hbox{$C\kern.4mm\raise.3mm\hbox{$^\infty$}\kern.15mm\kern-3.5mm_{\font\=cmr6\lower.15mm\hbox{\kern.1mm\0}}\kern1.9mm$}}

\def\RHB#1#2{\raise#1mm\hbox{$#2$}} 
\def\LHB#1#2{\lower#1mm\hbox{$#2$}} 
\def\fivemath#1{\hbox{\font\=cmmi5\#1}}
\def\fiveroman#1{\hbox{\font\=cmr5\#1\kern.1mm}}
\def\sixmath#1{\hbox{\font\=cmmi6\#1}}
\def\sixroman#1{\hbox{\font\=cmr6\#1\kern.1mm}}
\def\eightmath#1{\hbox{\font\=cmmi8\{#1}\kern.1mm}}
\def\eightroman#1{\hbox{\font\=cmr8\{#1}\kern.1mm}}
\def\erm#1{{\font\=cmr8\#1}}
\def\subtext#1{\raise.2mm\hbox{$_{_{\kern0.15mm\roman{#1}}}$}}
\def\subtexT#1{\raise.2mm\hbox{$_{_{\kern0.15mm\hbox{\font\=cmr5\#1}}}$}}
\def\sNor#1{\kern.25mm\lower.38mm\hbox{$_{#1}$}}
\def\sNorr#1{\kern-.2mm\lower.38mm\hbox{$_{#1}$}}
\def\sNoreset_#1{\kern.13mm\lower.83mm\hbox{\font\=cmmi6\C}\kern.32mm\lower.1mm\hbox{$_{^{\emptyset,#1}}$}}
\def\sbi#1{{_{\kern-0.1mm}}_{#1}} 
\def\ais#1_#2{{}_{\font\=cmmi6\lower.15mm\hbox{\kern-.1mm\#1\kern.15mm}}\lower.3mm\hbox{${_{\kern-0.3mm_{#2}}}$}} %
\def\aais#1_#2{\kern.1mm{}_{\font\=cmmi6\lower.25mm\hbox{\kern-.1mm\#1\kern.15mm}}\lower.4mm\hbox{${_{\kern-0.3mm_{#2}}}$}} %
\def\ai#1{{}_{\font\=cmmi6\lower.15mm\hbox{\kern-.1mm\#1\kern.15mm}}} 
\def\yi#1{^{\font\=cmmi6\raise.0mm\hbox{\kern-.1mm\#1\kern.15mm}}} 
\def\ear#1{{}_{\font\=cmr5\lower.15mm\hbox{\kern.1mm\#1}}} 
\def\ar#1{{}_{\font\=cmr6\lower.15mm\hbox{\kern.1mm\#1}}} 
\def\aar#1{_{\font\=cmr6\lower.15mm\hbox{\kern.1mm\#1}}} 
\def\yr#1{^{\font\=cmr6\raise.0mm\hbox{\kern.3mm\#1}}} 
\def\yrai^#1_#2{^{\kern.4mm\hbox{\font\=cmr6\{#1}}}_{\kern.2mm{#2}}}
\def\upparentes#1{^{\kern.2mm\raise.2mm\hbox{\font\=cmr6\\char'050}\kern.1mm{#1}\kern.1mm\raise.2mm\hbox{\font\=cmr6\\char'051}}} 
\def\lupar{\kern.2mm\lower1mm\hbox{$^{^(}$}} 
\def\rupar{\lower1mm\hbox{$^{^)}$}\kern-.15mm} 
\def\yyi#1{^{\font\=cmmi6\lower.6mm\hbox{\kern-.25mm\#1\kern-.05mm}}} 
\def\yyr#1{^{\font\=cmr6\lower.45mm\hbox{\kern-.25mm\#1\kern-.15mm}}} 
\def\yplus{\lower1mm\hbox{$^{^+}$}} 
\def\yminus{\lower1mm\hbox{$^{^-}$}} 
\def\aminus{{\kern.15mm\raise.3mm\hbox{$_{_-}$}\kern-.1mm}}%
\def\yvee{\LHB{.9}{^{^{\,\vee}}}\kern-.3mm} 
\def\ywed{\LHB{.9}{^{^{\,\wedge}}}\kern-.3mm} 
\def\yplk{{}^{{}_{{}^{'\!}}}} 
\def\adot{\kern.2mm\hbox{\font\=cmb10\\char'056}}%
\def\ydot{\kern.2mm\raise1.9mm\hbox{\font\=cmb10\\char'056}}
\def\yydot{\kern.2mm\raise1.35mm\hbox{\font\=cmb7\\char'056}\kern.2mm}
\def\yydott{\kern.2mm\raise1.35mm\hbox{\font\=cmb6\\char'056}\kern.2mm}
\def\ClT{{\rm Cl}\kern.25mm\lower.4mm\hbox{$_{\Cal T}$}\kern0.2mm} 
\def\IntT{\sp{\rm Int}\kern.2mm\lower.4mm\hbox{$_{\Cal T}$}\kern0.2mm} 
\def\Cl_taurd#1{\roman{Cl_{}}_{\kern0.37mm\hbox{\font\=cmmi8\\char'034}\kern-0.15mm{_{}}_{rd}\kern0.2mm#1\,}}
\def\Int_taurd#1{\roman{Int_{}}_{\kern0.37mm\hbox{\font\=cmmi8\\char'034}\kern-0.15mm{_{}}_{\roman{rd}}\kern0.2mm#1\,}}
\def\inc{\subseteq}
\def\iinc{\supseteq}
\def\exi#1{\exists\,#1\kern.2mm\,;}
\def\all#1{\forall\,#1\kern.2mm\,;}
\def\imply{\Rightarrow}
\def\equivv{\Leftrightarrow}

\def\spp{\kern0.07mm} 
\def\sp{\kern0.15mm} 
\def\ssp{\kern0.37mm} 
\def\snn{\kern-0.2mm} 
\def\sn{\kern-0.3mm} 
\def\ssn{\kern-0.63mm} 
\def\biggerlineskip#1 {\linebreak\nopagebreak\vskip-4.2mm\vskip.#1mm\nopagebreak\noindent}%
\def\Biggerlineskip#1 {\linebreak\nopagebreak\vskip-4.2mm\vskip#1pt\nopagebreak\noindent}%
\def\KP#1{\kern#1mm} 
\def\KN#1{\kern-#1mm} 
\def\nhskip#1mm{$\null$\kern#1mm}
\def\hyppy#1{$\phantom{}$\hskip#1}
\def\mhyppy#1{\null\kern#1mm}

\def\text#1{\hbox{\rm#1}}
\def\NS{\vskip1.7mm}

\def\VBOX/#1/#2/HEREend{\vbox{#2\vskip-#1mm}\vfill\null\eject}
\def\œ$#1${\hbox{$#1$}} 
\def\"{\"a} \def\"{\"o}
\def\q#1{``\kern0.37mm#1\kern0.37mm"}
\def\newProCla#1\par#2\par{\vskip1.7mm\noindent\bf#1\it#2\vskip1.7mm}
\def\Prooff{{\font\=cmssi10\P\kern.37mmr\kern.37mmo\kern.37mmo\kern.37mmf\kern.37mm. }\rm}

\def\QED{\hfill\hbox{$\ \sqcap\kern-2.45mm\sqcup$}}

\def\noin{\noindent}
\def\Newline{\kern-10mm\newline}
\font\rp=cmr8
\def\eps{\varepsilon}
\def\leu{\raise1.5mm\hbox{\font\=cmmi5\\char'074}\kern.2mm}%
\def\riu{\kern.2mm\raise1.5mm\hbox{\font\=cmmi5\\char'076}}%
\def\Symbol#1Ï{\kern.35mm\hbox{\font\=cmr10\\char'047}\kern.2mm#1\kern.35mm\hbox{\font\=cmr10\\char'047}}%
\def\Symboo#1Ï{\kern.35mm\text{`}\kern.2mm#1\kern.35mm\hbox{\font\=cmr10\\char'047}}
\def\RunMyHead#1#2#3#4{%
 \headline{\ifnum\pageno=\firstpage\hfil%
           \else{\ifodd\pageno{\rp#3\phantom\folio\hfil#4\hfil\phantom{#3}\folio}%
                 \else{\rp\folio\phantom{#2}\hfil#1\hfil\phantom\folio#2}%
                 \fi}%
           \fi}%
 \footline{\ifnum\pageno=\firstpage\hfil{\rp[\,\folio\,]}\hfil%
           \else\hfil%
           \fi}%
}%
\def\bulgin{\noindent$\bullet$ \ \kern.1mm} 
\def\bulgen{\noindent\kern-1.5mm$\bullet$\kern3.95mm} 
\def\subhead#1\par#2\par{\vskip4mm\smallbreak\null\smallskip\vbox{\noindent\bbf#1\hfill\kern1.5mm#2\hfill\phantom{#1}\vskip2.5mm\nopagebreak}\nopagebreak\noindent}
\def\subheadd#1\par#2\par#3\par{\vskip4mm\smallbreak\null\smallskip\vbox{\noindent\bbf#1\hfill#2\hfill\phantom{#1}\vskip1.5mm\centerline{#3}\vskip2.5mm\nopagebreak}\nopagebreak\noindent}
\def\insubsubhead#1\par{\vskip4mm$\null$\hskip2mm{\font\=cmss10\#1}\vskip2mm\noindent}%
\def\binsubsubhead#1#2\par{\vskip4mm{\bf#1.}\hskip5mm{\font\=cmss10\#2}\nopagebreak\vskip2mm\nopagebreak\noindent}%
\def\wave{\hbox{\font\†=cmsy10\†\hbox{\char'164}\kern-2.35mm\hbox{\char'165}\kern.55mm}}
\def\wavee{\hbox{\font\†=cmsy8\†\hbox{\char'164}\kern-2.0mm\hbox{\char'165}\kern.4mm}} 
\def\barmj{\kern.25mm\bar{\hbox{\font\=cmr10\\char'021}}\kern.4mm}

\def\sigrd{\sigma\kern-.2mm\lower.7mm\hbox{\font\=cmr6\r\font\=cmr5\d}\kern.6mm}
\def\ssigrd{\sigma\kern-.2mm\lower.7mm\hbox{\font\=cmr6\r\font\=cmr5\d}\kern-1.7mm\raise1.25mm\hbox{\font\=cmr6\2}\kern1mm}
\def\sssigrd{\sigma\kern-.2mm\lower.7mm\hbox{\font\=cmr6\r\font\=cmr5\d}\kern-1.7mm\raise1.25mm\hbox{\font\=cmr6\3}\kern1mm}
\def\taurd{\tau\kern-.4mm\lower.7mm\hbox{\font\=cmr6\r\font\=cmr5\d}\kern.6mm}
\def\tsigrd{\tau\sigma\kern-.2mm\lower.7mm\hbox{\font\=cmr6\r\font\=cmr5\d}\kern.6mm}
\def\tauR#1{\tau_{_{I\!\!R}}\kern-1.5mm^{#1}}
\def\RN{I\!\!R\kern.3mm^{\hbox{\font\=cmmi6\N}}} 
\def\QTN{Q\kern.1mm_{\lower.2mm\hbox{\font\=cmmi6\T}}^{\kern.2mm\hbox{\font\=cmmi6\N}}} 
\def\leLCSr{{\le}{}_{_{{\rm LCS}}}}
\def\leLCS-{{\le}{}_{_{{\rm LCS}}}\text{\sp-\sp}}
\def\Centerline#1\par#2\par#3{\noindent#1\phantom{#3}\hfill#2\hfill\phantom{#1}#3}

\begin{document}

\title[$\text{\sc Seip's differentiability as BGN}$]%
           {Seip's differentiability concepts as a particular case\vskip1mm
             of the Bertram\,--\,Gl\"ckner\,--\,Neeb construction}

\author[S. Hiltunen]{Seppo\ I\. Hiltunen}
\address{Helsinki University of Technology                             \vskip0mm$\hspace{2mm}$
           Institute of Mathematics, U311                              \vskip0mm$\hspace{2mm}$
           P.O.\ Box 1100                                              \vskip0mm$\hspace{2mm}$
           FIN-02015 HUT\vskip0mm
         FINLAND}
\email{shiltune\,@\,cc.hut.fi}



\begin{abstract}

From the point of view of unification of differentiation theory, it is of
interest to note that the general construction principle of Bertram, Gl\"ckner
and Neeb leading to a\linebreak $C^{\ssp k}$ differentiability concept from a
given $C^{\,0}$ one, besides subsuming the Keller\,--\,Bastiani 
$C_{\roman c}^{\ssp k}$ differentiabilities on real Hausdorff locally convex
spaces, also does the same to the \q{arc-generated} interpretation of the
Lipschitz theory of differentiation by Fr\"licher and Kriegl, and likewise to
the \q{compactly generated} theory of Seip's continuous differentiabilities.
In this article, we give the details of the proof for the assertion concerning
Seip's theory. We also give an example indicating that the premises in Seip's
various inverse and implicit function theorems may be too strong in order for
these theorems to have much practical value. Also included is a presentation
of the BGN\,--\,setting reformulated so as to be consistent with the
Kelley\,--\,Morse\,--\,G\"del\,--\,Bernays\,--\,von Neumann type approach to
set theory, as well as a treatment of the function space constructions and
development of their basic properties needed in the proof of the main result.

  \end{abstract}

\maketitle


        \insubsubhead{Introduction and some preliminaries}

According to the introductions in \cite{Se72} and \cite{Se79}\ssp, Seip's main
motivation for the development of his theory of differentiation was to
establish a setting with a \q{purely} topological basis where function spaces
also beyond Banach spaces could be considered as domains and ranges of
differentiable maps, and where also some kind of \q{cartesian closedness}
holds so that \q{exponential laws} would be available to facilitate proving
smoothness of given maps. Seip grounded his theory on the observation by
Gabriel and Zisman \cite[p.\ 47]{GZ} that considering the k-extension of the
compact open topology for the set of continuous functions between topological
spaces gives rise to a cartesian closed category.

A bit less imprecisely, the last assertion means the following. For a
topological space $X$ let $kX$ denote the space with the same underlying set
equipped with the finest topology such that the identity restricted on every
compact set is continuous \œ$X\to kX\sp$. Call $X$ {\it compactly generated\,}
if{}f it is Hausdorff with \œ$kX=X\sp$. Let \œ$X\sqcap Y=kP$ when $P$ is the
usual product space of $X$ and $Y\spp$, and let \œ$Y^{\,X}=kC$\linebreak when
$C$ is the compact open topological space of all continuous functions \œ$X\to
Y\spp$. Then \œ$u\mapsto\hat u\ssp$, where \œ$\hat u:(\sp x\ssp,y\sp)\mapsto
u\ssp(y)\ssp(x)\,$, defines a bijection \œ$\big(Z^{\,X}\sp\big)\sp^Y\to
Z\,^{X\,\sqcap\,Y}$ for any compactly generated $X\sp,Y\spp,Z\sp$.

In order to be still less imprecise, note that actually we do not have a
bijection as just written, but from the {\it underlying set\ssp} of $
\big(Z^{\,X}\sp\big)\sp^Y$ onto that of $Z\,^{X\,\sqcap\,Y}\sp$, that is
\œ$\ssp\sigrd\sn\big(\big(Z^{\,X}\sp\big)\sp^Y\sp\big) \to
\sigrd\sn\big(Z\,^{X\,\sqcap\,Y}\sp\big)\,$. However, it follows that the
same function \œ$u\mapsto\hat u$ is a homeomorphism \œ$
\big(Z^{\,X}\sp\big)\sp^Y\to Z\,^{X\,\sqcap\,Y}\sp$.

In Example \ref{Seip d:fm group} below, the use of this kind of exponential
laws \œ$\big(Z^{\,X}\sp\big)\sp^Y\approx Z\,^{X\,\sqcap\,Y}$ in Seip's theory
of differentiation is exemplified by giving a simple proof of Seip\,--\,%
smoothness \œ$E\sp\sqcap E\to kE$ of the map \œ$(\sp x\ssp,y\sp)\mapsto x\circ
(\sp\Iota+y\sp)\,$, where we have the locally convex test function space $E =
\Cal D\ssp(\Re\sp)\sp$ with $\ssp\Iota=\id\Re\sp\,$.

It is remarkable that the method of Example \ref{Seip d:fm group} applies
although the locally convex space $E\sp$ is neither a canonical function space
in Seip's theory, nor compactly generated by \cite[Theorem 6.1.4\ssp(iii)\ssp,
Proposition 6.2.8\ssp(ii)\ssp, pp.\ 190, 195]{FK}\ssp. Note that by a {\it%
canonical function space\ssp} in a theory of differentiation we mean any
prearranged object of the theory which has as its underlying set the set of
all order $i$ differentiable functions \œ$f:E\iinc U\to F$ for some fixed $
E\ssp,F\sp,U\spp,i\ssp$. Having intrinsic exponential laws in a theory
requires canonical function spaces. If there are no such, possible exponential
laws have to established in an ad hoc manner, as for example in (new)
particular cases of the general theory developed in \cite{BGN}\ssp.

However, it should be noted that although exponential laws may provide easy
proofs of smoothness of maps between spaces of smooth functions, the same does
not hold for maps between spaces of finite order differentiable functions.
Since practical inverse and implicit function theorems may require some kind
of use of Banach spaces, and since spaces of smooth functions seldom are such,
we see that the goals of on one hand possessing exponential laws, and on the
other hand having available usable inverse or implicit function theorems may
be somewhat contradictory.

Indeed, for example in \cite[pp.\ 55\,--\,57, 73, 80, 81, 92, 93]{Se72}\ssp,
inverse and implicit function theorems were provided. However, they were so
formulated that their proofs within the theory became almost trivial, this
having the consequence that verifying their presuppositions in practice
became almost impossible. For instance, if with \œ$I=[\,0\,,1\,]$ we consider
the diffeomorphism \œ$f:\Cinfty(I\sp)=G\to G$ defined by \œ$x\mapsto\varphi
\circ x$ for a fixed nonaffine smooth diffeomorphism \œ$\varphi:\Re\to\Re\sp\,
$, we shall see in Example \ref{scharf d:bar} below that $f$ {\it is not
scharf differenzierbar\ssp} at any constant point $x =I\snn\times\sn\{\sp\xi
\sp\}$ when $\xi\in\Re$ is such that $\varphi\ssp''(\xi)\not=0\,$.

In the literature, there do not seem to be any serious applications of Seip's
theory of differentiation. Besides the facts given above, one possible reason
for Seip's theory not having become popular may be the overwhelmingly
\q{categorical} style of presentation in \cite{Se72} and \cite{Se79} where
also things more simply expressible without any notions of category theory
have been stated in such terms. There also are some obscurities in the
notations for the various categories.

A good example of making simple matters obscure and complicated by jargonizing
them category-theoretically is \cite[Definition 8.6, p.\ 108]{Se72} where a {\it
norm\ssp} in a real locally convex Hausdorff space $E$ with topology $\sp
\Cal T\sp$ just means any \œ$n\in\Cal C\,^A$ satis\-fy\-ing \œ$z\in n\,(z)\inc
\ClT(\sp n\,(x)+n\,(\sp z-x\sp))$ for all \œ$x\ssp,z\in A$ when $A$ is the set
of vectors of $E$ and $\Cal C$ is the set of bounded absolutely convex closed
sets in $E\,$. A simple example of this kind of (rather useless\,?\sp)\ %
\q{norm} is \œ$A\owns x\mapsto\{\,t\,x:\minus 1\le t\le 1\ssp\}\,$. See also
Example \ref{"norm"} below for some more details on this matter.

Seip's theory of continuous differentiabilities in \cite{Se72} and \cite{Se79}
is one in the line of the theories of differentiation where some canonical
function spaces are incorporated in the theory wishing to obtain certain
intrinsic exponential laws. In these theories, there generally are no natural
inverse and implicit function theorems. Other examples in this line can be
found in \cite{Ba}\ssp, \cite{FB}\ssp, \cite{FK} and \cite{Nel90}\ssp.

In the other line, one has theories more adapted to handling inverse or
implicit function theorems, but function spaces generally have to be treated
in an ad hoc manner. One example is the Banach space calculus around the
classical continuous differentiabilities sketched in
\cite[pp.\ 147\,ff., 181\,--\,182]{Di}\ssp. Other examples are our
modification in \cite[pp.\ 237\,--\,241]{Hi1} of \cite{FB}\ssp, generalizing
\cite{Di}\ssp, and the particular cases in \cite{Helge-imp} of \cite{BGN}
which also generalizes \cite{Se72} and \cite{FK}\ssp, and a portion of
\cite{Ba}\ssp. Note also \cite[pp.\ 3\,--\,6, 22\,ff.]{Yam79}\ and
\cite[pp.\ 73\,ff., 140\,--\,143]{Ham82} where quite special differentia-
bility concepts are designed hoping to get certain \q{intrinsic} inverse and
implicit function theorems applicable e.g.\ to some maps of Fr\'echet function
  spaces.

In this article, we concentrate on proving that indeed Seip's theory can be
obtained as a particular case of the general construction in \cite{BGN}\ssp.
We hope to give the proof of the corresponding assertion associated with the
theory in \cite{FK} later.

To one acquainted with the developments in \cite{Ba}\ssp, \cite{Se72} and
\cite{BGN}\ssp, because of the manner the higher order differentiabilities are
constructed in the last, it should not be a surprise that portions of the
theories in the former are obtainable as particular cases of the construction
in \cite{BGN}\ssp.

Namely, one gets the higher order differentiabilities in \cite{BGN} by a
recursion using the difference quotient family associating with a continuous \œ$
f:E\to F$ the map \œ$
f\yBGN1\ssn:\text{``\ssp}E\times E\times\mathbb K\text{\,"}\to F$ given by \œ$
(\sp x\ssp,u\ssp,t\sp)\mapsto
  t^{\sp\mminus 1}\sp(\sp f\sp(\sp x+t\,u\sp) - f\sp(x))$ for $t$ invertible.
In \cite{Ba} and \cite{Se72}\ssp, one instead uses the derivative \œ$f\ssp'\sn
:E\to L$ given by \œ$f\ssp'(x)\,u = f\yBGN1(\sp x\ssp,u\ssp,0\sp)$ with $L$ a
space of continuous linear maps \œ$E\to F$ having the exponential law property
that a map \œ$E\to L$ is continuous if{}f the associated map \œ$
\text{``\ssp}E\times E\text{\sp\,"}\to F$ is such. Differentiability being
defined with the aid of a remainder condition, cf.\
\cite[D\'efinition 2.1, p.\ 39]{Ba} and
\cite[Definition 4.2, p.\ 60]{Se72}\ssp, requiring continuity of the map \œ$
(\sp x\ssp,u\ssp,t\sp)\mapsto f\yBGN1(\sp x\ssp,u\ssp,t\sp)-f\ssp'(x)\,u\ssp$,
it is expectable that the different approaches lead to the same concepts.

Although the results are not surprising, it is surprising that the proofs,
when properly presented, are not at all straightforward or trivial. This may
be considered as a consequence of the different modes of development in the
theories.

The contents of the present article consisting of a single section divided
into this introduction and four subsections is roughly described by the
following list of the titles of these subsections:\vskip.5mm

A \ \ A reformulation of the general BGN\,--\,setting \KP2\dotfill\KP4 p.\ \ \,\pageref{subsec A}\KP7

B \ \ Compactly generated topologies, vector and function spaces \hfill p.\ \pageref{subsec B}\KP7

C \ \ Riemann integration of curves in topologized vector spaces \hfill p.\ \pageref{subsec C}\KP7

D \ \ Seip's higher order differentiability classes \KP2\dotfill\KP4 p.\ \pageref{subsec D}\KP7\vskip.5mm

\noin In A we reformulate the basic setting in \cite{BGN} in order to get it
accordant with the logically economical and precise notational system we
followed in \cite{Hic} and to be complemented below. In B we give the basic
definitions and establish the basic facts concerning general compactly
generated (vector) spaces and spaces of continuous functions needed in the
sequel. In C we give a short account of the matter in the title to make things
precise. In D we establish the main result of this article, given as the
conjunction of Theorems \ref{Seip = BGN} and \ref{our equivv Seip} below.

At the end of D we have included two examples of which
Example \ref{Seip d:fm group} shows how the use of exponential laws in Seip's
theory can be utilized to get simple proofs of smoothness of certain kind of
maps. The purpose Example \ref{scharf d:bar} is to raise the question whether
any serious applications of Seip's various inverse and implicit function
theorems are possible.                                                   \NS

In the above introduction, we used some notations more or less informally,
which we abandon from now on. For example, we above let \q{$u\ssp(x)$} denote
the function value of $u$ at $x\ssp$. From now on, it is $u\fvalue x\ssp$.
Precisely, for any classes $\sp u\ssp,x\sp$ we define \œ$u\fvalue x =
\bigcap\ssp\{\,y:\all z\,(\sp x\ssp,z\sp)\in u\equivv y=z\,\}\,$ deviating
from \cite[Definition 68, p.\ 261]{Ky} in order to have \œ$\sp R\sp\fvalue x =
\Univ\ssp$ also e.g.\ if \œ$\ssp(\sp x\ssp,y\ssp)\ssp,(\sp x\ssp,z\ssp)\in R\sp
$ with \œ$\ssp y\not=z\ssp$. However, note that $\ssp u\fvalue x =
u\ssp(x)\subtext{Kelley}$ for any {\it function\ssp} $u$ and any class $x\ssp$.

Likewise, the symbol $\Symboo\text{\sp$\sqcap$\spp}Ï$ no longer refers to any
compactly generated product. Instead, we define \œ$E\sp\sqcap F =
((\sp\sigrd E\sp)\vstimes(\sp\sigrd F\sp)\ssp,
(\sp\taurd E\sp)\ttimesn(\sp\taurd F\sp))$ so that (see below) if $E\ssp,F$
are topologized $\sp\bold K\,$--\,modules, then \œ$E\sp\sqcap F$ is the
product of the underlying algebraic modules equipped the usual Tihonov product
of the respective topologies. In particular, this applies when $E\ssp,F$ are
real topological vector spaces.

Besides the preceding two ones, from now on we generally use also the other
notational conventions of \cite{Hic}\ssp. In addition to these, we give the
following

\begin{conventions}

Working in a set theoretical setting closely (but not exactly) parallelling
the presentation in \cite[Appendix, pp.\ 250\,--\,281]{Ky}\ssp, see also
\cite{Mo}\ssp, and defining \œ$\Univ=\{\,x:x=x\,\}$ and \œ$\{\sp x\sp\} =
\{\,z:\all y\,x\in y\imply x=z\,\}$ and \œ$\{\ssp x\ssp,y\ssp\}=\{\sp x\sp\}
\cup\{\ssp y\sp\}\,$,\linebreak we have [ \œ$\all z\,z\in\{\sp x\sp\}\equivv x
=z$ ] if \œ$x\in y$ for some $y\ssp$, i.e.\ if $x$ is a set, otherwise having
$\{\sp x\sp\}=\Univ\,$, i.e.\ for any set $z$ it holds that $z \in
\{\sp x\sp\}\,$.

Deviating from {\cite[Definition 48, p.\ 259]{Ky}}\ssp, we define \œ$
(\sp x\ssp,y\sp)=\{\sp\{\ssp x\ssp,y\ssp\}\ssp,\{\ssp y\ssp\}\sp\}\,$, and
call $P$ an {\it ordered pair\ssp} if and only if there are sets $x\ssp,y$
with \œ$P=(\sp x\ssp,y)\,$. If $x$ or $y$ is not a set, it follows that $
(\sp x\ssp,y)=\Univ\,$.

Adapting {\cite[Definitions 51, 52, p.\ 259\,]{Ky}}\ssp, we define\par\centerline{$
\sigrd z=\bigcup\bigcup\sp z\setminus\bigcup\bigcap\sp z\sp\cup\sp\bigcap
\bigcup\sp z$ \ and \ $\taurd z=\bigcap\bigcap\sp z\,$.}

\noin Then \œ$x=\sigrd P$ and \œ$y=\taurd P$ hold for any ordered pair \œ$
P=(\sp x\ssp,y)\,$. Furthermore, we have $\,\sigrd\Univ=\taurd\Univ=\Univ\,$,
\,since for example\vskip.3mm\centerline{$
\sigrd\Univ = \bigcup\bigcup\ssp\Univ\setminus\bigcup\bigcap\ssp\Univ\sp\cup\sp
  \bigcap\bigcup\ssp\Univ =
\bigcup\ssp\Univ\setminus\bigcup\ssp\emptyset\sp\cup\sp\bigcap\ssp\Univ =
\Univ\setminus\emptyset\sp\cup\sp\emptyset = \Univ\,$.} \vskip.3mm

We define \œ$(\sp x\ssp,y\ssp,z\sp)=((\sp x\ssp,y\sp)\ssp,z\sp)$ and \œ$
(\sp x\,;y\ssp,z\sp)=(\sp x\ssp,(\sp y\ssp,z\sp))$ and \œ$
(\sp x\ssp,y\,;u\ssp,v\sp)=$ $(\sp x\ssp,y\ssp,(\sp u\ssp,v\sp))$ and \œ$
\ssigrd Z=\sigrd(\sp\sigrd\spp Z\sp)$ and \œ$\tsigrd Z =
\taurd(\sp\sigrd\spp Z\sp)\,$. In an obvious manner, one may continue to
get a succession of definitions for instance for $(\sp x\ssp,y\ssp,z\ssp,u\sp)$
and $(\sp x\,;\spp y\ssp,z\ssp,u\sp)$ and $\sssigrd\spp Z$ and $
\tau\ssigrd\spp Z$ and $\tau\yr 2\sigrd\spp Z\sp$, etc.    \end{conventions}

With the aid of the preceding conventions, we can give a precise meaning to
the concept of map. Any set $\sp\Cal T$ with \œ$\{\,\bigcup\ssp\Cal A:\Cal A
\inc\Cal T\sp\,\}\cup\{\,U\cap V:U\spp,V\sn\in\Cal T\sp\,\}\inc\Cal T\ssp$ be-
ing called a {\it topology\ssp}, a {\it topological space\ssp} is any ordered
pair $X$ such that $\taurd\sp X$ is a topology with \œ$\sigrd\sp X =
\bigcup\ssp\taurd\sp X\sp$. Letting \œ$f\ssp[\,A\,]=f\ssp\image\!A=\{\,y :
\exi{x\in A}\,(\sp x\ssp,y\sp)\in f\,\}$ and \œ$f\ssp\images\ssn\Cal A =
\{\,f\ssp\image\ssn A:A \in\Cal A\,\}$ and \œ$\Cal T\lei A=\{\,U\cap A:U \in
\Cal T\sp\,\}\,$, by a {\it topological map\ssp} (or topological morphism) we
mean any ordered pair $\tilde f$ such that topologies $\Cal T$ and $\sp
\Cal U\sp$ and a function $f$ exist with \œ$f\inc(\ssp\bigcup\ssp\Cal T\ssp)
\snn\times\snn(\ssp\bigcup\,\Cal U\ssp)$ and \œ$f\sp\inve\images\Cal U \inc
\Cal T\lei(\dom\sn f\sp)$ and\Biggerlineskip2 $\tilde f=(\ssp\Cal T\spp,\ssp
\Cal U\sp,f\ssp)\,$. A topological map $\tilde f$ is called {\it global\ssp}
if{}f $\,\bigcup\ssp\ssigrd\sp\tilde f \inc\dom\taurd\sp\tilde f\sp$.

We say that $f$ is {\it continuous\ssp} \œ$\Cal T\to\Cal U\ssp$ if and only if
\œ$\sp\tilde f=(\ssp\Cal T\spp,\ssp\Cal U\sp,f\ssp)$ is a topological map.
Equivalently, we may also say that $\tilde f$ is continuous, or a continuous
map, or that $f:\Cal T\to\Cal U$ is continuous, or a continuous map.

Note that we above did not want to mix the concept of topological space in
that of topological map because this only would have made things more
complicated. We only introduced the concept of topological space because in
some (quite rare) connections we may be able to simplify wordings by using
such a concept.

Instead of topological maps, we shall below mainly consider module or vector
maps defined as follows. A {\it ring\ssp} (structure) on $\Ke\ssp$ is any
ordered pair $\ssp\bold K$ such that functions \œ$\ssp\mfrak a\ssp,\mfrak c:
\Ke^{\sp\times 2.}=\Ke\snn\times\sn\Ke\to\Ke$ exist with \œ$\ssp\bold K =
(\sp\mfrak a\ssp,\mfrak c\sp)$ and satisfying the usual ring postulates when
we let \œ$s+t=\mfrak a\fvalue(\sp s\ssp,t\sp)$ and \œ$s\,t =
\mfrak c\fvalue(\sp s\ssp,t\sp)\,$. The ring $\ssp\bold K$ is called {\it
commutative\ssp} if{}f $\ssp\mfrak c$ is such, meaning that \œ$
\{\ssp(\sp r\sp,s\ssp,t\sp):(\sp s\ssp,r\sp,t\sp)\in\mfrak c\,\}\inc\mfrak c\ssp
$ holds, and {\it unital\ssp} if and only if it possesses a unity $\ssp
\roman l\,$, meaning that $\ssp\roman l\ssp$ is a $\ssp\mfrak c\,
$--\,identity, i.e.\ that \œ$\emptyset\not=\{\ssp(\sp r\sp,s\ssp,t\sp) :
\{\ssp r\sp,s\ssp\}=\{\ssp t\ssp,\sp\roman l\,\}\sp\}\,|\,
\Ke^{\sp\times 2.}\inc\mfrak c\,$, and we also have \œ$\ssp\roman l\not=\smb O$
when we let $\smb O=\bzero{\,\hbox{\font\=cmbx6\K}} =
\bigcap\ssp\{\,s:(\sp s\ssp,s\ssp,s\sp)\in\sigrd\spp\bold K\,\}\,$. Let $
\vecss X=\rng\sigrd\sp X\sp$.

A {\it structured $\,\bold K\,$--\,module\ssp} on $S$ is any ordered pair $E$
such that $\ssp\bold K$ is a commutative unital ring and functions \œ$a :
S^{\sp\times 2.}\to S$ and \œ$c:(\spp\vecss\sp\bold K\ssp)\sn\times\sn S\to S$
exist with \œ$\sigrd E=(\sp a\ssp,c\sp)$ and satisfying the usual module
postulates when we let \œ$x+y=$ $a\fvalue(\sp x\ssp,y\sp)$ and \œ$t\,x =
c\fvalue(\sp t\ssp,x\sp)\,$. Defining \œ$\vecs E=\rng\ssigrd E\,$, now a $\ssp
\bold K\,${\it--\,vector map\ssp} is any ordered pair $\tilde f$ such that
structured $\,\bold K\,$--\,modules $E\ssp,F$ and a function $f$ exist
satisfying \œ$f\inc(\sp\vecs E\ssp)\snn\times\snn(\sp\vecs F\ssp)$ and \œ$
\tilde f=(E\ssp,F\spp,f\ssp)\,$. A vector map $\sp\tilde f$ we agree to say to
be {\it global\ssp} if and only if $\,\vecs\sp\ssigrd\sp\tilde f \inc
\dom\taurd\sp\tilde f\,$ holds.

If in some connection we wish to be perfectly explicit about the structured
(vec- tor space or) module with respect to whose linear structure the
algebraic operations in a given linear combination are to be taken, we shall
use a notational device giving for example \œ$\sp
(\sp x+t\,u\sp)\sbi{\fiveroman{svs\,}E} = \ssigrd E\sp\fvalue(\sp x\ssp,
\tsigrd E\sp\fvalue(\sp t\ssp,u\sp))$ instead of the ambiguous \q{\œ$\sp
x+t\,u\sp$}. In a similar fashion, we have for example \œ$\ssp
[\,\{\sp t\sp\}\,A+B+C\ssp\,]\sbi{\fiveroman{svs\,}E} = $ \œ$
\ssigrd E\sp\,[\,\sp\ssigrd E\sp\,[\,\sp\tsigrd E\sp\,[\,\{\sp t\sp\}\timesn
A\sp\,]\timesn B\sp\,]\timesn C\ssp\,]\,$. If in such a connection we also
have a function $\sp f$ with \œ$\dom\sn f\inc\vecs E\,$, in view of
\cite[Theorem 69, p.\ 261]{Ky} then writing \œ$\sp y =
f\sp\fvalue((\sp x+t\,u\sp)\sbi{\fiveroman{svs\,}E})\not=\Univ\sp$ is
equivalent to writing \œ$x\ssp,u\in\vecs E$ and \œ$t\in\Ke$ and $x+t\,u\in
\dom\sn f$ and $\sp y=f\sp\fvalue(\sp x+t\,u\sp)$ in the looser convention.

In the preliminaries to \cite{Hic} we gave the formal definition of the two
topology Tihonov product $\sp\Cal T\ttimes\Cal U\ssp$. We omit the definition
of the corresponding product $\,\prod\sn\LHB{.2}{\ar{ti}}\sp\bCal T\,
                                                                     $ \ of an
arbitrary small family $\sp\bCal T$ \ of topologies. The vector space product
of two algebraic modules $X\sp,Y$ is $X\vstimes Y\spp$. The corresponding
product of a small family $\biit X\ssp$ is $\,\vsprod_VS_\Ke\biit X\sp$ when
each $X\in\rng\biit X\sp$ is a $\ssp\bold K\,$--\,module with $\Ke =
\vecss\sp\bold K\,$.\vskip.3mm

We further let \œ$X\expnota^\Omega]_{vs} = \prod\sn\RHB{.5}{_{_{\roman{VS\sp
    \,\sixroman{dom}^2\sp\sixmath\tauu_{\snn r\snn d}\sp\sixmath X\,}}}}
(\sp\Omega\snn\times\sn\{\sp X\ssp\}\spp)\,$. If $X$ is a module (structure)
and $S$ is a submodule (set) therein, then \œ$X_{\sp|\ssp S} =
(\sp\sigrd\spp X\,|\,S^{\sp\times2.}\sp,
\taurd\spp X\,|\,(\sp\Univ\sn\times\sn S\sp))$ is the corresponding submodule
structure. We let $E_{\sp/\sp S} =
          (\sp\sigrd E_{\ |\ssp S}\ssp,\taurd E\lei S\sp)\,$.   \vskip.5mm

For convenience, we modify our notation
   \q{$\sp\seq{\,t\sbi z:z\in S\,}\sp$}  from \cite{Hic} as follows.  \vskip.3mm

We let \œ$\seq{\,\mfrak T:\mfrak y\ar 1\sp,\ldots\,\mfrak y\ssp\ai l\ssn:
\mfrak x\,\mfrak E\,}=\{\,\mfrak z:\exi{\mfrak x\ssp,\mfrak x\ar 1\sp,\ldots\,
\mfrak x\ai k}\,\mfrak z=(\ssp\mfrak x\,,\mfrak T\sp)\text{ and }\mfrak x\,\mfrak E
\,\}$ when $\mfrak x\ssp,\mfrak x\ar 1\sp,$ \œ$\ldots\,\mfrak x\ai k\ssp,
\mfrak y\ar 1\sp,\ldots\,\mfrak y\ssp\ai l\ssp,\mfrak z$ are distinct variable
symbols and $\mfrak T$ is a term and $\mfrak E$ is an expression such that the
expression $\mfrak x\,\mfrak E$ is a formula having the variable symbol $\mfrak x
$ in the first place, and the common free variables in $\mfrak T$ and $
\mfrak x\,\mfrak E$ distinct from \œ$\mfrak x\ssp,\mfrak y\ar 1\sp,\ldots\,
\mfrak y\ssp\ai l$ are exactly \œ$\mfrak x\ar 1\sp,\ldots\,\mfrak x\ai k\ssp$,
and $\mfrak z$ does not occur free in $\mfrak T$ or $\mfrak x\,\mfrak E\,$. In the
case where \œ$\mfrak y\ar 1\sp,\ldots\,\mfrak y\ssp\ai l$ is an empty list, we
let $\seq{\,\mfrak T:\mfrak x\,\mfrak E\,}=\seq{\,\mfrak T:\ :\mfrak x\,\mfrak E\,}
\,$. \hfill Thus for example the class \par$\mhyppy{12}
f=\seq{\,\roman e^{\,s\ssp+\ssp t}:x=(\sp s\sp,t\sp)\in\Re\times\snn\Re\ssp\,}$\par$\mhyppy{15}
 =\{\,z:\exi{x\ssp,s\ssp,t}\,z=(\sp x\ssp,\roman e^{\,s\ssp+\ssp t}\spp)\text{
 and }x=(\sp s\ssp,t\sp)\in\Re\times\snn\Re\sp\,\}$ \hfill            is the

\noin function defined on \œ$\dom\sn f=\Re\times\snn\Re\,$, and whose value at
\œ$x=(\sp s\ssp,t\sp)\in\Re\times\snn\Re$ is \œ$f\sp\fvalue x=
\roman e^{\,s\ssp+\ssp t}$, and which one might express by \œ$
f:\Re\times\snn\Re\owns(\sp s\ssp,t\sp)\mapsto\roman e^{\,s\ssp+\ssp t}$. For
\œ$g=$ \œ$\seq{\,\roman e^{\,t}:t:x=(\sp s\ssp,t\sp)\in\Re\times\snn\Re\ssp\,}
=\{\,z:\exi x\,z=(\sp x\ssp,\roman e^{\,t}\spp)\text{ and }x=(\sp s\ssp,t\sp)
\in\Re\times\snn\Re\sp\,\}\,$, we have $g=\{\ssp(\sp s\ssp,t\ssp,
\roman e^{\,t}\spp)\ssp\}\ssp$ if $\,s\ssp,t\in\Re\,$, and $g=\emptyset$ otherwise.

When compared to our old convention, for example, if $f$ is a function, then \œ$
f=\seq{\,f\sp\fvalue x:x\in\Univ\sp\,}\subtext{old}
 =\seq{\,f\sp\fvalue x:x\in\dom\sn f\sp\,}\subtext{old}=\seq{\,f\sp\fvalue x:f
:x\in\dom\sn f\sp\,}\,$, whereas\Biggerlineskip1 $
\seq{\,f\sp\fvalue x:x\in\dom\sn f\sp\,}=\{\,z:\exi{x\ssp,f}\,z=(\sp x\ssp,
f\sp\fvalue x\sp)$ and $x\in\dom\sn f\sp\,\}=\Univ^{\sp\times 2.}\sp$.

Recall that a {\it function\ssp} (or family) is any \œ$f \inc
\Univ^{\sp\times 2.}$ such that no $x\ssp,y\ssp,z$ exist with \œ$(\sp x\ssp,
y\sp)\ssp,(\sp x\ssp,z\sp)\in f$ and \œ$y\not=z\ssp$. We generally put \œ$
R\ssp\inve=\{\ssp(\sp y\ssp,x\sp):(\sp x\ssp,y\sp)\in R\sp\,\}\,$ and \œ$\,
g\spp\circ\snn f=\{\ssp(\sp x\ssp,z\sp):\exi y\,(\sp x\ssp,y\sp)\in f\text{
and }(\sp y\ssp,z\sp)\in g\,\}\,$. \hfill We shall also utilize the\vskip.3mm
\noin definitions $\mhyppy{2.3}
f\yvee = \{\ssp(\sp x\ssp,u\sp):\emptyset\not=u=\{\ssp(\sp y\ssp,z\sp):
                  (\sp x\ssp,y\ssp,z\sp)\in f\,\}\sp\}\,$,\vskip.3mm

$\mhyppy{15.3} \yvee\snn f = \{\ssp(\sp y\ssp,u\sp):\emptyset\not=u=
\{\ssp(\sp x\ssp,z\sp):(\sp x\ssp,y\ssp,z\sp)\in f\,\}\sp\}\,$,\vskip.3mm

$\mhyppy{15}  f\aar 1\KN1\ywed = \{\ssp(\sp x\ssp,y\ssp,z\sp):\exi u\,
(\sp x\ssp,u\sp)\in f\aar 1\text{ and }(\sp y\ssp,z\sp)\in u\,\}\,$,\vskip.3mm

$\mhyppy{14}  \ywed f\aar 1 = \{\ssp(\sp x\ssp,y\ssp,z\sp):\exi u\,
(\sp y\ssp,u\sp)\in f\aar 1\text{ and }(\sp x\ssp,z\sp)\in u\,\}\,$,\vskip.3mm

$\mhyppy{14.5} \roman{pr}\ar 1 = \{\ssp(\sp x\ssp,y\ssp,x\sp):x\ssp,y\in
\Univ\sp\,\}\,$, \hfill $\roman{pr}\ar 2 =
\{\ssp(\sp x\ssp,y\ssp,y\sp):x\ssp,y\in\Univ\sp\,\}\,$.\KP{9.5}\vskip.3mm

\noin The large {\it evaluation family\ssp} is \ $\roman{ev}=\{\ssp(\sp x\ssp,
u\ssp,y\sp):u\text{ function and }(\sp x\ssp,y\sp)\in u\,\}\,$. \hfill The 
evaluation family at $\sp x\sp$ is $\,\roman{ev}\sbi{\eightmath x} =
\{\ssp(\sp u\ssp,y\sp):u\text{ function and }(\sp x\ssp,y\sp)\in u\,\}\,$. \vskip.3mm

Note that \œ$0\adot=\emptyset$ and \œ$\sp 1\adot=\{\emptyset\}$ and \œ$\sp
2\adot=\{\,\emptyset\,,1\adot\ssp\}$ and \œ$\N=\No\ssn\setminus\sn 1\adot\ssp
$, and that \œ$k+1\adot=k\ssp\yplus=k\sp\cup\{\sp k\sp\}$ for \œ$k\in\No\ssp$,
and also that \œ$\infty=\No$ and \œ$\infty\ssp\yplus=\No\cup\{\infty\}\,$.
If $\,\bosy x\in\Univ\,^k$ and $\,\bosy y\in\Univ\,^l$ and $\,k\ssp,l\in\No\ssp
$, \,then $\,\bosy x\concc\bosy y =
   \bosy x\sp\cup\sp(\sp\bosy y\circ\seq{\,k+j:j\in l\sp\,}\inve\spp)\,$.

For simplicity, we allow the notational inconsistency that for extended real
numbers $\ssp s\le t\sp$ we have \œ$\,[\,s\ssp,t\,]=\{\,r:s\le r\le t\,\}\,
$, whereas for functions $\sp f\sp,\sp g\sp$ we define \œ$[\,f\sp,\sp g\,] =
\{\ssp(\sp x\,;y\ssp,z\sp):(\sp x\ssp,y\sp)\in f\text{ and }(\sp x\ssp,z\sp)
\in g\,\}\,$. Proper notations for these would be for example $\,
[\,s\ssp,t\,]\subtext i\sp$ and $\,[\,f\sp,\sp g\,]\subtext f\,$, \,respectively.

We generally (but not exclusively, for example \œ$\,x+y\snn\cdot\snn z \not=
(\sp x+y\sp)\snn\cdot\snn z\,$ usually) apply the rule for reduction of
parentheses given by the schema \œ$\,\mfrak T\ai k\ssp\roman b\ai k\ssp\ldots\,
\mfrak T\ar 1\ssp\roman b\ar 1\,\mfrak T\ar 0$ \œ$ =
(\sp\mfrak T\ai k\ssp\roman b\ai k\ssp\ldots\,\mfrak T\ar 1)\,\roman b\ar 1\,
\mfrak T\ar 0$ when the $\mfrak T\sbi\iota$ are terms and $\ssp
\roman b\sbi\iota$ are binary symbols. Also, we usually understand monadic
symbols to act prior to the binary ones. For example, above we have $\hfill
\sigrd z=((\ssp\bigcup\,(\ssp\bigcup\sp z\sp))\setminus(\ssp\bigcup\,(\ssp
\bigcap\sp z\sp)))\cup\sp(\ssp\bigcap\,(\ssp\bigcup\sp z\sp))\,$. \hfill
Further, usually\linebreak we let \œ$
\roman m\ai k\ldots\,\roman m\ar 0\,\mfrak T\,
\roman n\ar 0\ldots\,\roman n\ai{\,l} =
((\sp\roman m\ai k\ldots\ssp(\sp\roman m\ar 0\,\mfrak T\ssp))\,\roman n\ar 0)
\ldots\,\roman n\ai{\,l}$ when the $\roman m\sbi\iota$ and $\roman n\sbi\iota$
are monadic. For example \œ$\,\bigcup\ssp\Phi\sp\inve =
(\ssp\bigcup\ssp\Phi\sp)\inve = \{\ssp(\sp y\ssp,x\sp):\exi f\,
(\sp x\ssp,y\sp)\in f\in\Phi\,\}\,$, and \linebreak usually $\,
            \bigcup\ssp\Phi\inve\not=\bigcup\,(\sp\Phi\inve\spp)\,$.

\begin{remark}

If we wished to give systematically our rules for the reduction of parentheses
so that the use of such vague words as \q{usually} could be avoided, we should
give an exhaustive list of our function symbols, the symbols there being
grouped according to their intended \q{level}, like in
\cite[Theory of notation, 0.30\,--\,0.48, pp.\ 15\,--\,21]{Mo}\ssp, cf.\ the
\q{type} and \q{power} there. The preceding rules then would be applied
without exception when symbols on the same level are considered, whereas
symbols on a higher level would act prior to those on a lower one.

For example \,lev\œ$\,\Symboo{\sp{\cdot}\spp}Ï>\roman{lev\,}\Symboo{{+}\snn}Ï$
and \,lev\,\œ$\Symboo^{\sp\fiveroman{fct\ssp exp}}\snnÏ >
\roman{lev\,}\Symboo{\spp\image\sn}Ï$ so that we have \œ$\,r+s\cdotn t$ \œ$ =
r+(\sp s\cdotn t\sp)$ and \œ$f\ssp\image B^{\,k}=f\ssp\image(B^{\,k}\spp)=
f\sp\,[\,B^{\,k}\ssp]\,$, here $\Symboo^{\sp\fiveroman{fct\ssp exp}}\snnÏ$
being the invisible bin- ary function symbol the result of whose application
to $\Symboo BÏ\,\Symboo kÏ$ we write \q{$\sp B^{\,k}\sp$}.

Since we accept this kind of \q{nonlinear} expressions which are discarded in
\cite{Mo}\ssp, our theory of notation, if properly presented, is much more
complicated than the one in \cite[pp.\ 15\,--\,26]{Mo}\ssp. Note further that
although \,lev\,$\Symboo^{\sp\fiveroman{fct\ssp exp}}\snnÏ >
    \roman{lev\,}\Symboo{{+}\snn}Ï\sp$, we have $B\,^{k\ssp+\ssp l} =
B\,^{(k\ssp+\ssp l\spp)}\sp$, and generally $B\,^{k\ssp+\ssp l} \not=
B^{\,k}+l$ \,unless $\sp k=l=\emptyset\,$.                      \end{remark}

\begin{remark}\label{Weihe pairs}

Suppose, cf.\ \cite[p.\ 59]{Mo}\ssp, that one introduces the ordered {\it
Weihe pairs\ssp} \œ$(\sp x\ssp,y\sp)\subtext{we} = x\weco y =
(\sp\Pows x\sp)\timesn(\sp\Pows y\sp)=\{\ssp(\sp u\ssp,v\sp):u\inc x\text{ and }
v\inc y\,\}\,$ having the prop-\linebreak erty that the implication \,\œ$
[\ x\weco y=u\weco v\imply x=u\text{ and }y=v\ ]\,$ {\it holds for any classes\ssp}
$x\ssp,y\ssp,u\ssp,v\ssp$, and not only for such sets as in the case of our
(reversed) Wiener pairs, and further the ordered Weihe tuples by the recursion
schema \œ$\,(\sp x\ai k\ssp,\ldots\,x\ar 1\sp,x\ar 0)\subtext{we} = $ \œ$
x\ai k\weco\ldots\,x\ar 1\weco x\ar 0 =
(\sp x\ai k\weco\ldots\,x\ar 1)\weco x\ar 0\,$. Then one could manage with
very few binary symbols and with no higher order ones except the Weihe tuple
symbols themselves. Indeed, for example given a $\sp\ai k\ssp$--\,place
function symbol $\mfrak g$ and introducing the mon- adic function symbol $
\mfrak f$ by the definition\vskip.3mm\centerline{$
\mfrak f\,x = \{\,u:\all{x\ar 1\sp,\ldots\,x\ai k}\,x =
x\ar 1\weco\ldots\,x\ai k\imply u\in\mfrak g\,x\ar 1\,\ldots\,x\ai k\ssp\}\,$,}\vskip.3mm

\noin one could totally eliminate $\mfrak g$ from one's set theoretic edifice.

  \end{remark}


  \binsubsubhead A{A reformulation of the general BGN\,--\,setting}\label{subsec A}

In \cite[Definition 1.1, pp.\ 219\,--\,220]{BGN}\ssp, postulates are given for
a (generally large)\linebreak family \œ$[\,\sp\Cal C\yr 0\spp,\bosy\tau\,]:
(E\ssp,F\sp)\mapsto(\sp\Cal C\yr 0\fvalue(E\ssp,F\sp)\ssp,\bosy\tau\fvalue
(E\ssp,F\sp))$ where further \œ$\sp\Cal C\yr 0\fvalue(E\ssp,F\sp)$ is a\linebreak
(small) family \œ$\taurd E\owns U\mapsto\Cal C\yr 0\fvalue(E\ssp,F\sp)\fvalue
U\spp$, the last one being a set of functions\linebreak \œ$U\to\vecs F$ which
are continuous \œ$\taurd E\to\taurd F\sp$, and where $\bosy\tau\fvalue(E\ssp,
F\sp)$ is a topology for the set $(\sp\vecs E\sp)\snn\times\snn(\sp\vecs F\sp)\,
$. We reformulate these postulates in the next

\begin{definitions}\label{BGN axioms}

By a {\it topologized $\,\bold K\,$--\,module\ssp} we mean any structured $\ssp
\bold K\,$--\,module $E$ such that $(\sp\vecs E\ssp,\taurd E\ssp)$ is a
topological space. If $\ssp\bold K\sp$ is a division ring, or a field, we
speak of {\it topologized vector spaces\ssp}. Let now $\ssp\bold K\sp$ be a
commutative ring with unity $\ssp\roman l\,$, and let $\O$ be a class of
topologized $\bold K\,$--\,modules with \œ$\biit K = (\ssp\bold K\,,
\tau_{_{I\!\!K}}\snn)\in\O\ssp$. Letting \œ$F^{\sp/\sp E}=\{\,f:f\text{ a
function and }f\inc(\sp\vecs E\sp)\snn\times\snn(\sp\vecs F\sp)\ssp\}\,$, we
consider a class $\Cal C\ar 0$ having as its members some maps \œ$\tilde f =
(E\ssp,F\sp,f\sp)\in\O^{\sp\times 2.}\ssn\times\snn\Univ$ where \œ$f \in
F^{\sp/\sp E}\sp$. We require each \œ$\tilde f \in\Cal C\ar 0$ to be
continuous with open domain, which is expressed shortly by the inclusion $\,
\taurd\sp\tilde f\ssp\inve\sp\images\tau\yr 2\sigrd\sp\tilde f \inc
\tau\ssigrd\sp\tilde f\ssp$. \,Letting\vskip.4mm\centerline{$
\Cal C\ar 0\vProd = \{\ssp(E\ssp,F\sp,G\sp):E\ssp,F\in\domm\Cal C\ar 0$ and $\sp
\Cal C\ar 0\text{\ssp-\,prod}\ssp\subtext{mcl}\sp(E\ssp,F\sp,G\sp)\ssp\}\,$,}\vskip.4mm

\noin where $\sp\Cal C\ar 0\text{\ssp-\,prod}\ssp\subtext{mcl}\sp(E\ssp,F\sp,
G\sp)$ denotes the condition that\vskip.4mm\centerline{$
\sigrd G=(\sp\sigrd E\ssp)\vstimes(\sp\sigrd F\ssp)$ and $
(\spp G\sp,E\ssp,\roman{pr}\ar 1\ssp|\,\vecs G\sp)\ssp,
(\spp G\sp,F\sp,\roman{pr}\ar 2\ssp|\,\vecs G\sp)\in\Cal C\ar 0$}\nopagebreak\vskip.2mm\nopagebreak\centerline{%
and $\sp\all{f\sp,g\ssp,H}\,
(\spp H\sp,E\ssp,f\sp)\ssp,(\spp H\sp,F\sp,g\sp)\in\Cal C\ar 0\imply
(\spp H\sp,G\sp,[\,\sp f\sp,\spp g\sp\,]\ssp)\in\Cal C\ar 0\,$,}\vskip.4mm

\noin we now say that $\Cal C\ar 0$ is a {\it productive class\ssp} on $\O$
over $\biit K$ if{}f in addition to the above assumptions, we also have \œ$\ssp
\O^{\sp\times 2.}\inc\dom\sn(\sp\Cal C\ar 0\vProd\spp)\,$. A
productive class $\sp\Cal C\ar 0$ on $\O$ over $\biit K$ we say to be a
\erm{BGN}{\it\,--\,class\ssp} if{}f with the inversion \œ$\Iota =
\{\ssp(\sp s\ssp,t\sp):(\sp s\ssp,t\ssp,\sp\roman l\ssp) \in
\taurd\sp\bold K\sp\,\}\,$, for all $E\ssp,F\sp,G\in\O$ and for all $
f\sp,g\ssp,y\ssp,P\spp,U\spp$, we  also have the following

\begin{enumerate}

\item\label{locality}
      \ \ $P=(E\ssp,F\sp)$ and $f$ is a function and\nopagebreak\par\nopagebreak\hyppy{-.9mm}%
  [ $\spp\all{\smb Z}\,\exi h\,\smb Z\in f\imply\smb Z\in h\inc f$ and $(\spp
  P\spp,h\sp)\in\Cal C\ar 0\ssp\,]\imply(\spp P\spp,f\sp)\in\Cal C\ar 0\,$,\vskip.3mm

\item\label{compos rule}
      \ \ $(E\ssp,F\sp,f\sp)\ssp,(\spp F\sp,G\sp,g\sp)\in\Cal C\ar 0\imply
  (E\ssp,G\sp,g\circ f\sp)\in\Cal C\ar 0\,$,\vskip.3mm

\item\label{id and open const in C_0}
      \ \ $U\in\taurd E$ and $y\in\vecs F\imply(E\ssp,E\ssp,\idv E\sp)\ssp,
  (E\ssp,F\sp,\spp U\snn\times\sn\{\sp y\sp\}\sp)\in\Cal C\ar 0\,$,\vskip.3mm

\item\label{inversion rule}
     \ \ $(\biit K\ssp,\biit K\ssp,\sp\Iota\ssp)\in\Cal C\ar 0\,$,\vskip.3mm

\item\label{determination axiom}
      \ \ $f\sp,g\in\Cal C\ar 0\KN1\image\snn\{\ssp(\biit K\ssp,F\ssp)\ssp\}$
      and $\dom\sn f=\dom g$ and $\sp f\circ\Iota\circ\Iota\inc g\sp
      \imply\sp f=g\,$,\vskip.3mm

\item \ \ $
   (\sp\Cal C\ar 0\vProd\KN1\fvalue(E\ssp,E\sp)\ssp,E\ssp,\ssigrd E\ssp)\,,\sp
   (\sp\Cal C\ar 0\vProd\KN1\fvalue(\biit K\ssp,E\sp)\ssp,E\ssp,\tsigrd E\ssp)
   \in\Cal C\ar 0\,$.                      \end{enumerate} \end{definitions}

\begin{lemma}\label{for product family}
                                        
If $\ssp\Cal C\ar 0$ is a productive class on $\sp\O$ over $\sp\biit K\,,$
it then holds that \œ$\dom\Cal C\ar 0\inc$ $\dom\sn(\sp\Cal C\ar 0\vProd\spp)=
\O^{\sp\times2.}\sp,$ and that $\,\Cal C\ar 0\vProd$ is a function $\sp
\O^{\sp\times2.}\to\O\ssp$.                                      \end{lemma}

\begin{proof} Let $\sp\Cal C\ar 0$ be a productive class on $\sp\O$ over $\sp
\biit K\,$. To see that $\sp\Cal C\ar 0\vProd$ is a function, letting \œ$
(E\ssp,F\sp,G\sp)\ssp,(E\ssp,F\sp,H\sp)\in\Cal C\ar 0\vProd\,$, then \œ$
\sigrd G=(\sp\sigrd E\ssp)\vstimes(\sp\sigrd F\ssp)=$ $\sigrd H\sp$. If we
take \œ$f=\roman{pr}\ar 1\ssp|\,\vecs G$ and \œ$g=\roman{pr}\ar 2\ssp|\,\vecs
G\sp$, then \œ$\idv G=[\,\sp f\sp,\spp g\sp\,]=\idv H$ with \œ$
(\spp G\sp,H\sp,\idv G\sp)\ssp,(\spp H\sp,G\sp,\idv G\sp)\in\Cal C\ar 0\,$,
whence by the assumption that every \œ$\tilde f\in\Cal C\ar 0$ is continuous
we get \œ$\taurd H\inc\taurd G\inc\taurd H\sp$, hence \œ$G=H\sp$. To get \œ$
\dom\sn(\sp\Cal C\ar 0\vProd\spp)=$\Biggerlineskip1 $\O^{\sp\times2.}\sp$,
note that \œ$\sp\O^{\sp\times2.}\inc\dom\sn(\sp\Cal C\ar 0\vProd\sp)\inc
(\sp\domm\Cal C\ar 0)^{\sp\times2.}\inc\O^{\sp\times2.}\sp$. We trivially have
$\rng\snn(\sp\Cal C\ar 0\vProd\spp)\inc\domm\Cal C\ar 0\inc\O\sp$ and $
\dom\Cal C\ar 0\inc\O^{\sp\times2.}\sp$.                         \end{proof}

\begin{proposition}\label{product axioms}
                                        
Let $\ssp\Cal C\ar 0$ be a \ssp\erm{BGN}\,--\,class on $\sp\O$ over $\sp
\biit K\,$. Then \œ$\ssp\dom\Cal C\ar 0=\O^{\sp\times2.}\sp,$ and for all $\ssp
E\ssp,F\sp,G\in\O$ and for all $\ssp f\sp,\spp g$ with $\,\varPi =
\Cal C\ar 0\vProd\KN1\fvalue(\spp F\sp,G\sp)$ it holds that\vskip.6mm{\rm
(a)} \ \ $(\spp\varPi\sp,F\sp,\roman{pr}\ar 1\ssp|\,\vecs\varPi\sp)\,,\sp
    (\spp\varPi\sp,G\sp,\roman{pr}\ar 2\ssp|\,\vecs\varPi\sp)\in\Cal C\ar 0\,,$\vskip.4mm{\rm
(b)} \ \ $(E\ssp,F\sp,f\sp)\ssp,(E\ssp,G\sp,g\sp)\in\Cal C\ar 0\imply
    (E\ssp,\varPi\sp,\spp[\,\sp f\sp,\spp g\,\sp]\ssp)\in\Cal C\ar 0\,$.

  \end{proposition}

\begin{proof} Lemma \ref{for product family} gives \œ$\dom\Cal C\ar 0\inc
\O^{\sp\times2.}\sp$, and \œ$\sp\O^{\sp\times2.}\inc\dom\Cal C\ar 0$ follows
from (\ref{id and open const in C_0}) of Definitions \ref{BGN axioms} above.
If \œ$\sp F\sp,G\in\O$ and \œ$\sp\varPi =
\Cal C\ar 0\vProd\KN1\fvalue(\spp F\sp,G\sp)\,$, since by
Definitions \ref{BGN axioms} we have \œ$\sp(\spp F\sp,G\sp)\in\O^{\sp\times2.}
\inc\dom\sn(\sp\Cal C\ar 0\vProd\spp)\,$, and since by
Lemma \ref{for product family} we know that $\sp\Cal C\ar 0\vProd$ is a
function \œ$\sp\O^{\sp\times2.}\to\O\ssp$, it follows that \œ$\sp(\spp F\sp,
G\sp,\varPi\sp)\in\Cal C\ar 0\vProd\,$, and consequently that $\sp
\Cal C\ar 0\text{\ssp-\,prod}\ssp\subtext{mcl}\sp(\spp F\sp,G\sp,\varPi\sp)$
holds. This directly gives the asserted (a) and (b) above.       \end{proof}

Note that for the validity of (a) and (b) in Proposition \ref{product axioms}
above it in fact suffices that $\sp\Cal C\ar 0$ only is a productive class on
$\sp\O$ over $\sp\biit K\,$. However, without the other BGN postulates,
especially (\ref{compos rule}) of Definitions \ref{BGN axioms} above, these
are of little use. For this reason, we gave preference to the preceding
formulation.

Note also that by \œ$\sp\Iota\circ\Iota\inc\roman{id}\,$ we could have given
(\ref{determination axiom}) of Definitions \ref{BGN axioms} the following
longer equivalent formulation: we have \œ$f=g$ whenever \œ$\sp
\{\ssp(\biit K\ssp,F\ssp)\ssp\}\timesn\{\ssp f\sp,\spp g\ssp\}\inc\Cal C\ar 0$
and \œ$\dom\sn f=\dom g$ and \œ$\sp f\sp\fvalue t = g\fvalue t$ for all $\sp
\tsigrd\sp\biit K\,$--\,invertible $t\ssp$. This corresponds to the {\it
determination axiom\ssp} \cite[\erm{III}\sp, p.\ 220]{BGN}\ssp. Our
(\ref{inversion rule}) corresponds to \cite[\erm I\sp.4, p.\ 220]{BGN}\ssp,
and (\ref{locality}) to \cite[\erm I\sp.5, p.\ 220]{BGN}\ssp. For a further
comparison, see also Remark \ref{comparison of our's to orig. BGN} below.

\begin{definitions}

For any classes $\sp\Cal C\ar 0\,,\sp\biit K\,,\sp\tilde f$ generally letting\vskip.5mm\noin
(1) \ \ $\Cal C\ar 0\vDifq^1_K = \Cal C\ar 0\KN{.5}^{\times2.}\cap
\{\ssp(\sp\tilde f\sp,\tilde h\sp) :
\all{E\ssp,F\sp,f}\,\exi{E\ar 2\ssp,G\sp,h}\,\all{h\yr 0}\,\tilde f=$ \par $\null\hfill
(E\ssp,F\sp,f\sp)\imply\tilde h = (\spp G\sp,F\sp,\spp h\ssp)$ and $
(E\ssp,E\ssp,E\ar 2)\ssp,(E\ar 2\ssp,\biit K\ssp,G\sp)\in\Cal C\ar 0\vProd$\KP{11}\par\hyppy{10mm}
and \sp[ $h\yr 0 = \{\ssp(\sp x\ssp,u\ssp,t\ssp,y\sp) :
f\sp\fvalue((\sp x+t\,u\sp)\sbi{\fiveroman{svs\,}E}) =
(\sp f\sp\fvalue x+t\,y\sp)\sbi{\fiveroman{svs\,}F}\not=\Univ\sp\,\}$ \par$\null\hfill
\imply\dom h\yr 0\inc\dom h$ and $h\inc h\yr 0\ssp\,]\ssp\,\}\,$,\KP8\vskip.5mm\noin
(2) \ \ $\Cal C\ar 0\vbarDelta_K\tilde f = \bigcap\ssp\{\,\tilde h :
\all{\tilde h\ar 1}\,\tilde h = \tilde h\ar 1\sp\equivv\sp(\sp\tilde f\sp,
\tilde h\ar 1)\in\Cal C\ar 0\vDifq^1_K\ssp\}\,$,\vskip.5mm\noin  
(3) \ \ $\CalDBGN^k(\sp\Cal C\ar 0\ssp,\biit K\ssp) = \{\,\tilde f :
   k\in\infty\ssp\yplus$ and $\exi{\biit f}\,(\sp\emptyset\,,\tilde f\,)
     \in\biit f\in\Cal C\ar 0{^{k\ssp+\sp 1.}}$ and \par$\null\hfill
\all{i\in k}\,(\sp\biit f\ssp\fvalue i\ssp,\biit f\ssp\fvalue i\ssp\yplus)
              \in \Cal C\ar 0\vDifq^1_K\,\}\,$,\KP8\vskip.7mm

\noin we then say that $\tilde f$ is $\ssp k\,^{\fiveroman{th}}$ order $\sp
\Cal C\ar 0\ssp$--$\subtext{\,BGN\,}${\it differentiable\ssp} over $
\bosy K\ssp$ if and only if we have \œ$\sp
\tilde f\in\CalDBGN^k(\sp\Cal C\ar 0\ssp,\bosy K\ssp)\,$. We may say {\it
first\ssp} instead of $1\adot\sn^{\fiveroman{th}}\,$, and {\it second\,}
instead of $2\adot\sn^{\fiveroman{th}}\,$, \biggerlineskip5 etc. For $\sp
\tilde f\in\Univ\,$, we may call any \œ$\sp\tilde f\yr 1 \in
\Cal C\ar 0\vDifq^1_K\KN{.5}\image\snn\{\sp\tilde f\ssp\}$ a first order
\erm{BGN}{\it\,--\,difference \biggerlineskip6 quotient\ssp} map for $\sp
\tilde f$ in $\sp\Cal C\ar 0$ over $\biit K\,$. In case \œ$\sp
\Cal C\ar 0\vbarDelta_K\tilde f\not=\Univ\,$, we may call $\sp
\Cal C\ar 0\vbarDelta_K\tilde f$ \biggerlineskip5 {\it the\ssp} first order
\erm{BGN}\,--\,difference quotient map for $\sp\tilde f$ in $\sp\Cal C\ar 0$
over $\biit K\,$.                                          \end{definitions}

\begin{lemma}

For any $\sp\Cal C\ar 0\ssp,\biit K\,,\tilde f$ either \œ$\sp
\Cal C\ar 0\vbarDelta_K\tilde f\in\Cal C\ar 0$ or \œ$\sp
\Cal C\ar 0\vbarDelta_K\tilde f=\Univ\,$. Further- more $\,
\Cal C\ar 0\vbarDelta_K\tilde f\in\Cal C\ar 0$ holds if{}f a unique $\sp
\tilde f\yr 1$ exists \,with $\,(\sp\tilde f\sp,\tilde f\yr 1) \in
\Cal C\ar 0\vDifq^1_K\,$.                                        \end{lemma}

\begin{proof} Letting \œ$\sp\Cal H=\{\,\tilde h:\all{\tilde h\ar 1}\,\tilde h
= \tilde h\ar 1\sp\equivv\sp(\sp\tilde f\sp,\tilde h\ar 1) \in
\Cal C\ar 0\vDifq^1_K\ssp\}\,$, one quickly\linebreak verifies that either
there is a set $\tilde h$ with \œ$\Cal H=\big\{\ssp\tilde h\ssp\big\}$ or \œ$
\Cal H=\emptyset\,$. In the former case we have \œ$
\Cal C\ar 0\vbarDelta_K\tilde f=\bigcap\ssp\Cal H=\tilde h \in
\rng\snn(\sp\Cal C\ar 0\vDifq^1_K\spp) \inc
\rng\snn(\sp\Cal C\ar 0\KN{.5}^{\times2.}) \inc \Cal C\ar 0\,$, and in the\linebreak
latter $\,\Cal C\ar 0\vbarDelta_K\tilde f=\bigcap\ssp\emptyset=\Univ\,$. The
assertions now easily follow.                                    \end{proof}

If $\sp\tilde f$ is first order $\sp\Cal C\ar 0\ssp$--$\subtext{\,BGN\,}
$differentiable over $\biit K\,$, then $\sp\tilde f$ has at least one first
order \erm{BGN}\,--\,difference quotient map in $\sp\Cal C\ar 0$ over $\biit K\,
$. In the case where $\sp\Cal C\ar 0$ is a \erm{BGN}\,--\,class over $\biit K\,
$, there is only one such, namely $\sp\Cal C\ar 0\vbarDelta_K\tilde f\sp$, as
follows from the next

\begin{proposition}

If $\,\Cal C\ar 0$ is a \,\erm{BGN}\,--\,class over $\biit K\sp$ and $\ssp
\tilde f$ is first order $\ssp\Cal C\ar 0\ssp$--$\subtext{\,BGN\,}$differ-
entiable over $\biit K\,,$ then $\ssp\Cal C\ar 0\vbarDelta_K\tilde f$ is the
unique first order \erm{BGN}\,--\,difference quotient map for $\sp\tilde f$ in
$\sp\Cal C\ar 0$ over $\biit K\,$.                         \end{proposition}

\begin{proof}

Indeed, supposing that $\sp\tilde f^{\,\iota\,}$ for $^{\,\iota\ssp=\sp
\sixroman{1\sp,\ssp 2}}$ are first order difference quotient maps for $\sp
\tilde f$ in $\sp\Cal C\ar 0$ over $\biit K\,$, we first get \œ$\sp
\sigrd\spp\tilde f\yr 1=\sigrd\spp\tilde f\yr 2$ since $\Cal C\ar 0\vProd$ is
a function.\linebreak Putting \œ$\sp h\sbi\iota=\taurd\spp\tilde f^{\sp\iota}
$, from \œ$\dom h\yr 0\inc\dom h\sbi\iota$ and \œ$h\sbi\iota\inc h\yr 0$, we
obtain \œ$\dom h\ar 1=$\linebreak \œ$\dom h\yr 0=\dom h\ar 2\,$, whence it
suffices to prove that for any \œ$\sp\smb W=(\sp x\ssp,u\ssp,t\sp) \in
\dom h\yr 0$ we have the equality $\sp h\ar 1\KN1\fvalue\smb W =
h\ar 2\KN1\fvalue\smb W\sp$. Let $\sp f=\taurd\spp\tilde f\sp$.

To obtain \œ$\sp h\ar 1\KN1\fvalue\smb W=h\ar 2\KN1\fvalue\smb W\sp$, putting
\œ$\sp \gamma = \seq{\,(\sp x\ssp,u\ssp,s\sp):s\in\vecs\bosy K\ssp\,}\,$, we
consider \œ$\sp\gamma\sbi\iota=$\linebreak $h\sbi\iota\snn\circ\gamma\,$. Then
\œ$\dom\gamma\ar 1=\dom\gamma\ar 2\,$, and using item (b) of
Proposition \ref{product axioms} and (\ref{compos rule}) and
(\ref{id and open const in C_0}) of Definitions \ref{BGN axioms} one deduces
that we also have \œ$\sp(\biit K\ssp,F\sp,\gamma\sbi\iota)\in\Cal C\ar 0$ for
$\,\sbi{\iota\ssp=\sp\sixroman{1\sp,\ssp 2}}\,$. For every \œ$\sp
s\in\dom\gamma\sbi\iota$ we have \œ$\sp\gamma\sbi\iota\KN1\fvalue s =
h\sbi\iota\KN1\fvalue(\sp x\ssp,u\ssp,s\sp)\,$, hence \œ$\sp
(\sp x\ssp,u\ssp,s\ssp,\gamma\sbi\iota\KN1\fvalue s\sp) \in
h\sbi\iota\inc h\yr 0\,$, and consequently \œ$\sp f\sp\fvalue(\sp x+s\,u\sp) =
f\sp\fvalue x+s\,(\sp\gamma\sbi\iota\KN1\fvalue s\sp)\,$, whence for
invertible $s$ we get \œ$\gamma\ar 1\KN1\fvalue s =$ \œ$
\Iota\fvalue s\,\sp(\sp f\sp\fvalue(\sp x+s\,u\sp)-f\sp\fvalue x\sp) =
\gamma\ar 2\KN1\fvalue s\,$. By (\ref{determination axiom}) of
Definitions \ref{BGN axioms} we get \œ$\sp\gamma\ar 1=\gamma\ar 2\,$, whence
finally $\,h\ar 1\KN1\fvalue\smb W=h\ar 1\KN1\fvalue(\sp\gamma\fvalue t\sp) =
h\ar 1\snn\circ\gamma\fvalue t=\gamma\ar 1\KN1\fvalue t =
\gamma\ar 2\KN1\fvalue t=h\ar 2\snn\circ\gamma\fvalue t =
h\ar 2\KN1\fvalue(\sp\gamma\fvalue t\sp)=h\ar 2\KN1\fvalue\smb W\sp$.
  \end{proof}

Our definition above of $\ssp k\ssp\yr{th}$ order $\sp\Cal C\ar 0\ssp$--$
\subtext{\,BGN\,}$differentiability over $\biit K\ssp$ precisely captures the
content of \cite[Definitions 2.1, 4.1, Remark 4.3, pp.\ 222, 228]{BGN} because
of the recursion rule given by the following

\begin{proposition}\label{BGN recurs}

Let $\ssp\Cal C\ar 0$ be a \erm{BGN}\,--\,class over $\sp\biit K\,$. For all $\ssp
\tilde f\ssp,\sp k\ssp$ then\vskip.5mm\centerline{$
\tilde f\in\CalDBGN^{k\ssp+\sp 1.}(\sp\Cal C\ar 0\ssp,\biit K\ssp)\sp\equivv\sp
\tilde f\in\CalDBGN^{1.}(\sp\Cal C\ar 0\ssp,\biit K\ssp)$ and $\,
\Cal C\ar 0\vbarDelta_K\tilde f \in
\CalDBGN^k(\sp\Cal C\ar 0\ssp,\biit K\ssp)\,$.}            \end{proposition}

\begin{proof} Having the setting fixed, to simplify the notations, we agree to
let \œ$\Cal D\sbi k=$ $\CalDBGN^k(\sp\Cal C\ar 0\ssp,\biit K\ssp)$ and $\ssp
\Delta\,\tilde f=\Cal C\ar 0\vbarDelta_K\tilde f$ and\vskip.4mm\centerline{$
\roman S\ssp(\sp\tilde f\sp,k\sp)=\{\,\biit f:(\sp\emptyset\,,\tilde f\,)\in
\biit f\in\Cal C\ar 0{^{k\ssp+\sp 1.}}$ and $\sp
\all{i\in k}\,\biit f\ssp\fvalue i\ssp\yplus =
              \Delta\ssp(\sp\biit f\ssp\fvalue i\sp)\ssp\}\,$.}\vskip.4mm

\noin Now letting \œ$\ssp\tilde f\in\Cal D\sbi{k\ssp+\sp 1.}\,$, there is \œ$
\biit f\in\roman S\ssp(\sp\tilde f\sp,k+1\adot)\,$. With \œ$\Iota = \seq{\,\sp
i\ssp\yplus\sn:i\in k+1\adot\ssp}\,$, then \œ$\ssp\biit f\sp\,|\,2\adot \in
\roman S\ssp(\sp\tilde f\sp,1\adot)$ and \œ$\ssp\biit f\circ\Iota \in
\roman S\ssp(\sp\Delta\,\tilde f\sp,k\sp)\,$, whence we get \œ$\tilde f \in
\Cal D\sbi{1.}$ and \œ$\ssp\Delta\,\tilde f\in\Cal D\sbi k\,$.\linebreak
Conversely, if \œ$\tilde f\in\Cal D\sbi{1.}$ and \œ$\ssp\Delta\,\tilde f \in
\Cal D\sbi k\,$, there are \œ$\ssp\biit f\ar 0\in\roman S\ssp(\sp\tilde f\sp,
1\adot)$ and \œ$\ssp\biit f\ar 1\in\roman S\ssp(\sp\Delta\,\tilde f\sp,k\sp)\,
$. Then $\,\biit f\ar 0\cup(\ssp\biit f\ar 1\snn\circ\Iota\inve\spp) \in
\roman S\ssp(\sp\tilde f\sp,k+1\adot)\,$, whence we see that $\ssp\tilde f \in
\Cal D\sbi{k\ssp+\sp 1.}\,$.                                     \end{proof}

Note that even if the productive class $\sp\Cal C\ar 0$ over $\biit K\sp$ is
not a \erm{BGN}\,--\,one, we still have  a recursion that \œ$\sp\tilde f\in
\CalDBGN^{k\ssp+\sp 1.}(\sp\Cal C\ar 0\ssp,\biit K\ssp)$ if and only if some
first order \erm{BGN}\,--\,difference quotient map $\sp\tilde f\yr 1$ for $\sp
\tilde f$ in $\sp\Cal C\ar 0$ over $\biit K$ satisfies $\sp\tilde f\yr 1 \in
\CalDBGN^k(\sp\Cal C\ar 0\ssp,\biit K\ssp)\,$, put concisely\vskip.3mm\centerline{$
\tilde f\in\CalDBGN^{k\ssp+\sp 1.}(\sp\Cal C\ar 0\ssp,\biit K\ssp)\,\equivv\,
\tilde f\in\Univ\,$ and $\,\Cal C\ar 0\vDifq^1_K\KN{.5}\image\snn\{\sp
\tilde f\ssp\}\cap\CalDBGN^k(\sp\Cal C\ar 0\ssp,\biit K\ssp)\not=\emptyset\,$.}\vskip.3mm

\begin{remark}\label{comparison of our's to orig. BGN}

The transition from our reinterpretation of the BGN setting to the original
one and vice versa is accomplished as follows. Given our $\sp\Cal C\ar 0\,$,
upon putting\vskip.6mm

$\mhyppy{6.2} \bosy\tau =
\seq{\,\taurd(\sp\Cal C\ar 0\vProd\KN1\fvalue P\ssp):P \in
\dom\Cal C\ar 0\ssp}\subtext{old}$ \hfill and\KP9\vskip.3mm

\centerline{$\sp\Cal C\yr 0 = \seq{\,\seq{\,\sp\Cal C\ar 0\KN1\image\snn\{\sp
P\ssp\}\cap\{\,f:\dom\sn f=U\ssp\,\}:U\sn\in\tsigrd\sp P\sp\,}\subtext{old}
                                    :P\in\dom\Cal C\ar 0\ssp}\subtext{old}\,$,}\vskip.6mm

\noin then $\sp\Cal C\yr 0$ is the family associating with each pair \œ$P =
(E\ssp,F\ssp)\in\O^{\sp\times2.}=\dom\Cal C\ar 0$ the small family \œ$\sp
\taurd E\owns U\mapsto\Cal C\yr 0\fvalue\snn P\sp\fvalue U$ such that $\sp
\Cal C\yr 0\fvalue\snn P\sp\fvalue U$ is the set of functions which in
\cite[Definition 1.1\sp(b)\ssp, p.\ 219]{BGN} is denoted by \q{\spp$C^0(U,F)
$\sp}, where thus the dependence on $E$ is suppressed.

Likewise $\sp\bosy\tau\sp$ maps the pair \œ$\sp Q=(E\ar 1\sp,E\ar 2) \in
\dom\Cal C\ar 0$ to the topology $\sp\bosy\tau\sp\fvalue\snn Q\sp$ for the
underlying set \œ$\sp(\sp\vecs E\ar 1)\timesn(\sp\vecs E\ar 2)$ of the
algebraic product module \œ$\sp(\sp\sigrd E\ar 1)\vstimes(\sp\sigrd E\ar 2)$
which in \cite[Definition 1.1\sp(c)\ssp, p.\ 219]{BGN} is denoted by \q{\sp$
\Cal T(E_1\times E_2)$\sp}.

Conversely, given only the \erm{BGN}\,--\,family $\sp\Cal C\yr 0\sp$, we get\vskip.4mm

\centerline{$\sp\Cal C\ar 0 = \{\ssp(\spp P\spp,f\sp):\exi{\bosy c}\,
(\spp P\spp,\bosy c\ssp)\in\Cal C\yr 0\text{ and }\sp f \in
\bosy c\sp\fvalue(\dom\sn f\sp)\not=\Univ\sp\,\}\,$.}\vskip.4mm

\noin Hence $\sp\bosy\tau\sp$ is not needed to get $\sp\Cal C\ar 0\,$, only
its existence is needed to get the properties of $\sp
\Cal C\ar 0\vProd\,$. Note further that our $\O=$ the $
\Cal M$ in \cite[Definition 1.1\sp(a)\ssp, p.\ 219]{BGN}\ssp.

The difference quotient map $\sp\Cal C\ar 0\vbarDelta_K\tilde f\sp$ for \œ$\sp
\tilde f=(E\ssp,F\sp,f\sp)$ with \œ$\dom f=U$ in\linebreak the case where we
have that \œ$\sp\Cal C\ar 0\vbarDelta_K\tilde f\not=\Univ\sp$ is the one
determined by the mapping data \œ$
 f\yBGN1:\text{``\ssp}U\times E\times\mathbb K\text{\,"} \iinc
  U\ssp\LHB1{^{^[}}\ssn{\yr 1}\snn\LHB1{^{^]}}\to F$ in
\cite[Definition 2.1, p.\ 222]{BGN} associated with the mapping data $\sp
f:E\iinc U\to F$ corresponding to the map $\sp\tilde f\sp$.     \end{remark}


  \binsubsubhead B{Compactly generated topologies, vector and function spaces}\label{subsec B}

Here we develop those properties of compactly generated real vector spaces
which are needed in proving that Seip's higher order differentiability
concepts can be obtained as a particular case of the general construction
given above.

The most important particular result obtained of which we shall need to make
explicit use several times in the sequel is Proposition \ref{lin expo} below
giving a certain exponential law property of the space $\Lincg(E\ssp,F\ssp)$
of continuous linear maps \œ$E\to F$ topologized by the k\ssp-\ssp extension
of the compact open topology.

First, to make matters precise, we begin by introducing the following

\begin{definitions}\label{topo defs}

Balancing between \cite[p.\ 146]{Ky} and the conventional
\cite[Definition 6.1, p.\ 237]{Dugu}\ssp, we here agree on saying that a
topology $\sp\Cal T$ is {\it locally compact\ssp} if{}f for all $x\ssp,V$ with
\œ$x\in V\in\Cal T$ there are $K\sp,U$ such that \œ$x\in U\in\Cal T$ and \œ$\sp
U\inc K\inc V$ and $K$ is $\Cal T\,$--\,compact in the \q{non-Bourbakian}
sense \cite[p.\ 135]{Ky}\ssp.

Our {\it Kelleyfication\ssp} of a topology $\ssp\Cal T\sp$ is\vskip.3mm$\mhyppy4
\keltop\Cal T=\{\,U:U\inc\bigcup\ssp\Cal T$ and $\sp\all K\,\exi V$\par$\null\hfill
K$ is $\Cal T\,$--\,compact $\imply V\in\Cal T$ and $U\cap K=V\cap K\sp\,\}\,$.\KP{10}\vskip.3mm

\noin Our definition of the {\it k\ssp-\ssp extension\ssp} $k\ar{ex}\sp\Cal T\sp
$ of $\ssp\Cal T\sp$ is obtained from the preceding by putting there
\q{$\sp K$ is $\Cal T\,$--\,closed compact} in place of
\q{$\sp K$ is $\Cal T\,$--\,compact}.

A topology $\Cal T\spp$, and the topological space $\sp
(\ssp\bigcup\ssp\Cal T\spp,\Cal T\,)$ we say to be {\it almost compactly
generated\,} if{}f also \œ$\Cal T=\keltop\Cal T\sp$ holds, and we say {\it
compactly generated\,} if{}f also $\sp\Cal T$ is Hausdorff with \œ$\sp\Cal T =
\keltop\Cal T\spp$. We let \œ$\Cal T\ktimes\Cal U =
\keltop(\sp\Cal T\ttimes\Cal U\ssp)\,$, and also put \œ$\kappatv=$ $
\{\ssp(X\sp,\Cal T\sp\,;X\sp,\keltop\Cal T\ssp):X=X$ and $\sp\Cal T$ is a
topology$\sp\,\}\,$.                                       \end{definitions}

Observe that according to the preceding definitions {\it every locally compact
topology is almost compactly generated\,}, and every Hausdorff locally compact
topology is compactly generated. In view of
\cite[Theorem 5.17, p.\ 146]{Ky}\ssp, {\it every compact topology which is
also Hausdorff or regular is locally compact.}

Since \œ$\sp U\in\Cal T\sp$ whenever $\Cal T$ is a first countable topology
and $\sp U$ is sequentially $\Cal T\,$--\,open, it follows, cf.\
\cite[Theorem 7.13, p.\ 231]{Ky}\ssp, that {\it first countable topologies are
almost compactly generated\,} and {\it
             metrizable topologies are compactly generated\,}.

Note the following slight difference between our definitions above and
\cite[pp.\ 230\,--\,231]{Ky}\ssp. A topological space $\sp
(\ssp\Omega\,,\Cal T\,)$ is a {\it k\ssp-\ssp space\ssp} in the sense of
\cite{Ky} if{}f we have \œ$\sp \Cal T=k\ar{ex}\sp\Cal T\spp$, and this implies
that $\sp\Cal T\sp$ is almost compactly generated since we gener- ally have $\,
\keltop\Cal T\inc k\ar{ex}\sp\Cal T\spp$, \,equality here holding if also $\sp
\Cal T\sp$ is Hausdorff.

\begin{example}

We assume that \œ$\,[\,0\,,1\,]\inc\Omega\inc\Re\sp\,$, and we put \œ$\ssp A =
\Omega\setminus{\ssp]}\,0\,,1\,{[\,}$ and $\sp\Cal T = \{\,\bigcup\ssp\Cal A :
\Cal A\inc\Cal B\,\}$ where we have\vskip.3mm\centerline{$
\Cal B = \big\{\,{\sp]}\,s\ssp,t\,{[\sp}:0<s<t<1\ssp\big\}\cup\big\{\,{\sp
]}\,s\ssp,1\,{[\ssp}\cup\{\sp t\sp\}:0<s<1$ and $t\in A\,\big\}\,$.}\vskip.3mm

\noin Then $\Cal T$ is a first countable locally compact topology for $\ssp
\Omega\sp$ which is neither Hausdorff nor regular. Hence $\sp\Cal T\sp$ is
almost compactly generated, but if the set $A$ is infinite, then $\sp
(\ssp\Omega\,,\Cal T\,)$ is not a k\ssp-\ssp space. Indeed, then \œ$\,
\keltop\Cal T\not=k\ar{ex}\sp\Cal T\sp$ holds since we have \œ$\sp
\{\sp t\sp\}\in k\ar{ex}\sp\Cal T\sp\setminus\keltop\Cal T\sp$ for every \œ$
t\in A\,$, noting that for every $\Cal T\,$--\,closed compact set $K\sp$ there
is $\sp s\sp$ with $0<s<1$ and $\,{\sp
    ]}\,s\ssp,1\,{[\ssp}\cap\spp K=\emptyset\,$.               \end{example}

\begin{definitions}\label{k vec topo}

The class $\ssp\cgVS(\biit R\ssp)$ of real {\it compactly generated vector
space\ssp}s we let have as its members exactly the real topologized vector
spaces \œ$\sp E=(X\sp,\Cal T\ssp)$ such that $\sp\Cal T$ is compactly
generated and $\sp\sigrd\spp X$ is continuous \œ$\sp\Cal T\ktimes\Cal T \to
\Cal T$ and $\sp\taurd\spp X$ is continuous \œ$\sp
\tau_{_{I\!\!R}}\snn\ttimes\Cal T\to\Cal T\spp$. The class $\ssp
\scgVS(\biit R\ssp)$ of {\it Seip\,--\,convenient\ssp} spaces has as its
members exactly the \œ$\sp E\in\cgVS(\biit R\ssp)$ such that $\sp\Cal U\sp$
exists with $(\sp\sigrd E\ssp,\sp\Cal U\ssp)$ a sequentially complete
Hausdorff locally convex space and \œ$\sp\taurd E=\keltop\,\Cal U\ssp$. We
let \œ$E\skcap F = ((\sp\sigrd E\sp)\vstimes(\sp\sigrd F\sp))\ssp,
             (\sp\taurd E\sp)\ktimes(\sp\taurd F\sp))\,$. With every \œ$
\bosy E\in\cgVS(\biit R\ssp)\,^I$ and any set $I$ we associate the compactly
generated product vector space\vskip.3mm\centerline{$
\vsprod_cgVS_\Re\bosy E =
(\ssp\vsprod_VS_\Re(\ssp\roman{pr}\ar 1\snn\circ\bosy E\,)\ssp,
\keltop\prod{_{\sn_{\bold{ti}\,}}}(\ssp\roman{pr}\ar 2\circ\bosy E\,))\,$.}\vskip.5mm

\noin We also let $\,\taukv=\seq{\,\sp\leLCSr(\biit R\ssp)\text{\,-\ssp}\sup\,
(\sp\kappatv\KN1\inve\image\sn\{\sp F\ssp\}\cap\LCS(\biit R\ssp)) :
   F\in\cgVS(\biit R\ssp)\sp\,}\,$.                        \end{definitions}

Note that for \œ$F\in\scgVS(\biit R\ssp)$ the class \œ$
\kappatv\KN1\inve\image\sn\{\sp F\ssp\}\cap\LCS(\biit R\ssp)$ need not be a
singleton. E.g., if with \œ$\sp
E=\vscoprod_LCS_\Re(\Re\times\sn\{\sp\biit R\,\}\spp)$ we let \œ$\sp E\ar 1 =
(\sp\sigrd E\ssp,\Cal T\!\subtext{box})$ where $\sp\Cal T\!\subtext{box}$ is\Biggerlineskip2
the \q{box topology} having \œ$\big\{\,\vecs E\sp\cap\spp\prodc\bosy V :
\bosy V\in\{\,{\ssp]}\minus n^{\sp\mminus 1}\spp,n^{\sp\mminus 1}\,{[\sp} :
n\in\Zepp\ssp\}\,^{\Re}\,\big\}$ as a\Biggerlineskip1 filter basis for \œ$
\ymp E\ar 1 = \Nbh(\Re\times\sn\{0\}\ssp,\Cal T\!\subtext{box})\,$, then \œ$\sp
E\not=E\ar 1$ with \œ$\rajou E=\rajou E\ar 1$\linebreak whence compactness is
the same for $\taurd E$ and $\taurd E\ar 1\ssp$, and so \œ$\,
\kappatv\KN1\fvalue E=\kappatv\KN1\fvalue E\ar 1\ssp$.            \vskip.3mm

When \œ$F\in\scgVS(\biit R\ssp)\,$, the space $\taukv\KN1\fvalue F$ is the
strictest locally convex space looser than $F\sp$. Put otherwise, the topology
$\sp\taurd(\sp\taukv\KN1\fvalue F\ssp)$ is the finest locally convex topology
on the vector space $\sigrd F$ which is coarser than $\taurd F\sp$. In
\cite[p.\ 38\,ff.]{Se72}\ssp, the space $\taukv\KN1\fvalue F$ is denoted by
\q{\ssp LK$\,(\spp F\ssp)\sp$}, and it is sequentially complete, noting that
if \œ$E\,,E\ar 1 \in
\kappatv\KN1\inve\image\sn\{\sp F\ssp\}\cap\LCS(\biit R\ssp)$ with $E$
sequentially complete and \œ$\ssp E\le E\ar 1\ssp$, then also $E\ar 1$ is
sequentially complete.

\begin{lemma}\label{basic comp gen property}

If $\,(\ssp\Cal T\spp,\ssp\Cal U\ssp,f\ssp)$ is a global topological map,
so is $\,(\ssp\keltop\Cal T\spp,\sp\keltop\,\Cal U\ssp,f\ssp)\,$.

  \end{lemma}

\begin{proof} Under the premise, if \œ$\ssp V\sn\in\keltop\,\Cal U$ and $K$
is $\sp\Cal T\,$--\,compact, we have to prove that there is \œ$\sp U \sn\in
\Cal T$ with \œ$\sp f\sp\inve\image\spp V\cap K=U\cap K\sp$. For \œ$K\aar 1 =
f\sp\image\snn K$ now $K\aar 1$ is $\ssp\Cal U\,$--\,compact, whence there is
\œ$\ssp V\aar 1\in\Cal U$ with \œ$\sp V\aar 1\snn\cap K\aar 1=V\cap K\aar 1\ssp
$. Putting \œ$\sp U=f\sp\inve\image\spp V\aar 1\ssp$, we have \œ$\sp U\sn \in
\Cal T\spp$, and also $\sp f\sp\inve\image\spp V\cap K=U\cap K\sp$ holds.
  \end{proof}

\begin{lemma}\label{about cg product}

If $\ssp\Cal T$ and $\,\Cal U\sp$ are any topologies$\ssp,$ then \œ$\ssp
\Cal T\ktimes\Cal U=(\ssp\keltop\Cal T\,)\ktimes(\ssp\keltop\,\Cal U\,)\,$.
For any small family $\ssp\bosy{\Cal T}\ssp$ of topologies it holds that\vskip.3mm\centerline{$
\keltop\Tihprod\bosy{\Cal T}=\keltop\Tihprod\{\ssp(\sp i\ssp,\keltop\Cal T\,):
(\sp i\ssp,\Cal T\,)\in\bosy{\Cal T}\sp\,\}\,$.}                 \end{lemma}

\begin{proof} Let \œ$\sp\Cal T=\Cal T\!\ar 1$ and \œ$\,\Cal U=\Cal T\ssn\ar 2$
be topologies, and put \œ$\ssp\Cal W =
(\ssp\keltop\Cal T\,)\ktimes(\ssp\keltop\,\Cal U\,)$ and \œ$P =
\bigcup\ssp\Cal W\sp$. Since \œ$\sp\Cal T\ktimes\Cal U =
\keltop(\ssp\Cal T\ttimes\Cal U\,)$ and \œ$\ssp\Cal W =
\keltop((\ssp\keltop\Cal T\,)\ttimesn(\ssp\keltop\,\Cal U\,))\,$, by
Lemma \ref{basic comp gen property} above we trivially have \œ$\sp
\Cal T\ktimes\Cal U\inc\Cal W\sp$. To get the converse, it suffices that $
\id P$ is continuous \œ$\sp\Cal T\ktimes\Cal U\to\Cal W\sp$. This in turn
follows if for arbitrarily given \œ$\sp\Cal T\ttimes\Cal U\,$--\,compact $K$
and \œ$\smb P\in K\sp$, we have that $\id K$ is continuous \œ$\sp
\Cal T\ttimes\Cal U\to\Cal W$ at the point $\smb P\,$. That is, given \œ$
\smb P\in W\in\Cal W\sp$, there should (x) exist some \œ$\sp
U\in\Cal T\ttimes\Cal U\sp$ such that $
           \smb P\in U$ and $\ssp U\cap K\inc W\sp$ hold.

To obtain (x) above, letting \œ$B\ar 1=\dom K$ and \œ$B\ar 2=\rng K$ and \œ$\sp
C=B\ar 1\timesn B\ar 2\,$, for $\sbi{\iota\ssp=\sp\sixroman{1\sp,\ssp 2}}$
then $B\sbi\iota$ is $\sp\Cal T\ssn\sbi\iota\ssp$--\,compact, hence $\sp
\keltop\Cal T\ssn\sbi\iota\ssp$--\,compact, and consequently $C$ is \œ$
\keltop\Cal T\ttimes\keltop\,\Cal U\,$--\,compact. Hence, there is \œ$\sp
W\!\ar 1\in\keltop\Cal T\ttimesn\keltop\,\Cal U$ with \œ$\smb P\in W\cap C=$ \œ$
W\!\ar 1\cap\sp C\ssp$, whence further there are \œ$\sp V\ssn\sbi\iota \in
\keltop\Cal T\ssn\sbi\iota$ with \œ$\smb P\in V\!\ar 1\timesn V\ssn\ar 2 \inc
W\!\ar 1\ssp$. Then there are \œ$\sp U\sn\sbi\iota\in\Cal T\ssn\sbi\iota$ with
\œ$V\ssn\sbi\iota\cap B\sbi\iota=U\sn\sbi\iota\cap B\sbi\iota\,$. Putting \œ$\sp
U=U\aar 1\timesn U\aar 2\,$, we have \œ$\sp U\in\Cal T\ttimes\Cal U\sp$ with \œ$
\smb P\in V\!\ar 1\timesn V\ssn\ar 2\cap\sp C=U\aar 1\timesn U\aar 2\cap\sp C
\inc U$ and \œ$\ssp U\cap K=U\cap C\cap K\inc
V\!\ar 1\timesn V\ssn\ar 2\cap\sp C$ $\inc$ $W\!\ar 1\cap\sp C=W\cap\sp C\inc
W\spp$, \,hence (x) as we wished.

We leave it as an exercise to the reader to naturally generalize the preceding
idea in order to establish the latter assertion. As a hint we only mention
that with the family \œ$\ssp\bosy K=\{\ssp(\sp i\ssp,S\sp):\emptyset \not= S =
\{\,\xi:\exi x\,(\sp i\ssp,\xi\sp)\in x\in K\sp\,\}\sp\}\sp$ one now considers
\œ$\ssp C=\prodc\bosy K\ssp$, and one finally obtains \œ$U=\prodc\bosy U$ for
some \œ$\sp\bosy U\in\prodc\bosy{\Cal T}\sp$ with the property that $\sp\{\,i:
\bosy U\ssp\fvalue i\not=\bigcup\,(\ssp\bosy{\Cal T}\ssp\fvalue i\ssp)\ssp\}\sp
$ is a finite set.                                               \end{proof}

\begin{lemma}\label{open and clos comp gen}

Let $\ssp\Cal T$ be a compactly generated topology, and let \œ$A\inc\Omega =
\bigcup\ssp\Cal T\spp$. If in addition $\ssp A\in\Cal T$ or $\,
\Omega\sn\setminus\sn A\in\Cal T\spp,$ then $\Cal T\lei A$ is a compactly
generated topology.                                              \end{lemma}

\begin{proof} First letting \œ$\ssp\Omega\sn\setminus\sn A\in\Cal T\spp$,
arbitrarily given \œ$B\inc A$ such that for every $\sp\Cal T\lei A\,
$--\,compact $K\aar 1$ there is \œ$\sp U\sn\in\Cal T$ with \œ$\sp
U\cap K\aar 1=B\cap K\aar 1\ssp$, as for \œ$\sp V =
\Omega\sn\setminus\sn A\cup B$ we have \œ$\sp V\cap A=B\ssp$, it suffices that
\œ$\sp V\sn\in\Cal T\spp$. As $\sp\Cal T$ is compactly generated, for this it
suffices that for each given $\sp\Cal T\,$--\,compact $K$ there is \œ$\sp U\sn
\in\Cal T$ with \œ$\sp U\cap K=V\cap K\sp$. Putting \œ$K\aar 1=K\cap A\,$, now
$K\aar 1$ is $\sp\Cal T\lei A\,$--\,compact, whence there is \œ$\sp U\aar 1\in
\Cal T$ with \œ$\sp U\aar 1\snn\cap K\aar 1=B\cap K\aar 1\ssp$, hence \œ$\sp
U\aar 1\snn\cap K\cap A=B\cap K\sp$. For \œ$\sp U =
\Omega\sn\setminus\sn A\sp\cup U\aar 1$ then \œ$\sp U\snn\in\Cal T\sp$ and $\ssp
U\cap K = \Omega\snn\setminus\sn A\cup\spp U\aar 1\snn\cap K =
           \Omega\snn\setminus\sn A\cup B\cap K = V\cap K\sp$.

Next assuming that \œ$A\in\Cal T\spp$, arbitrarily given \œ$B\inc A$ such that
for every $\sp\Cal T\lei A\,$--\,compact $K\aar 1$ there is \œ$\sp U\sn \in
\Cal T$ with \œ$\sp U\cap K\aar 1=B\cap K\aar 1\ssp$, for arbitrarily fixed $\sp
\Cal T\,$--\,com- pact $K$ we must prove (k) that there is \œ$\sp U\sn \in
\Cal T$ with \œ$\sp U\cap K=B\cap K\sp$. For this consi- dering \œ$\ssp \Cal U
=\Cal T\sp\cap\{\,U:\exi x\,x\in U\cap K\inc B\cap K\sp\,\}\,$, if for
arbitrarily fixed \œ$\sp x\in B\cap K$ we prove (x) that there is \œ$\sp U\sn
\in\Cal T$ with \œ$\sp x\in U\cap K\inc B\cap K\sp$, we are done since to get
(k) we may take $\sp U=\bigcup\ssp\Cal U\ssp$.

To establish (x)\ssp, putting \œ$K\ar 2=K\setminus\sn A\,$, then $K\ar 2$ is $\sp
\Cal T\,$--\,compact whence noting that $\sp\Cal T$ is Hausdorff, there are
disjoint \œ$U\aar 1\ssp,\sp U\aar 2\in\Cal T$ with \œ$x\in U\aar 1$ and \œ$
K\aar 2\inc U\aar 2\ssp$. Putting \œ$K\aar 1=K\setminus\sn U\aar 2\ssp$, now $
K\aar 1$ is $\sp\Cal T\lei A\,$--\,compact whence there is \œ$\sp U\aar 0 \in
\Cal T$ with \œ$\sp U\aar 0\cap K\aar 1=B\cap K\aar 1\ssp$. Taking \œ$U =
U\aar 0\cap\spp U\aar 1\ssp$, we have obtained (x) since by \œ$
U\aar 1\snn\cap\spp U\aar 2$ $=\emptyset\,$ we have \ \ $x\in
B\cap K\cap\spp U\aar 1 = B\cap K\cap\spp U\aar 1\ssn\setminus U\aar 2 =
B\cap K\aar 1\cap\spp U\aar 1 = U\aar 0\cap K\aar 1\cap\spp U\aar 1$

$\mhyppy{21} = U\cap(K\setminus\sn U\aar 2)\inc U\cap K\inc
U\aar 0\cap\spp U\aar 1\snn\cap K =
U\aar 0\cap\spp U\aar 1\snn\cap K\setminus\sn U\aar 2$

$\mhyppy{21} = U\aar 0\cap\spp U\aar 1\snn\cap K\aar 1 =
B\cap K\aar 1\snn\cap\spp U\aar 1\inc B\cap K\aar 1\inc B\cap K\sp$.
  \end{proof}

\begin{corollary}\label{clos subsp in cgVS}

Let \œ$\ssp E=(X\sp,\Cal T\,)\in\cgVS(\biit R\ssp)$ and \œ$\,F=E_{\sp/\sp S}$
where $\ssp S\sp$ is a vector subspace in $\sp X$ which is $\sp\Cal T\,
$--\,closed. Then \œ$\ssp F\in\cgVS(\biit R\ssp)\,$. If in addition $\ssp E$
is Seip\,--\,convenient$\ssp,$ so is $\ssp F\sp$ too.        \end{corollary}

\begin{proof} To get \œ$\sp F\in\cgVS(\biit R\ssp)\,$, noting that by
Lemma \ref{open and clos comp gen} we have \œ$\sp\taurd F=\Cal T\lei S$
compactly generated, it suffices (a) that $\ssigrd F$ is continuous \œ$
\taurd F\ktimes\taurd F\to\Cal T\spp$,\linebreak and (b) that \œ$\tsigrd F =
\taurd\sp X\sp\,|\,(\Re\timesn S\sp)$ is continuous \œ$
\tau_{_{\Re}}\ttimes\taurd F\to\Cal T\spp$. We imme- diately get (b) from
continuity of \œ$\taurd\sp X:\tau_{_{\Re}}\ttimes\Cal T\to\Cal T\spp$, and we
deduce (a) as follows. Letting \œ$\Iota=\id(\spp S^{\sp\times 2.})\,$, since \œ$
\ssigrd F=\sigrd\sp X\circ\Iota\ssp$, it suffices that $\Iota$ is con- tinuous
\œ$\taurd F\ktimes\taurd F\to\Cal T\ktimes\Cal T\spp$, and by
Lemma \ref{basic comp gen property} this follows if $\Iota$ is continuous $
\Cal T\ttimes\Cal T\sp\lei(\spp S^{\sp\times 2.})=\taurd F\ttimes\taurd F \to
\Cal T\ttimes\Cal T\spp$. This is trivial.

Assuming that $\sp E$ is Seip\,--\,convenient, there is $\sp\Cal U\sp$ such
that \œ$(X\sp,\sp\Cal U\ssp)\in\LCS(\biit R\ssp)$ is sequentially complete
with \œ$\sp\Cal T=\keltop\,\Cal U\ssp$. Then also \œ$
(\sp\sigrd F\sp,\sp\Cal U\sp\lei S\sp)\in\LCS(\biit R\ssp)\,$, and is
sequentially complete. To prove that \œ$\sp\Cal T\lei S =
\keltop\sp(\ssp\Cal U\sp\lei S\sp)\,$, noting that we have \œ$\sp
\keltop\sp(\ssp\Cal U\sp\lei S\sp)\inc\Cal T\lei S$ by
Lemma \ref{basic comp gen property} since $\sp\Cal T\lei S$ is compactly
generated, it suffices that $\id S$ is continuous \œ$
\keltop\sp(\ssp\Cal U\sp\lei S\sp)\to\Cal T\spp$, but this is immediate by
Lemma \ref{basic comp gen property} from continuity of $
\id S:\Cal U\sp\lei S\to\Cal U\ssp$.                             \end{proof}

For sets of functions we obtain topologies, and topologized linear structures
with the aid of the following basic

\begin{constructions}\label{defi fct spaces}

$\null$

\vskip.4mm\noin(1) \ \
$\basTop\Cal B=
\{\,\bigcup\ssp\Cal A:\Cal A\inc\{\,\bigcap\ssp\Cal B\aar 1\sn:
\Cal B\aar 1\inc\Cal B$ and $
\Cal B\aar 1$ is finite$\sp\,\}\sp\}\,$,

\vskip.4mm\noin
(2) \ \ $\tFunso S\sbi{\Cal G\,}\Cal T = \basTop\{\,S\cap\{\,x:x\image\snn B
            \inc V\ssp\,\}:B\in\Cal G$ and $\ssp V\sn\in\Cal T\ssp\,\}\,$,

\vskip.4mm\noin
(3) \ \ $\tConco(\ssp\Cal T\spp,\ssp\Cal U\,) = \{\,U\sn:\all{C\ssp,\Cal K}\,
        C = (\ssp\bigcup\ssp\Cal U\,)^{\,\bigcup\sp\Cal T}\cap
            \{\,x:x\inve\images\Cal U\inc\Cal T\ssp\,\}$\nopagebreak\par\nopagebreak$\null\hfill$
and $\Cal K=\{\,K:K$ is $\sp\Cal T\,$--\,compact$\sp\,\}\imply U\sn\in
\tFunso C\sbi{\Cal K}\,\sp\Cal U\ssp\,\}\,$,\KP{8}

\vskip.4mm\noin%
(4) \ \ $\roman{vF}\!\ar{so\,}S\sbi{\Cal G}\ssp F = \big((\sp\sigrd
         F\ssp)\expnota^\roman{dom}\ssp\bigcup\sp S]_{vs}{_{|\ssp S}}
                        \,,\sp \tFunso S\sbi{\Cal G\,}\taurd F\,\big)\,$,

\vskip.4mm\noin
(5) \ \ $\roman{vF}\!\ar{su\,}S\sbi{\Cal G}\ssp F = \bigcap\ssp\{\ssp((\sp
        \sigrd F\ssp)\expnota^\Omega]_{vs}{_{|\ssp S}}\ssp,\Cal T\,):F\sp,S:
        \Omega=\dom\sn\bigcup\ssp S$ and\nopagebreak\par\nopagebreak$\mhyppy{28}
   \Cal T =
   \{\,U:U\inc S$ and $\sp\all{x\in U}\,\exi{B\in\Cal G\ssp,\sp V\in\ymp F}$\nopagebreak\par\nopagebreak$\null\hfill
   \all{u\in S}\,u\image\snn B \inc
              V\imply\ssigrd F\circ[\,x\ssp,u\,] \in U\sp\,\}\sp\}\,$.\KP{8}

\end{constructions}

One immediately observes that $\ssp\basTop\Cal B\sp$ is a topology for $\ssp
\bigcup\ssp\Cal B\sp$ whenever \œ$\sp\Cal B\in\Univ\,$. Consequently $\ssp
\tFunso S\sbi{\Cal G\,}\Cal T$ is a topology whenever $S$ is a set. However,
it may be quite pathological or uninteresting, unless $\ssp\Cal T$ is a
topology and \œ$\,\bigcup\ssp\Cal G\inc\Omega$ and \œ$S \inc $ $
(\ssp\bigcup\ssp\Cal T\,)\,^\Omega\sp$ for some $\sp\Omega\,$, and possibly
some further conditions hold.

In particular, if $\ssp\Cal T$ and $\,\Cal U$ are topologies, then $\ssp
\tConco(\ssp\Cal T\spp,\ssp\Cal U\,)$ is the {\it compact open\ssp} topology
for the set $\bigcup\ssp\tConco(\ssp\Cal T\spp,\ssp\Cal U\,)$ of functions $x$
with $(\ssp\Cal T\spp,\ssp\Cal U\ssp,x\sp)$ a global topological map.
Informally $\ssp\roman{vF}\!\ar{su\,}S\sbi{\Cal G}\ssp F$ may be called the
set $S$ of functions into $F$ turned into the space with topology that of the
uniform convergence on the members of $\ssp\Cal G\sp$.

\begin{proposition}\label{Jar fct space}

Let \œ$\ssp F\in\LCS(\biit R\ssp)$ and \œ$\,\Omega=\bigcup\ssp\Cal G\in\Univ$
and \œ$\ssp X=(\sp\sigrd F\ssp)\expnota^\Omega]_{vs}\ssp$. Let $\ssp S$ be a
vector subspace in $\sp X$ such that \œ$\ssp\{\,x\image\snn B:x\in S\text{ and }\sp
B\in\Cal G\,\}\inc\rajou F\sp,$ and such that for any \œ$\sp A\,,B\in\Cal G$
there is \œ$\,C\in\Cal G$ with \œ$\ssp A\cup B\inc C\ssp$. Let \œ$\ssp E =
\roman{vF}\!\ar{su\,}S\sbi{\Cal G}\ssp F\sp$. Then \œ$\ssp E\in\LCS(\biit R\ssp)$ with
\œ$\ssp\sigrd E=X_{\sp|\ssp S}\,$. If $\ssp\Cal V$ is any filter basis for $\sp
\ymp F\sp,$ then the class $\ssp\{\,S\cap\{\,x:x\image\snn B\inc V\sp\,\}:B\in
\Cal G\text{ and }\ssp V\sn\in\Cal V\,\}$ is a filter basis for $\sp\ymp E\,$.

  \end{proposition}

\begin{proof} See \cite[Proposition 2.10.1, pp.\ 43\,--\,44]{Jr}\ssp.
  \end{proof}

Adapting the proof of \cite[Theorem 7.11, p.\ 230]{Ky}\ssp, we obtain

\begin{proposition}\label{comp open = comp unif}

Let \œ$\ssp F\in\LCS(\biit R\ssp)\,,$ and let $\sp(\sp\Omega\,,\Cal T\,)$ be
any topological space with \œ$\ssp\Cal K=\{\,K:K\sp\text{ is $\sp\Cal T\,
$--\,compact}\sp\,\}\,$. Also let \œ$\ssp\Cal F=\tConco(\ssp\Cal T\spp,
\taurd F\ssp)\,,$ and let $\sp S$ be a vector subspace in $\sp
(\sp\sigrd F\ssp)\expnota^\Omega]_{vs}$ with $\sp S\inc\bigcup\ssp\Cal F\sp$.
Then $\sp\taurd\roman{vF}\!\ar{su\,}S\sbi{\Cal K}\ssp F =
            \Cal F\spp\lei S\ssp$.                         \end{proposition}

\begin{proof} Let \œ$\ssp\Cal U=\taurd\roman{vF}\!\ar{su\,}S\sbi{\Cal K}\ssp F
$ and \œ$\ssp\Cal P\aar 2=\Cal K\timesn\taurd F\sp$. Also putting \œ$
N\ssp(\sp x\,;K\sp,V\ssp)=$ \œ$\{\,\ssigrd F\circ[\,x\ssp,u\,]:u\in S\text{
and }u\image\snn K\inc V\ssp\,\}\,$, then by Proposition \ref{Jar fct space}
above, for \œ$x\in S$ the class $\{\,N\ssp(\sp x\,;K\sp,V\ssp):K\in\Cal K$ and
$\ssp V\sn\in\ymp F\sp\,\}$ is a filter basis for $
                    \Nbh(\sp x\ssp,\ssp\Cal U\,)\,$.

To prove that \œ$\ssp\Cal U\inc\Cal F\spp\lei S\ssp$, arbitrarily given \œ$\sp
x\in U\in\Cal U\ssp$, there are some \œ$K\in\Cal K$ and \œ$\ssp V\sn\in\ymp F$
with \œ$\sp N\ssp(\sp x\,;K\sp,V\ssp)\inc U\spp$, and we should establish (x)
existence of a finite \œ$\sp\Cal P\inc\Cal P\aar 2$ such that for \œ$\sp N =
S\cap\{\,y:\all{P\in\Cal P}\,y\image\sigrd P\inc\taurd P\sp\,\}$ we have $\sp
x\in N\inc N\ssp(\sp x\,;K\sp,V\ssp)\,$. Let $\sp
             \tilde x=(\ssp\Cal T\spp,\taurd F\sp,x\sp)\,$.

To get (x)\ssp, we first take a $\sp\taurd F\,$--\,closed \œ$\ssp V\aar 1 \in
\ymp F$ with \œ$\ssp V\aar 1\snn-V\aar 1\snn-V\aar 1\sn\inc V\spp$, and put \œ$\ssp
W=\{\ssp(\sp\zeta\,,\xi\sp):\zeta-\xi\in V\aar 1\ssp\}\,$. Now $\sp
x\image\snn K$ is $\sp\taurd F\,$--\,compact by continuity of $\ssp\tilde x$
and $\,\Cal T\,$--\,compactness of $\sp K\sp$, whence there is a finite \œ$
A\inc K$ with \œ$x\image\snn K\inc W\circ x\sp\image\ssn A\,$. \linebreak
Taking \œ$\Cal P=\{\ssp(K\cap(\ssp x\inve\snn\circ W\circ
x\sp\image\snn\{\sp\smb Q\spp\}\sp)\,,\sp\Int_taurd F(\ssp W\circ W\circ
x\sp\image\sn\{\sp\smb Q\spp\}\sp)):\smb Q\in A\sp\,\}\,$, we\linebreak have
$\sp\Cal P\sp$ finite with \œ$\sp\Cal P\inc\Cal P\aar 2$ since by continuity
of $\ssp\tilde x$ and $\ssp\taurd F\,$--\,closedness of $\,V\aar 1$ the set $
x\inve\snn\circ W\circ x\sp\image\snn\{\sp\smb P\ssp\}$ is $\ssp
\Cal T\,$--\,closed for every $\smb P\in A\,$. $\null\hfill$ To see that $\sp
x\in N\sp$, just note that if $\sp\smb P\in K\cap(\ssp x\inve\snn\circ W\circ
x\sp\image\snn\{\sp\smb Q\spp\}\sp)$ with $\smb Q\in A\,$, we then have

\centerline{$ x\fvalue\smb P\in x\,[\,x\inve\snn\circ W\circ
x\sp\image\snn\{\sp\smb Q\spp\}\,]\inc W\circ x\sp\image\snn\{\sp\smb Q\spp\}
\inc\Int_taurd F(\ssp W\circ W\circ x\sp\image\sn\{\sp\smb Q\spp\}\sp)\,$.}

To prove that \œ$\sp N\inc N\ssp(\sp x\,;K\sp,V\ssp)\,$, arbitrarily fixing \œ$\sp
y\in N$ and \œ$\smb P\in K\sp$, it suffices to show that \œ$\sp
(\sp y-x\sp)\fvalue\smb P\in V\spp$. Now, by \œ$\smb P\in K \inc
x\inve\snn\circ W\circ x\sp\image\ssn A\,$, there is \œ$\smb Q\in A$\linebreak
with $\smb P\in x\inve\snn\circ W\circ x\sp\image\snn\{\sp\smb Q\spp\}\,$,
consequently $\smb P\in K\cap(\ssp x\inve\snn\circ W\circ
x\sp\image\snn\{\sp\smb Q\spp\}\sp)\,$. Letting\vskip.3mm

\centerline{$ P = (K\cap(\ssp x\inve\snn\circ W\circ
x\sp\image\snn\{\sp\smb Q\spp\}\sp)\,,\sp
\Int_taurd F(\ssp W\circ W\circ x\sp\image\sn\{\sp\smb Q\spp\}\sp))\,$,}\vskip.3mm

\noin then \œ$P\in\Cal P\sp$ holds, and since by \œ$\sp y\in N$ we have \œ$y
\image\sigrd P\inc\taurd P\sp$, we get \œ$y\fvalue\snn\smb P\in$ \œ$
\Int_taurd F(\ssp W\circ W\circ x\sp\image\sn\{\sp\smb Q\spp\}\sp) \inc
W\circ W\circ x\sp\image\sn\{\sp\smb Q\spp\}\,$. We also have \œ$
x\fvalue\snn\smb P\in W\circ x\sp\image\sn\{\sp\smb Q\spp\}\,$, and hence \œ$
x\fvalue\snn\smb Q\in W\sp\inve\snn\circ x\sp\image\sn\{\sp\smb P\ssp\}\,$.
Consequently, we get \œ$(\sp x\fvalue\smb P\ssp,y\fvalue\smb P\ssp) \in
W\circ W\circ W\sp\inve\sp$, that is $\sp(\sp y-x\sp)\fvalue\smb P =
y\fvalue\smb P-x\fvalue\smb P\in V\aar 1\snn-V\aar 1\snn-V\aar 1\sn\inc V\spp$.

Conversely, to prove that \œ$\ssp\Cal F\spp\lei S\inc\Cal U\ssp$, arbitrarily
given a nonempty \œ$K\in\Cal K$ and \œ$U\in\taurd F$ and \œ$x\in S$ with \œ$
C=x\image\snn K\inc U\spp$, it suffices to find \œ$\ssp V\in\ymp F$ with\linebreak
\œ$N\ssp(\sp x\,;K\sp,V\ssp)\inc S\sp\cap\{\,y:y\image\snn K\inc U\sp\,\}\,$.
Noting that this immediately follows if we have \œ$C+V\inc U\spp$, it suffices
to establish the latter. For this, we consider \œ$\Cal N\aar 1=$ \œ$
C\timesn\ymp F\sp\cap\{\ssp(\sp\xi\,,V\ssp):\xi+V+V\inc U\sp\,\}\,$. Since $
(\sp\ssigrd F\sp,\taurd F\ssp)$ is a topological group, we have \œ$
C\inc\dom\Cal N\aar 1\ssp$, and since $C$ is $\taurd F\,$--\,compact, a finite
\œ$\Cal N\inc\Cal N\aar 1$ exists with \œ$\sp C\inc\{\,\xi+\zeta:\exi V\,
(\sp\xi\,,V\ssp)\in\Cal N\text{ and }\zeta\in V\sp\,\}\,$. Taking now \œ$\sp
V=\bigcap\sp\rng\Cal N\sp$, we have $\ssp V\in\ymp F$ with $
C+V\inc\{\,\xi+\zeta+\zeta\ar 1\sn:\exi V\,(\sp\xi\,,V\ssp)\in\Cal N\text{ and }
\zeta\,,\zeta\ar 1\in V\sp\,\}\inc U\spp$.                       \end{proof}

Observe above that when establishing the inclusions in both directions between 
$\ssp\Cal U\ssp$ and $\sp\Cal F\spp\lei S\ssp$, we needed to know that $
(\sp\ssigrd F\sp,\taurd F\ssp)$ is a topological group. The preceding proof
hence does not give a corresponding result if we assume $F$ to be a compactly
generated vector space instead of a topological one. In fact, if the
topologized linear structure of \œ$\sp F\in\cgVS(\biit R\ssp)$ determines a
uniformity in the natural manner,
             then $F$ necessarily is a topological vector space.

\begin{proposition}\label{ev conti T ktim C to U}

Let $\,\Cal T$ and $\sp\,\Cal U\sp$ be any topologies such that {\ssp\rm(a)}
or {\ssp\rm(b)} below holds. Then $\,\roman{ev}\sp\,|\,\bigcup\ssp\Cal W\sp$
is continuous $\,\Cal W\to\Cal U\ssp$.\vskip.3mm{\ssp\rm
(a)} \ \ $\Cal T\sp$ is locally compact and $\,
           \Cal W=\Cal T\ttimes\tConco(\ssp\Cal T\spp,\ssp\Cal U\ssp)\,,$\par{\ssp\rm
(b)} \ \ $\Cal T\sp$ is Hausdorff and $\,\Cal W =
   \Cal T\ktimes\tConco(\ssp\Cal T\spp,\ssp\Cal U\ssp)\,$. \end{proposition}

\begin{proof} We give the deduction assuming (b) leaving it as an exercise to
the reader to extract therefrom the almost trivial case when (a) is assumed
instead. Put- ting \œ$\Cal F=\tConco(\ssp\Cal T\spp,\ssp\Cal U\ssp)\,$,
arbitrarily given a \œ$\sp\Cal T\ttimes\Cal F\,$--\,compact $K\aar 2$ and \œ$\ssp
V\sn\in\sp\Cal U\ssp$, there should be some \œ$\sp W\sn\in\Cal T\ttimes\Cal F$
with \œ$\roman{ev}\inve\image\spp V\cap K\aar 2=W\cap K\aar 2\ssp$. Supposing
we have (x) that for any \œ$\smb Z\in\roman{ev}\inve\image\spp V\cap K\aar 2$
there are \œ$N\in\Cal T$ and \œ$U\sn\in\Cal F$ with $\smb Z\in N\timesn U$ and\linebreak
\œ$\roman{ev}\sp\,[\,N\snn\times\snn U\spp\cap\spp K\aar 2\,]\inc V\spp$, we
are done since $\sp W=\bigcup\ssp\Cal A\sp$ will do for\vskip.3mm\centerline{$
\Cal A=\{\,N\snn\times\snn U:(N\sp,U\ssp)\in\Cal T\snn\times\snn\Cal F$ and $\,
\roman{ev}\sp\,[\,N\snn\times\snn U\spp\cap\spp K\aar 2\,]\inc V\ssp\,\}\,$.}\vskip.3mm

To get (x) we arbitrarily fix \œ$\smb Z=(\sp\smb P\ssp,x\sp) \in
\roman{ev}\inve\image\spp V\cap K\aar 2\ssp$, and first note that by\linebreak
continuity of \œ$x:\Cal T\to\Cal U$ there is $N\aar 0$ with \œ$\smb P \in
N\aar 0\in\Cal T$ and \œ$x\image\snn N\aar 0\inc V\spp$. For \œ$K\aar 0=$ $
\dom K\aar 2$ and \œ$K\aar 1=K\aar 0\ssn\setminus\sn N\aar 0\,$, then $K\aar 0$
and $K\aar 1$ are $\sp\Cal T\,$--\,compact. Now consider- ing \œ$\sp
\Cal N\aar 1=\Cal T^{\ssp\times2.}\cap\spp\{\ssp(N\aar 1\sp,N\aar 2):\exi{
\smb Q}\,(\sp\smb P\ssp,\smb Q\spp)\in N\aar 1\sn\times\sn N\aar 2\text{ and }
N\aar 1\cap N\aar 2=\emptyset\,\}\,$, since $\sp\Cal T$ is Hausdorff and \œ$
\smb P\not\in K\aar 1\ssp$, we have \œ$K\aar 1\inc\bigcup\sp\rng\Cal N\aar 1\ssp
$, whence by compactness there is a finite \œ$\Cal N\inc\Cal N\aar 1$ with \œ$
K\aar 1\inc\bigcup\sp\rng\Cal N\sp$. Taking \œ$N=\bigcap\sp\dom\Cal N$ and \œ$
K=K\aar 0\ssn\setminus\snn\bigcup\sp\rng\Cal N$ and \œ$U =
\bigcup\ssp\Cal F\spp\cap\spp\{\,u:u\image\snn K\inc V\ssp\,\}\,$, then \œ$
\smb P\in N\in\Cal T$ and $K$ is $\sp\Cal T\,$--\,compact with \œ$
N\cap K\aar 0\inc K\inc N\aar 0\,$, and consequently \œ$U\sn\in\Cal F$ with \œ$
\smb Z\in N\snn\times\snn U\spp$. In addition, we also have \œ$\ssp
\roman{ev}\sp\,[\,N\snn\times\snn U\spp\cap\spp K\aar 2\,]\inc V\sp$ since if
\œ$\sp\smb Z=(\sp\smb Q\ssp,u\sp)\in N\snn\times\snn U\spp\cap\spp K\aar 2\,$,
then we get $\smb Q\in N\cap\dom K\aar 2=N\cap K\aar 0\inc K$ and $
u\image\snn K\inc V\spp$, whence $\ssp\roman{ev}\fvalue\smb Z=u\fvalue\smb Q
\in V\spp$.                                                      \end{proof}

The preceding might be compared to one half of \cite[Lemma 1.13, p.\ 13]{Se72}\ssp.

\begin{proposition}\label{conti expo law}

Let $\ssp(\sp\Omega\,,\ssp\Cal U\,)\ssp,(\sp T\spp,\Cal T\,)\ssp,
(\sp T\aar 1\sp,\Cal T\!\ar 1)$ be any topological spaces such that {\ssp\rm%
(a)} or {\ssp\rm(b)} below holds$\ssp,$ and let \œ$\sp f \in
T\aar 1\sp^{\Omega\ssp\times\ssp T}\sp$. Then $\sp f$ is continuous \œ$\,
\Cal W\to\Cal T\!\ar 1$ if and only if $\sp f\yvee$ is continuous $\,
\Cal U\to\tConco(\ssp\Cal T\spp,\Cal T\!\ar 1)\,$.\vskip.3mm{\ssp\rm
(a)} \ \ $\Cal T\sp$ is locally compact and $\,
           \Cal W=\Cal U\ttimes\Cal T\spp,$\par{\ssp\rm
(b)} \ \ $\Cal T\sp$ is Hausdorff and $\,\keltop\,\Cal U=\Cal U\sp$ and $\,
           \Cal W = \Cal U\ktimes\Cal T\spp$.              \end{proposition}

\begin{proof} Again leaving to the reader the simpler case where (a) holds,
under the assumption (b) putting \œ$
\Cal F=\tConco(\ssp\Cal T\spp,\Cal T\!\ar 1)\,$, we first assume that $f$ is
continuous \œ$\,\Cal W\to\Cal T\!\ar 1\ssp$, and proceed to prove that $f\yvee
$ is continuous \œ$\,\Cal U\to\Cal F$ as follows. Ar- bitrarily given a $\sp
\Cal U\,$--\,compact $K\sp$, it suffices that $f\yvee\ssp|\,K$ is continuous \œ$
\Cal U\lei K\to\Cal F\sp$. For this, arbitrarily fixing \œ$\smb P\in K\sp$,
and a $\Cal T\,$--\,compact $B$ and \œ$\ssp V\in\Cal T\!\ar 1$ such that \œ$\sp
f\,[\,\{\sp\smb P\ssp\}\sn\times\sn B\sp\,]\inc V\spp$, it suffices to find \œ$
N\in\Nbh(\sp\smb P\ssp,\ssp\Cal U\lei K\ssp)$ with \œ$\sp
f\,[\,N\snn\times\sn B\sp\,]\inc V\spp$. To establish this, we consider\vskip.3mm

$\mhyppy7 \Cal N\aar 1 = \{\ssp(N\sp,U\ssp):\exi x\,N \in
                          \Nbh(\sp\smb P\ssp,\ssp\Cal U\lei K\ssp)$ and\nopagebreak\par\nopagebreak%
$\mhyppy{41.3} U\sn\in\Nbh(\sp x\ssp,\Cal T\lei B\sp)$ and $\sp
              f\,[\,N\snn\times\snn U\ssp\,]\inc V\ssp\,\}\,$.\vskip.5mm

Letting \œ$K\aar 2=K\snn\times\sn B\ssp$, as $K\aar 2$ is $\ssp
\Cal U\ttimes\Cal T\,$--\,compact, by continuity of \œ$
f:\Cal W\to\Cal T\!\ar 1\ssp$,\linebreak we obtain continuity of \œ$\ssp
f\,|\,K\aar 2\sn:(\ssp\Cal U\ttimes\Cal T\,)\lei K\aar 2\to\Cal T\!\ar 1\ssp$.
Using this, we see that\linebreak \œ$B\inc\bigcup\sp\rng\Cal N\aar 1\ssp$,
whence by compactness of $B$ there is a finite \œ$\Cal N\inc\Cal N\aar 1$ with
\œ$B\inc$ $\bigcup\sp\rng\Cal N\ssp$. Taking \œ$N=\bigcap\sp\dom\Cal N\ssp$,
we have \œ$N\in\Nbh(\sp\smb P\ssp,\ssp\Cal U\lei K\ssp)\,$, and one quickly
checks that also $\sp f\,[\,N\snn\times\sn B\sp\,]\inc V\sp$ holds.

Conversely, next assuming that $f\yvee$ is continuous \œ$\,\Cal U\to\Cal F\sp
$, we get continuity of \œ$f:\Cal W\to\Cal T\!\ar 1$ as follows. Putting \œ$
\Iota=\{\ssp(\sp\smb P\ssp,x\,;x\ssp,\smb P\ssp):\smb P\in\Omega\text{ and }
x\in T\sp\,\}$ and \œ$\Eps =$ $
\roman{ev}\,|\,(\ssp T\timesn\bigcup\ssp\Cal F\,)\,$, we have $
f=\Eps\circ(\id T\risti2 f\yvee\spp)\circ\Iota$ with continuous maps\vskip.4mm

$\mhyppy{24.5}\Iota:\Cal U\ttimes\Cal T\to\Cal T\ttimes\Cal U\sp$, \,hence
$\,\Cal W=\Cal U\ktimes\Cal T\to\Cal T\ktimes\Cal U\sp$, \hfill and\vskip.3mm

$\mhyppy{8.5}(\id T\risti2 f\yvee\spp):\Cal T\ttimes\Cal U \to
\Cal T\ttimes\Cal F\sp$, \,hence $\,\Cal T\ktimes\Cal U \to
\Cal T\ktimes\Cal F\sp$, \hfill and\vskip.5mm

\noin $\Eps:\Cal T\ktimes\Cal F\to\Cal T\!\ar 1$ by
Proposition \ref{ev conti T ktim C to U} above.
                         The assertion follows.                  \end{proof}

\begin{corollary}\label{coro 1 of conti exp law}

Let $\ssp(\sp\Omega\,,\ssp\Cal U\,)\ssp,(\sp T\spp,\Cal T\,)\ssp,
(\sp T\aar 1\sp,\Cal T\!\ar 1)$ be any topological spaces such that $\,\Cal U$
is almost compactly generated and $\,\Cal T\sp$ is locally compact$\ssp,$ and
let \œ$\sp f \in T\aar 1\sp^{\Omega\ssp\times\ssp T}$ and \œ$\,\Cal P =
(\ssp\Cal U\ttimes\Cal T\spp,\Cal T\!\ar 1)\,$. If \œ$\,
(\sp\Cal P\spp,f\,|\,(K\snn\times\Univ\sp))$ is a topological map for all $\,
\Cal U\,$--\,compact $\ssp K\sp,$ then also $\,(\sp\Cal P\spp,f\ssp)$ is a
topological map.                                             \end{corollary}

\begin{proof} Assuming the premise and putting $\Cal F=\tConco(\ssp\Cal T\spp,
\Cal T\!\ar 1)\,$, by Proposition \ref{conti expo law}\ssp(a) it suffices that
$(\ssp\Cal U\sp,\spp\Cal F\sp,f\yvee)$ is a topological map, which in turn
follows if $(\ssp\Cal U\lei K\sp,\spp\Cal F\sp,f\yvee\ssp|\,K\ssp)$ is a such
for all $\ssp\Cal U\,$--\,compact $K\sp$. Again by
Proposition \ref{conti expo law}\ssp(a) it suffices that $
(\sp f\yvee\ssp|\,K\ssp)\ywed=f\sp\,|\,(K\snn\times\Univ\sp)$ is continuous $
(\ssp\Cal U\lei K\ssp)\ttimes\Cal T\to\Cal T\!\ar 1\ssp$, but this is
immediate from the premise.                                      \end{proof}

\begin{corollary}\label{2 about cg product}

If $\,\Cal T$ and $\,\Cal U\sp$ are any topologies of which one is compactly
genera- ted and the other is locally compact$\ssp,$ then $\ssp
\Cal T\ktimes\Cal U = \Cal T\ttimes\Cal U\ssp$ holds.        \end{corollary}

\begin{proof} Letting \œ$P=\bigcup\ssp\Cal T\timesn\bigcup\,\Cal U$ and \œ$
\Iota\ar 0=\{\ssp(\sp x\ssp,y\,;\sp y\ssp,x\sp):(\sp x\ssp,y\sp)\in P\sp\,\}\,
$, we have\linebreak \œ$
(\ssp\Cal T\ttimes\Cal U\ssp,\spp\Cal U\ttimes\Cal T\spp,\spp\Iota\ar 0)$ a
homeomorphism, hence also \œ$
(\ssp\Cal T\ktimes\Cal U\ssp,\sp\Cal U\ktimes\Cal T\spp,\spp\Iota\ar 0)$ by
Lemma \ref{basic comp gen property} above. Consequently, if we prove the
assertion for $\ssp\Cal U$ compactly generated and $\,\Cal T\sp$ locally
compact, it also holds with the roles reversed. Now, letting the roles be as
suggested, it suffices that $\id P$ is continuous \œ$\Cal T\ttimes\Cal U \to
\Cal T\ktimes\Cal U\sp$, equivalently that $\Iota\ar 0\KN1\inve$ is continuous
\œ$\ssp\Cal U\ttimes\Cal T \to \Cal T\ktimes\Cal U\sp$, and for this by
Corollary \ref{coro 1 of conti exp law} it suffices that for every $\,\Cal U\,
$--\,compact $K$ we have \œ$\Iota\ar 0\KN1\inve\ssp|\,(K\snn\times\Univ\sp)$
continuous \œ$\ssp\Cal U\ttimes\Cal T \to \Cal T\ktimes\Cal U\sp$,
equivalently \œ$\id(\ssp\bigcup\ssp\Cal T\timesn K\ssp)$ continuous \œ$
\Cal T\ttimes(\ssp\Cal U\lei K\ssp) \to \Cal T\ktimes\Cal U\sp$. For this,
since by \cite[Theorem 5.17, p.\ 146]{Ky} the topology $\ssp\Cal U\lei K$ is
locally compact in the sense of our Definitions \ref{topo defs} above, further
by Corollary \ref{coro 1 of conti exp law} it suffices that for every $\sp
\Cal T\,$--\,com- pact $B$ we have $\id(B\timesn K\ssp)$ continuous $
\Cal T\ttimes\Cal U \to \Cal T\ktimes\Cal U\sp$, but this is trivial.
  \end{proof}

Our Corollary \ref{2 about cg product} above may be compared to the slightly
weaker result in \cite[Theorem 4.4, p.\ 263]{Dugu} with less formal and six
lines longer presentation of proof where both topologies are required to be
Hausdorff.

\begin{proposition}\label{basic Seip conv char}

It holds that \œ$\,\kappatv\KN1\image\LCS(\biit R\ssp)\inc\cgVS(\biit R\ssp)\,
$. For every sequentially complete $\sp E\in\LCS(\biit R\ssp)$ it holds that $\,
\kappatv\KN1\fvalue E\in\scgVS(\biit R\ssp)\,$.            \end{proposition}

\begin{proof} Arbitrarily given \œ$\sp E=(\sp a\ssp,c\,,\Cal T\ssp) \in
\LCS(\biit R\ssp)\,$, for \œ$\ssp\Cal U=\keltop\Cal T$ we have to prove that $\sp
a$ is continuous \œ$\ssp\Cal U\ktimes\Cal U\to\Cal U$ and $\ssp c\ssp$ is
continuous \œ$\sp\tau_{_{I\!\!R}}\snn\ttimes\Cal U\to\Cal U\ssp$. Using
Lemmas \ref{basic comp gen property} and \ref{about cg product} and
Corollary \ref{2 about cg product} above, one easily deduces these.
  \end{proof}

In particular, from the preceding proposition we see that the Seip convenient
spaces are exactly the $\sp\kappatv\KN1\fvalue E\sp$ for some sequentially
complete real Hausdorff locally convex space $E\,$. Again using
Lemma \ref{about cg product} and noting that $\ssp\tau_{_{\Re}}$ and $\,
\tau_{_{\Re}}\snn\ttimes\tau_{_{\Re}}$ are compactly generated, and that
sequential completeness is preserved when taking products of small families of
locally convex spaces, we get the following

\begin{corollary}\label{Seip vs products}

Let \œ$\,\O=\scgVS(\biit R\ssp)\,$. Then \œ$\,\bosy R\in\O\ssp,$ and for any $\,
\bosy E\,,E\,,F\sp,I$ the implications $\,E\,,F\in\O\imply E\skcap F\in\O$ and
$\,\bosy E\in\O\,^I\imply\vsprod_cgVS_\Re\bosy E\in\O$ hold. \end{corollary}

We shall need the spaces of continuous linear maps given in the following

\begin{constructions}\label{fct spaces of lin. maps}

$\null$\vskip.4mm\noin
(1) \ \ $\Cal L\,(E\ssp,F\ssp)=\{\,\ell:\ell$ linear $\sigrd E\to\sigrd F$ and
        $\ssp\ell\ssp\inve\images\taurd F\inc\taurd E\sp\,\}\,$,\vskip.4mm\noin
(2) \ \ $\Cal L\ar{co}(E\ssp,F\ssp)=\bigcap\ssp\{\,\roman{vF}\!\ar{so\,}
        \Cal L\,(E\ssp,F\ssp)\sbi{\Cal G}\ssp F:E:\Cal G=\{\,K:K$ is $\sp
                      \taurd E\,$--\,compact$\sp\,\}\sp\}\,$,\vskip.4mm\noin
(3) \ \ $\Link(E\ssp,F\ssp)=\bigcap\ssp\{\,\roman{vF}\!\ar{su\,}
        \Cal L\,(E\ssp,F\ssp)\sbi{\Cal G}\ssp F:E:\Cal G=\{\,K:
            K$ is $\sp\taurd E\,$--\,compact$\sp\,\}\sp\}\,$,\vskip.4mm\noin
(4) \ \ $\Lincg(E\ssp,F\ssp)=\bigcap\ssp\{\,\kappatv\KN1\fvalue
         \roman{vF}\!\ar{so\,}C\sbi{\Cal G}\ssp F_{\sp/\ssp S}:F:
              C = \sp\bigcup\ssp\tConco(\sp\taurd E\ssp,\taurd F\ssp)$ and \par$\null\hfill
         \Cal G = \{\,K:K$ is $\sp\taurd E\,$--\,compact$\sp\,\}$ and $
              S = \Cal L\,(E\ssp,F\ssp)\ssp\}\,$.\KP{5.7}

\end{constructions}

\begin{lemma}\label{lin is clos}

Let \œ$\ssp\Cal T=\tConco(\sp\taurd E\ssp,\taurd F\ssp)\sp$ where \œ$\ssp F\in
\cgVS(\biit R\ssp)$ and $\,E\sp$ is any real topologized vector space. Then $\ssp
\Cal L\,(E\ssp,F\ssp)$ is $\ssp\keltop\Cal T\,$--\,closed.       \end{lemma}

\begin{proof} For \œ$\ssp\roman E\,x\,y\,t =
\seq{\,(\sp u\fvalue(\sp x+y\sp) - u\fvalue x-u\fvalue y\ssp,
u\fvalue(\sp t\,x\sp)-t\,(\sp u\fvalue x\sp)):u\in\bigcup\ssp\Cal T\ssp\,}\,$,\linebreak
we have \œ$\Cal L\,(E\ssp,F\ssp) =
\bigcap\ssp\{\ssp(\sp\roman E\,x\,y\,t\sp)\inve\image\snn
\{\ssp(\sp\bnull F\ssp,\bnull F)\ssp\}:x\ssp,y\in\vecs E\text{ and }t \in
\Re\sp\,\}\,$. For fixed \œ$z\in\vecs E\,$, directly by definition, we see
that $\ssp\roman{ev}\sbi{\eightmath z}\ssp|\,\bigcup\ssp\Cal T\sp$ is
continuous \œ$\sp\Cal T\to\taurd F\sp$. For \œ$x\ssp,y\in\vecs E$ and \œ$t \in
\Re\sp\,$, letting \œ$f=\seq{\,(\sp u\fvalue(\sp x+y\sp)\ssp,u\fvalue x\ssp,
u\fvalue y\ssp,u\fvalue(\sp t\,x\sp)):u\in\bigcup\ssp\Cal T\ssp\,}$ and \œ$\,
\Cal U\sp=\taurd F\ttimes\taurd F\ttimes\taurd F\ttimes\taurd F\sp$, it
follows that $f$ is continuous \œ$\Cal T\to\Cal U\ssp$, and by
Lemma \ref{basic comp gen property} hence \œ$\keltop\Cal T\to\keltop\,\Cal U\ssp
$. Taking into account that \œ$F\in\cgVS(\biit R\ssp)\,$, and using
Lemma \ref{about cg product} one deduces that \œ$\ssp\roman E\,x\,y\,t$ is
continuous $\ssp\keltop\Cal T\to\taurd F\ttimes\taurd F\sp$. The assertion of
the lemma now follows.                                           \end{proof}

\begin{proposition}\label{space of lin map properties}

Let $\ssp E\ssp,F\in\cgVS(\biit R\ssp)$ and $\ssp F\!\ar 1\in\LCS(\biit R\ssp)\,
$. Then \begin{enumerate}\vskip.4mm

\item \ \ $\Lincg(E\ssp,F\ssp)=\kappatv\KN1\fvalue\Cal L\ar{co}(E\ssp,F\ssp)
            \in\cgVS(\biit R\ssp)\,,$\vskip.4mm

\item \ \ $\Link(E\ssp,F\!\ar 1) =
           \Cal L\ar{co}(E\ssp,F\!\ar 1)\in\LCS(\biit R\ssp)\,,$\vskip.4mm

\item \ \ $\Lincg(E\ssp,\kappatv\KN1\fvalue F\!\ar 1) =
            \kappatv\KN1\fvalue\Link(E\ssp,F\!\ar 1)\in\cgVS(\biit R\ssp)\,,$\vskip.4mm

\item \ \ $F\in\scgVS(\biit R\ssp) \imply
             \Lincg(E\ssp,F\ssp)\in\scgVS(\biit R\ssp)\,$.

\end{enumerate} \end{proposition}

\begin{proof} For (1) letting \œ$\sp L=\Lincg(E\ssp,F\ssp)$ and \œ$\,\Cal U =
\tConco(\sp\taurd E\ssp,\taurd F\ssp)\,$, by an elementary set theoretic
verification the proof of the asserted equality is reduced to getting \œ$\ssp
\keltop\,\Cal U\sp\lei\vecs L=\keltop(\ssp\Cal U\sp\lei\vecs L\sp)\,$. For \œ$\ssp
\keltop(\ssp\Cal U\sp\lei\vecs L\sp)\inc\keltop\,\Cal U\sp\lei\vecs L\,$, we
note that trivially $\idv L$ is continuous \œ$\sp\keltop\,\Cal U\sp\lei\vecs L
\to\Cal U\sp\lei\vecs L\,$. By Lemmas \ref{lin is clos} and
\ref{open and clos comp gen} and \ref{basic comp gen property} it is also
continuous \œ$\sp
\keltop\,\Cal U\sp\lei\vecs L\to\keltop(\ssp\Cal U\sp\lei\vecs L\sp)\,$,
whence the inclusion. For the converse using
Lemma \ref{basic comp gen property} note the continuity implications \œ$\ssp
\Cal U\sp\lei\vecs L\to\Cal U\ssp$, hence $\ssp
\keltop(\ssp\Cal U\sp\lei\vecs L\sp)\to\keltop\,\Cal U\ssp$, hence $\ssp
\keltop(\ssp\Cal U\sp\lei\vecs L\sp)\to\keltop\,\Cal U\sp\lei\vecs L\,$.

To get \œ$\sp L\in\cgVS(\biit R\ssp)$ in (1)\ssp, also let \œ$\sp G =
\kappatv\KN1\fvalue\roman{vF}\!\ar{so\,}\bigcup\ssp\tConco(\sp\taurd E\ssp,
\taurd F\ssp)\sp\sbi{\Cal K\,}F\sp$ where \œ$ \Cal K =
\{\,K:K$ is $\sp\taurd E\,$--\,compact$\sp\,\}\,$. If \œ$
G\in\cgVS(\biit R\ssp)$ holds, since \œ$L=G_{\sp/\ssp\upsilon_s\sp L}\,$, we
immediately get \œ$\sp L\in\cgVS(\biit R\ssp)$ by Lemma \ref{lin is clos} and
Corollary \ref{clos subsp in cgVS} above. It hence suffices to establish \œ$
G\in\cgVS(\biit R\ssp)$ which in turn follows if $
(\ssp\Cal T\ktimes\Cal T\spp,\Cal T\spp,\ssigrd\spp G\ssp)$ and $
(\sp\tau_{_{I\!\!R}}\snn\ttimes\Cal T\spp,\Cal T\spp,\tsigrd\spp G\ssp)$ are
topological maps when $\sp\Cal T=\taurd\spp G=\keltop\,\Cal U\ssp$.

To prove continuity of \œ$\sp\ssigrd\spp G:\Cal T\ktimes\Cal T\to\Cal T\sp$,
by Proposition \ref{conti expo law} and Lemma \ref{basic comp gen property} it
suffices that $\sp\ssigrd\spp G\ywed$ is continuous \œ$
\taurd\sp(\spp G\skcap G\skcap E\ssp)\to\taurd F\sp$. For \œ$\sp\Eps =
\roman{ev}\,|\,\vecs(E\sp\sqcap G\ssp)$ and \œ$\sp\Iota\ar 2 =
\{\ssp(\sp u\ssp,v\ssp,x\ssp,(\sp x\ssp,u\,;x\ssp,v\sp)):u\ssp,v\in\vecs G\text{
and }x\in\vecs E\sp\,\}\,$, observing that we\linebreak have \œ$\sp
\ssigrd\spp G\ywed =
\ssigrd\spp F\spp\circ(\sp\Eps\risti2\Eps\sp)\circ\Iota\ar 2\,$, since $\Eps$
is continuous \œ$\taurd\sp(E\skcap G\ssp)\to\taurd F\sp$ by
Proposition \ref{ev conti T ktim C to U} above, it suffices to note that \œ$
(\spp G\skcap G\skcap E\,,\spp E\skcap G\skcap(E\skcap G\ssp)\ssp,\spp
\Iota\ar 2)$ and $(\spp F\skcap F\sp,F\sp,\ssigrd\spp F\ssp)$
                       are continuous linear maps.

For continuity of \œ$\sp\tsigrd\spp G:\tau_{_{I\!\!R}}\snn\ttimes\Cal T \to
\Cal T\sp$, as above it suffices that $\sp\tsigrd\spp G\ywed$ is con- tinuous
\œ$\sp\taurd(\biit R\sp\skcap G\skcap E\ssp)\to\taurd F\sp$, now noting
that $\sp\tau_{_{I\!\!R}}\snn\ttimes\Cal T\sp$ is compactly gener- \linebreak
ated by Corollary \ref{2 about cg product} above. Since by a twofold
application of Lemma \ref{basic comp gen property} with\linebreak $\Iota\ar 1
=$ \œ$\{\ssp(\sp t\ssp,u\ssp,x\ssp,(\sp t\,;x\ssp,u\sp)):t\in\Re\text{ and }u
\in\vecs G\text{ and }x\in\vecs E\sp\,\}$ we have $\sp\Iota\ar 1$ continuous
\œ$\sp\taurd(\biit R\sp\skcap G\skcap E\ssp)\to
\taurd(\biit R\sp\skcap(E\skcap G\ssp))$ with \œ$\sp\tsigrd\spp G\ywed =
\tsigrd\spp F\spp\circ(\sp\id\Re\risti2\Eps\sp)\circ\Iota\ar 1\ssp$,
continuity of $\sp(\sp\tau_{_{I\!\!R}}\snn\ttimes\Cal T\spp,\Cal T\spp,
\tsigrd\spp G\ssp)\sp$ then follows.

We get (2) by applying elementary set theoretic manipulations to
Constructions \ref{defi fct spaces} and \ref{fct spaces of lin. maps} above,
and using Propositions \ref{Jar fct space} and \ref{comp open = comp unif}
with $F\aar 1$ in place of $F$ and also taking $\ssp\Omega=\vecs E$ and $\sp
S=\Cal L\,(E\ssp,F\aar 1)\,$.

For (3) letting \œ$\sp L\ar 1=\Lincg(E\ssp,\kappatv\KN1\fvalue F\aar 1)$ and \œ$\ssp
L\ar 0=\Cal L\ar{co}(E\ssp,F\aar 1)$ and \œ$\ssp L=\kappatv\KN1\fvalue L\ar 0\,
$, and noting (2) we should have \œ$\sp L\ar 1=L\sp$ since we get the
membership \œ$\sp L\ar 1\in\cgVS(\biit R\ssp)$ from (1) by
Proposition \ref{basic Seip conv char} above. Since by
Lemma \ref{basic comp gen property} we have \œ$\sp\vecs L\ar 0\inc\vecs L\ar 1\ssp
$, we triv- ially get \œ$\sp\sigrd L=\sigrd L\ar 0=\sigrd L\ar 1\ssp$, and
hence it remains to establish \œ$\sp\taurd L\ar 1=\taurd L\,$. To get the
simpler \œ$\sp\taurd L\inc\taurd L\ar 1\ssp$, we should have (u) that $\idv L$
is continuous \œ$\sp\taurd L\ar 1\to\taurd L\,$. Now, for \œ$\sp\Cal T\aar 0 =
\tConco(\sp\taurd E\ssp,\taurd F\aar 1)$ and \œ$\sp\Cal T\aar 1 =
\tConco(\sp\taurd E\ssp,\keltop\taurd F\aar 1)$\linebreak we trivially have $
\id\bigcup\ssp\Cal T\aar 1$ continuous \œ$\sp
\keltop\Cal T\aar 1\to\Cal T\aar 0\,$, and hence $\idv L$ continuous\linebreak
\œ$\sp\taurd L\ar 1=\keltop\Cal T\aar 1\lei\vecs L\to\Cal T\aar 0\,$, whence
also \œ$\sp\taurd L\ar 1\to\Cal T\aar 0\lei\vecs L=\taurd L\ar 0\,$. By
Lemma \ref{basic comp gen property} the required (u) immediately follows.

For \œ$\sp\taurd L\ar 1\inc\taurd L\sp$ in (3) we show that $\idv L$ is
continuous \œ$\sp\taurd L\to\taurd L\ar 1\ssp$. Having \œ$\sp\taurd L\ar 1 =
\keltop\Cal T\aar 1\lei\vecs L$ and $\sp\taurd L$ compactly generated (again
by Lemma \ref{basic comp gen property} above) it suffices to establish
continuity \œ$\sp\taurd L\to\Cal T\aar 1\ssp$. By
Proposition \ref{conti expo law}\ssp(b) above, this fol-\linebreak lows if we
show that $\idv L\ywed$ is continuous \œ$\sp\taurd(L\skcap E\ssp) \to
\keltop\taurd F\aar 1$ which in turn follows from continuity \œ$\sp
\taurd(L\skcap E\ssp)\to\taurd F\aar 1\ssp$. Noting that $\idv L\ywed$ is
defined by \œ$(\ssp\ell\,,x\sp)\mapsto\ell\ssp\fvalue x\ssp$, and putting \œ$\sp
P= \vecs(E\sp\sqcap L\sp)$ and \œ$\ssp Q= \vecs E\times\bigcup\ssp\Cal T\aar 0$
and \œ$\ssp\Eps=\roman{ev}\,|\,Q\ssp$, it suffices (e) that $\ssp\Eps\,|\,P$
is continuous \œ$\sp\taurd(E\skcap L\sp)\to\taurd F\aar 1\ssp$. Now, from
Proposition \ref{ev conti T ktim C to U}\ssp(b) we know that $\ssp\Eps$ is
continuous \œ$\sp\taurd E\ktimes\Cal T\aar 0\to\taurd F\aar 1\ssp$, and hence
that $\ssp\Eps\,|\,P$ is continuous \œ$\sp\taurd E\ktimes\Cal T\aar 0\lei P\to
\taurd F\aar 1\ssp$. Consequently (e) follows if $\ssp\id P$ is continuous \œ$\sp
\taurd(E\skcap L\sp)\to\taurd E\ktimes\Cal T\aar 0\lei P\sp$, following from
continuity \œ$\sp\taurd(E\skcap L\sp)\to\taurd E\ttimes\Cal T\aar 0\,$, which
indeed is the case by \ $\taurd E\ttimes(\sp\Cal T\aar 0\lei\vecs L\sp) =
\taurd E\ttimes(\sp\taurd L\ar 0)\inc$ \par $\mhyppy{40}
\taurd E\ttimes(\sp\keltop\taurd L\ar 0) = \taurd E\ttimes(\sp\taurd L\sp)
                       \inc\taurd(E\skcap L\sp)\,$.\vskip.3mm

For (4) assuming \œ$F\in\scgVS(\biit R\ssp)\,$, a sequentially complete \œ$
F\aar 1\in\LCS(\biit R\ssp)$ exists with \œ$F=\kappatv\KN1\fvalue F\aar 1$
whence by (3) we have \œ$\Lincg(E\ssp,F\ssp)=\kappatv\KN1\fvalue\Link(E\ssp,
F\aar 1)\,$. By a suitable adaptation of the proof in
\cite[Theorem 7.12, p.\ 231]{Ky} or \cite[Satz 1.43, p.\ 22]{Se72}\ssp, the
locally convex space $\Link(E\ssp,F\aar 1)$ is sequentially complete. Hence
the assertion follows from Proposition \ref{basic Seip conv char} above.
  \end{proof}

\begin{proposition}\label{lin expo}

Let \œ$\ssp E\ssp,F\sp,G\in\cgVS(\biit R\ssp)$ and \œ$\,U\sn\in\taurd E$ and \œ$\sp
f\aar 1\in\big((\sp\vecs G\ssp)^{\,\upsilon_s\sp F\,}\big){^{\,U}}$ be such
that $\sp f\aar 1\KN1\fvalue x$ is linear \œ$\ssp\sigrd F\to\sigrd G$ for all
\œ$\ssp x\in U\spp$. Then $\sp f\aar 1\KN1\ywed$ is continuous $\ssp
\taurd(E\spp\skcap F\ssp)\to\taurd G$ if and only if $\ssp f\aar 1$ is
continuous $\ssp\taurd E\to\taurd\Lincg(\spp F\sp,G\ssp)\,$.

  \end{proposition}

\begin{proof} Let \œ$\sp L=\Lincg(\spp F\sp,G\ssp)$ and \œ$\ssp\Cal T =
\tConco(\sp\taurd F\sp,\taurd G\ssp)$ and \œ$\,W = U\timesn\vecs F\sp$. By
Lemma \ref{open and clos comp gen} then \œ$\sp\taurd(E\spp\skcap F\ssp)\lei W$
and $\ssp\taurd E\sp\lei\sp U$ are compactly generated. In view of
Lemma \ref{basic comp gen property} it follows that \œ$\sp
\taurd(E\spp\skcap F\ssp)\lei W = \taurd E\sp\lei\sp U\ktimes\taurd F\ssp$.
Using also Proposition \ref{conti expo law}\ssp(b) in $\ssp
\overset{_{\sp!}}\equivv\ssp$ below, we get the equivalences\vskip.7mm $\mhyppy{18}
f\aar 1\KN1\ywed$ is continuous $\ssp\taurd(E\spp\skcap F\ssp)\to\taurd G$\vskip.3mm$\mhyppy{12}
\equivv\, f\aar 1\KN1\ywed$ is continuous $\ssp
                               \taurd(E\spp\skcap F\ssp)\lei W\to\taurd G$\vskip.3mm $\mhyppy{12}
\equivv\, f\aar 1\KN1\ywed$ is continuous $\ssp
                           \taurd E\sp\lei\sp U\ktimes\taurd F\to\taurd G$\vskip.3mm $\mhyppy{12}
\overset{_{\sp!}}\equivv\, f\aar 1$ \ is continuous $\ssp
                                            \taurd E\sp\lei\sp U\to\Cal T$\vskip.3mm $\mhyppy{12}
\equivv\, f\aar 1$ \ is continuous $\ssp
                                     \taurd E\sp\lei\sp U\to\keltop\Cal T$\vskip.3mm $\mhyppy{12}
\equivv\, f\aar 1$ \ is continuous $\ssp
               \taurd E\to\keltop\Cal T\lei\vecs L=\taurd L\,$.  \end{proof}

\begin{definition}

The class of continuous bilinear maps of Seip\,--\,convenient spaces is\vskip.5mm

$\mhyppy{1.4}
\RHB{1.5}{\fiveroman{bi}\sn}\Cal L\sp\subtext{scgVS}\spp(\biit R\ssp)=
\{\ssp(E\spp\skcap F\sp,G\sp,b\ssp):E\ssp,F\sp,G\in\scgVS\sp(\biit R\ssp)$ and

$\mhyppy{9} b$ function $\vecs(E\sp\sqcap F\sp)\to\vecs G$ and
              $b$ continuous $\taurd(E\spp\skcap F\sp)\to\taurd G$

\hyppy{29mm} and $\sp\all{x\ssp,y}\,[\ x\in\vecs E\imply
b\,(\sp x\ssp,\cdot\sp)$ linear $\sigrd F\to\sigrd G\ ]$

\hyppy{40.5mm} and $[\ y\in\vecs F\imply
b\,(\ssp\cdot\ssp,y\sp)$ linear $\sigrd E\to\sigrd G\ ]\sp\,\}\,$.

  \end{definition}

\begin{proposition}\label{ev conti}

Let \œ$\ssp E\ssp,F\in\cgVS(\biit R\ssp)$ and \œ$\,G =
E\spp\skcap\Lincg(E\ssp,F\ssp)$ and \œ$\ssp\Eps=\roman{ev}\,|\,\vecs G\sp$.
Then $\ssp\Eps$ is continuous \œ$\sp\taurd G\to\taurd F\sp$. If in addition \œ$\ssp
E\ssp,F\in\scgVS(\biit R\ssp)\,,$ the mem- bership $\ssp
(\spp G\sp,F\sp,\Eps\sp) \in
\RHB{1.5}{\fiveroman{bi}\sn}\Cal L\sp\subtext{scgVS}\spp(\biit R\ssp)\sp$
also holds.                                                \end{proposition}

\begin{proof} Let \œ$\sp L=\Cal L\ar{co}(E\ssp,F\ssp)$ and \œ$\ssp\Cal T =
\tConco(\sp\taurd E\ssp,\taurd F\ssp)\,$. By
Proposition \ref{ev conti T ktim C to U}\ssp(b) we have \œ$\,
\roman{ev}\,|\,(\sp\vecs E\times\bigcup\ssp\Cal T\,)$  continuous \œ$\sp
\taurd E\ktimes\Cal T\to\taurd F\sp$, and hence $\sp\Eps$ continuous \œ$\sp
\taurd E\ktimes\Cal T\spp\lei\sp\vecs G\to\taurd F\sp$. For the first
assertion it hence suffices (c) that $\sp\idv G\sp$ is continuous \œ$\sp
\taurd G\to\taurd E\ktimes\Cal T\spp\lei\sp\vecs G\sp$. Recalling
Lemma \ref{about cg product} and and using
(1) of Proposition \ref{space of lin map properties} we see that \œ$\,
  \taurd\sp G = \keltop(\sp\taurd E\ttimes\keltop\taurd L\sp) =
  \keltop(\sp\taurd E\ttimes\taurd L\sp) = $ $
  \keltop(\sp\taurd E\ttimes(\sp\Cal T\lei\vecs L\sp))\,$, whence by
Lemma \ref{basic comp gen property} then (c) follows. \hfill For the latter
as- sertion just note also
(4) of Proposition \ref{space of lin map properties} above.      \end{proof}

\begin{lemma}\label{seq conv in L_co(F,G)}

Let \œ$\ssp F\in\cgVS(\biit R\ssp)$ and \œ$\ssp G\in\LCS(\biit R\ssp)\,$. If
\œ$\,\bosy\ell\to\ell$ in top $\sp\taurd\Link(\spp F\sp,G\ssp)\,,$ then also $\,
\bosy\ell\to\ell$ in top $\sp\taurd\Lincg(\spp F\sp,\kappatv\KN1\fvalue G\ssp)
\,$. \hfill {\rm(\sp sequential convergences\,!\sp)}\ \ssp\null  \end{lemma}

\begin{proof} Let \œ$H=\kappatv\KN1\fvalue G$ and \œ$L =
\Cal L\ar{co}\spp(\spp F\sp,H\sp)$ and \œ$\sp\Cal T =
\taurd\Cal L\ar{co}\spp(\spp F\sp,G\ssp)\,$. By
Propo- sition \ref{basic Seip conv char} and (1) and (2) of
Proposition \ref{space of lin map properties} then \œ$\sp\Cal T =
\taurd\Link(\spp F\sp,G\ssp)$ and \œ$\ssp\kappatv\KN1\fvalue L=$\linebreak $
\Lincg(\spp F\sp,H\ssp)\,$. Since the topologies $\sp\taurd\sp L$ and $\sp
\keltop\taurd\sp L$ have the same compact sets, also the {\it sequential\ssp}
convergence \œ$\,\bosy\ell\to\ell$ means the same in both topologies. Hence,
assuming (c) that \œ$\,\bosy\ell\to\ell$ in top $\sp\Cal T\spp$, it
suffices to prove that \œ$\,\bosy\ell\to\ell$ in top $\sp\taurd\sp L\,$. For
this, arbitrarily given a $\taurd F\,$--\,compact $K$ and $\sp V$ with \œ$\ssp
\ell\ssp\image\snn K\inc V\sn\in\taurd H\sp$, we should establish (e)
existence of \œ$\smb N\in\No$ such that for \œ$\,\roman C\,\smb N\sp = \ssp
\bigcup\,(\ssp\bosy\ell\,\image(\No\ssn\setminus\sn\smb N\sp))\image\snn K\sp$
we have \œ$\,\roman C\,\smb N\inc V\spp$. To get this, putting \œ$B =
\rng\bosy\ell\ssp\cup\sp\{\ssp\ell\ssp\}$ and \œ$A =
\bigcup\ssp B\sp\image\snn K\spp$, then $B$ is $\sp\Cal T\,$--\,compact.
Considering \œ$\Eps=\roman{ev}\,|\,(K\snn\times\bigcup\ssp\Cal T\,)\,$, like
in Proposition \ref{ev conti T ktim C to U}\ssp(a) one verifies that $\Eps$ is
continuous \œ$\taurd F\ttimes\Cal T\to\taurd G\sp$. Also noting that \œ$A =
\Eps\,[\,K\snn\times B\sp\,]\,$,\linebreak it follows that $A$ is $\taurd G\,
$--\,compact. Consequently, there is \œ$U\in\taurd G$ such that \œ$U\cap A =
V\cap A\,$. As also \œ$\ssp\ell\ssp\image\snn K\inc A$ holds, we have \œ$\ssp
\ell\ssp\image\snn K\inc V\cap A = U\cap A\inc U\spp$, and hence by (c) there
is \œ$\smb N\in\No$ with \œ$\,\roman C\,\smb N\inc U\spp$. Also having \œ$\,
\roman C\,\smb N\inc A\,$, we get $\,\roman C\,\smb N\inc U\cap A =
V\cap A\inc V\spp$, \,as required in (e) above.                  \end{proof}


  \binsubsubhead C{Riemann integration of curves in topologized vector spaces}\label{subsec C}

In some proofs below, we shall need to integrate a (continuous) curve on the
unit interval $[\,0\,,1\,]$ taking values in a compactly generated or locally
convex space. To be able to treat these in a unified manner accordant with our
general approach in this paper, we here give a short account on this matter.

\begin{definition}\label{defi Riem int}

Let $I$ be a nontrivial real compact interval, i.e.\ for some \œ$a\ssp,b \in
\Re$\linebreak with \œ$a<b\,$, we have \œ$I = [\,a\ssp,b\,] =
\{\,t:a\le t\le b\,\}\,$. Let \œ$E=(X\sp,\Cal T\ssp)$ be a real topologized
vector space. We then say that $\gamma$ on $I$ in $E$ is {\it Riemann
integrable\ssp} to $x$ if and only if $\gamma$ is a function with \œ$I \inc
\dom\gamma$ and \œ$\rng\gamma\inc\vecs E\,$, and \œ$x\in\vecs E$ is such that
for every \œ$\sp V\sn\in\Nbh(\sp x\ssp,\Cal T\ssp)$ there is \œ$\delta\in\Rep$
with the property that for all \œ$k\in\N$ and for all \œ$\biit t \in
I\,^{k\ssp+\sp 1.}$ and \œ$\biit s\in I\,^k$ satisfying\vskip.3mm\centerline{$
I=[\,\biit t\fvalue\emptyset\,,\biit t\fvalue k\,]$ \ and \ $k\inc\{\,i :
\biit t\fvalue i\le\biit s\fvalue i\le\biit t\fvalue(\sp i\ssp\yplus)$ and $
\biit t\fvalue(\sp i\ssp\yplus)-(\sp\biit t\fvalue i\sp)<\delta\sp\,\}\,$,}\vskip.4mm

\noin we have \hfill $\sigrd E\text{\,-\sp}\sum_{\,i\ssp\in\ssp k}\sp
(\sp\biit t\fvalue(\sp i\ssp\yplus)-(\sp\biit t\fvalue i\sp))\,
(\sp\gamma\circ\biit s\spp\fvalue i\sp)\in V\spp$. \hfill We also define\vskip.4mm\centerline{$
E\text{\,-\sp}\int_{\,a}^{\,b}\gamma=\bigcap\ssp\{\,x:\all z\,x = z \equivv
\gamma$ on $[\,a\ssp,b\,]$ in $E$ is Riemann integrable to $z\,\}\,$,}\vskip.4mm

\noin and \œ$E\text{\,-\sp}\int_{\,a}^{\,b}\mathfrak F\ssp\,\roman d\,
\mathfrak s = E\text{\,-\sp}\int_{\,a}^{\,b}\sp\seq{\,\mathfrak F:\mathfrak s
\in[\,a\ssp,b\,]\ssp\,}\subtext{old}$ whenever $\mathfrak F$ is any term, and
$\mathfrak s$ is any\biggerlineskip4 variable. In a connection where
\q{\sn$E$\spp} may be considered as being implicitly understood, we may drop
\q{\snn$E$\,-\sp} from the notation.                        \end{definition}

We do not formulate as an explicit proposition the obvious fact that if $\sp
\gamma$ on $I$ in $F$ is Riemann integrable to $x$ and to $y\ssp$, then \œ$x=y
$ in the case where $\taurd F$ is a Hausdorff topology. Even without any
Hausdorff assumption, either we have that $\sp\gamma$ on $[\,a\ssp,b\,]$ in $F
$ is Riemann integrable to $F\text{\,-\sp}\int_{\,a}^{\,b}\gamma\,$, or else
$F\text{\,-\sp}\int_{\,a}^{\,b}\gamma=\Univ\,$. The latter in particular holds
if $\gamma$ is integrable to no $x$ or to several such.

\begin{lemma}

Let \œ$\ssp(\spp F\aar 1\sp,F\ssp)\in\kappatv\ssp|\,\LCS(\bosy R\sp)\,$. Then
$\sp\gamma$ on $I$ in $F$ is Riemann inte- grable to $\ssp x\sp$ if and only
if the same holds in the space $\sp F\aar 1\ssp$.                \end{lemma}

\begin{proof} By \œ$\ssp\taurd F\aar 1\inc\taurd F\ssp$, integrability in $F$
trivially implying that in $F\aar 1\ssp$, it suf- fices to show that a
contradiction follows if $\gamma$ on $I$ in $F\aar 1$ is Riemann integrable to
$x$ but not in $F\sp$. Indeed, then there is \œ$\sp V\sn\in\Nbh(\sp x\ssp,
\taurd F\sp)$ such that for every \œ$i\in\No\ssp$, considering \œ$\delta =
(\sp i\ydot\snn+1\sp)^{\sp\mminus 1}\sp$, there are $\biit t\ssp,\biit s$ with
the properties in Definition \ref{defi Riem int} such that for the Riemann sum
\œ$x\ar 1 = \sigrd F\text{\,-\sp}\sum_{\,l\ssp\in\ssp k}\sp
(\sp\biit t\fvalue(\ssp l\,\yplus)-(\sp\biit t\fvalue l\ssp))\,
(\sp\gamma\circ\biit s\spp\fvalue l\ssp)$ we\linebreak have \œ$x\ar 1 \not\in
V\spp$. By the axiom of choice we obtain \œ$\biit x \in
(\sp\vecs F\aar 1\ssn\setminus\sn V\ssp)\potNo$ with $\biit x\fvalue i$
corresponding to the $x\ar 1$ in the preceding for every \œ$i\in\No\ssp$. By
integrability in $F\aar 1\ssp$, we have \œ$\biit x\to x$ in top $\taurd
F\aar 1\ssp$. Putting \œ$K=\rng\biit x\sp\cup\{\sp x\sp\}\,$, it holds that $K$
is $\taurd F\aar 1\ssp$--\,compact, hence also $\taurd F\,$--\,compact. From \œ$
\taurd F\lei K=\taurd F\aar 1\snn\lei K\sp$, it follows that for some \œ$
V\aar 1\in\taurd F\aar 1$ we have \œ$x\in V\aar 1\snn\cap K\inc V\cap K\sp$,
and consequently for large $i$ it holds that $\biit x\fvalue i\in
V\aar 1\snn\cap K\inc V\cap K\inc V$ in contradiction with $\biit x\fvalue i
\not\in V$ for all $i\ssp$.                                      \end{proof}

\begin{proposition}\label{conti is int}

Let \œ$\sp F\in\LCS(\biit R\ssp)$ be sequentially complete. Let $I\sn$ be
a notrivial real compact interval and let \œ$\ssp\gamma\in(\sp\vecs F\ssp)\,^I
$ be continuous \œ$\sp\tau_{_{I\!\!R}}\to\taurd F\sp$. Then there is $x$ such
that $\sp\gamma$ on $I$ in $F$ is Riemann integrable to $x\ssp$.

  \end{proposition}

\begin{proof} Letting \œ$I=[\,a\ssp,b\,]\,$, and considering \œ$\biit x\in(\sp
\vecs F\ssp)\potNo$ defined by\vskip.7mm\centerline{$
i\mapsto\sigrd F\text{\,-\sp}\sum_{\,l\ssp\in\ssp i^+}
(\sp i\ydot\snn+1\sp)^{\sp\mminus 1}\sp(\sp b-a\sp)\,(\sp\gamma\fvalue(\sp a +
l\sp\,(\sp i\ydot\snn+1\sp)^{\sp\mminus 1}\sp(\sp b-a\sp)))\,$,}\vskip.7mm

\noin using continuity of $\sp\gamma$ and compactness of $\sp I$ and local
convexity of $\sp F\sp$, one first verifies that $\biit x$ is a Cauchy
sequence in $F\sp$, whence by the assumed sequential completeness there is $\sp
x$ with \œ$\biit x\to x$ in top $\taurd F\sp$. By the same token, one verifies
that $\gamma$ inte- grates to $x\ssp$, the details again being left to the
reader.                                                          \end{proof}

\begin{lemma}\label{int-mvT}

Let $\ssp U\sn$ be closed and convex in $\ssp E\sp$ with \œ$\ssp E \in
\LCS(\bosy R\sp)\,$. If also $\ssp\gamma$ on $\ssp[\,0\,,1\,]$ in $\sp E$ is
Riemann integrable to $\sp y\sp$ with $\sp\rng\gamma\inc U\spp,$ then $\,
y\in U\spp$.                                                     \end{lemma}

\begin{proof} Put \œ$I=[\,0\,,1\,]\,$, and first note the following. If \œ$
\varphi\in{\ssp]}\,0\,,\plusinfty\,{[\,\ssp}^I$ is inte- \Biggerlineskip1
grable in the sense of Lebesgue, then \œ$
\int_{\,I}\varphi\ssp\,\roman d\,\mu\subtext{Leb}>0\,$. Riemann integrability
implying that of Lebesgue, if \œ$\sp r\in\Re\sp$ and \œ$\sp
\chi\in{\ssp]}\,r\ssp,\plusinfty\,{[\,\ssp}^I$ is Riemann integrable, consid-
ering \œ$\varphi=\seq{\,\sp\chi\fvalue t-r:t\in I\sp\,}\,$, we get \œ$
 \bosy R\ssp\text{\,-}\int_{\,0}^{\,1}\chi =
  \int_{\,I}\chi\ssp\,\roman d\,\mu\subtext{Leb} =
   r+\int_{\,I}\varphi\ssp\,\roman d\,\mu\subtext{Leb}>r\ssp$. \Biggerlineskip1
Now, under the premise of the lemma, if we have \œ$y\not\in U\spp$, by
Hahn\,--\,Banach, there \linebreak are \œ$\sp r\in\Re\sp$ and \œ$\,\ell \in
\Cal L\ssp(E\ssp,\biit R\sp)$ with \œ$\ssp\ell\ssp\image\snn U \inc
{\sp]}\,r\ssp,\plusinfty\,{[\sp}$ and \œ$\,\ell\ssp\fvalue y=r\ssp$, \,and the
prece- ding applied to $\sp\chi=\ell\circ\gamma\sp$ gives $\,
r = \ell\ssp\fvalue y = \bosy R\ssp\text{\,-}\int_{\,0}^{\,1}\chi > r\ssp$,
\,a contradiction.                                               \end{proof}

Note that the preceding Hahn\,--\,Banach argument of the proof fails if the
convex set $\sp U$ is required to be neither closed nor open. Call a set \œ$\sp
C\inc\vecs E$ {\it countably convex\ssp} if{}f for any \œ$\bosy x\in C\potNo$
and \œ$\ssp\bosy s\in[\,0\,,1\,]\potNo$ and any $x$ with \œ$\sum\ssp\bosy s=1$
and \œ$
\Seq{\ssp\sum_{\,l\ssp\in\ssp i^+}\bosy s\fvalue l\,(\sp\bosy x\fvalue l\ssp)
 : i\in\No\ssp}\to x$ in top $\taurd E$ we have \œ$x\in C\ssp$. {\it Does the
conclu- sion in\ssp} Lemma \ref{int-mvT} {\it remain valid if $\,U\sn$ is
required to be merely countably convex\,\sp}?

For example, let \œ$\ssp\roman e\,i = \{\ssp(\sp\emptyset\,,\minus 1\sp)\ssp,
(\sp i\ssp,1\sp)\ssp\}\cup((\N\setminus\sn\{\sp i\sp\})\timesn\{0\})$ for \œ$
i\in\N\ssp$, and put \œ$\ssp\roman e\,\emptyset =
\{\ssp(\sp\emptyset\,,1\sp)\ssp\}\cup(\N\timesn\{0\})\,$, and in the Fr\'echet
space \œ$E=(X\sp,\Cal T\,)=\bosy R\expnota^\iNo]_{tvs}$
 con- sider the countably convex set\vskip.3mm\centerline{$
U = \{\sp\{\ssp(\sp\emptyset\,,2\,\bosy s\fvalue\emptyset-\sum\ssp\bosy s\sp)
\ssp\}\cup(\sp\bosy s\,|\,\N\sp):\bosy s\in[\,0\,,1\,]\potNo$ and $\ssp
                 0<\sum\ssp\bosy s\le 1\ssp\}$}\vskip.3mm

\noin for which we have \œ$\bnull E\in\ClT\ssp U\setminus\sn U\spp$. If $\sp
\gamma$ on $[\,0\,,1\,]$ in $E$ is Riemann (or even \q{sca- larwise} Lebesgue)
integrable to $x$ with $\rng\gamma\inc U\spp$,
              is then $x=\bnull E$ impossible\,? \vskip.3mm

Note that there does not exist \œ$\ssp\ell\in\Cal L\,(E\ssp,\bosy R\ssp)$ with
\œ$\ssp\ell\,\image U\inc\Rep$ since for every such \œ$\ssp\ell\sp$ there is
\œ$\ssp y\in\vecss\vscoprod_VS_\Re(\No\timesn\{\sp\bold R\sp\})$ with \œ$\,
\ell\ssp=\sp\Seq{\,\sum_{\,i\ssp\in\ssp\iNo}x\fvalue i\,(\sp y\fvalue i\sp)
 : x\in\vecs E\sp\,}\,$, and hence otherwise for a sufficiently large $i\in\N
$ we would get\vskip.3mm\centerline{$
0 < \ell\ssp\fvalue\roman e\,\emptyset = y\fvalue\emptyset =
     \minus\ell\ssp\fvalue\roman e\,i < 0\,$, \,a contradiction.}

\begin{definitions}\label{defi curve dif:ty}

Let \œ$E=(X\sp,\Cal T\ssp)$ be a real topologized vector space, and consider \œ$
c:J\to\vecss X\sp$ with \œ$t\ar 0\in J\inc\Re$ such that also $0$ is a $
\tau_{_{\Re\,}}$--\,limit point of the set \œ$
\Re\sp\cap\{\,t:t\ar 0+t\in J\,\}\setminus\{0\}\,$. We then say that $\sp c$
is {\it differentiable to $x$ in $E$ at\ssp} $t\ar 0$ if and only if also \œ$
x\in\vecss X$ and for every \œ$\sp V\sn\in\Nbh(\sp x\ssp,\Cal T\ssp)$ there is
\œ$\delta\in\Rep$ with the property that for all \œ$t\in\Re$ with \œ$
t\ar 0+t\in J$ and \œ$0<|\ssp t\ssp|<\delta$ we have the mem- bership $\,
t^{\sp\mminus 1}\sp(\sp c\fvalue(\sp t\ar 0+t\sp)-c\fvalue t\ar 0)\in V\spp$.
We let\vskip.3mm\centerline{$
\cder_E c\ssp=\sp\{\ssp(\sp t\ssp,x\sp):\all z\,x=z\equivv c$ is
differentiable to $z$ in $E$ at $t\,\}\,$,}\vskip.3mm

\noin when $E$ is as above and \œ$\,c\in(\sp\vecs E\ssp)^{\,\roman{dom}\,c}$
with \œ$\dom c\inc\Re\sp\,$, otherwise putting\linebreak $\cder_E c=\Univ\,$,
and say that $\ssp c$ is {\it differentiable in} $\,E\ssp$ if{}f $\,
                  \dom c\inc\dom\cder_E c\not=\Univ\,$.

Let \œ$I=[\,0\,,1\,]$ and \œ$\sp Q=I^{\sp\times 2.}$ and \œ$\,\Cal U\sp =
\tau_{_{\Re}}\ttimes\tau_{_{\Re}}\spp\lei\sp Q\ssp$. Any \œ$\sp
c\in(\sp\vecs E\ssp)^{\,I}$ we call a {\it standard differentiable curve\ssp}
in $E$ in case also \œ${\,]}\,0\,,1\,{[\,}\inc\dom\cder_E c\not=\Univ\ssp$
holds and\linebreak $c$ is continuous \œ$\tau_{_{\Re}}\to\taurd E$ at the
points $0$ and $1\ssp$. By a {\it standard family of continuously
differentiable curves\ssp} in $\ssp E\sp$ we mean any \œ$\ssp
\Gamma\in\big((\sp\vecs E\sp)^{\,I\sp}\big){\sp^I}$ such that $\ssp\Gamma\sp$
is\linebreak continuous \œ$\ssp\tau_{_{\Re}}\to\taurd E\expnota^I]_{ti}$ and
such that for \œ$\ssp g = \sp\seq{\,
\cder_E(\sp\Gamma\sp\fvalue s\sp)\fvalue t:\smb S=(\sp s\ssp,t\sp)\in Q\,}\ssp
$ we also have that $\sp g\sp$ is continuous $\ssp\Cal U\to\taurd E$ with
$\sp Q\inc\dom g\,$.                                       \end{definitions}

\begin{proposition}[mean value theorem]\label{mvT}

Let $\sp U\sn$ be a closed and convex set in $E$ with \œ$\ssp E \in
\LCS(\bosy R\sp)\,$. Let $\,c\ssp$ be a standard differentiable curve in $\ssp
E\sp$ with \œ$\,c\fvalue 0=\bnull E$ $\in U$ and $\,\sp
\cder_E c\sp\image{]}\,0\,,1\,{[}\inc U\spp$. Then $\,c\fvalue 1\in U\sp$
holds.                                                     \end{proposition}

\begin{proof} To proceed by reductio ad absurdum, let \œ$c\fvalue 1\not\in U\spp
$. By Hahn\,--\,Banach, there is \œ$\ell\in\Cal L\ssp(E\ssp,\biit R\sp)$ with
\œ$\ell\ssp\image\snn U\inc{]}\minusinfty\,,1\,{[}$ and \œ$
\ell\ssp\fvalue(\sp c\fvalue 1)=1\ssp$. The classical mean value theorem when
applied to the function \œ$\sp\ell\circ c$ gives some $t$ with \œ$0<t<1$ and $
1 =\ell\ssp\fvalue(\sp c\fvalue 1)=\cder_{\bosy R}(\sp\ell\circ c\sp)\fvalue t
  = \ell\ssp\fvalue(\cder_E c\fvalue t\sp)\in\ell\ssp\image\snn U
\inc {]}\minusinfty\,,1\,{[}\ssp\,$, \,a contradiction.          \end{proof}

\begin{proposition}\label{diff under int}

Let \œ$\ssp I=[\,0\,,1\,]\,,$ and let \œ$\ssp E\in\LCS(\bosy R\sp)$ be
sequentially complete. Let $\,\Gamma\sp$ be a standard family of continuously
differentiable curves in $\ssp E\,,$ and let \œ$\,c\ssp = \sp
\Seq{\,E\text{\,-\sp}\int_{\,0}^{\,1}
          \Gamma\sp\fvalue s\fvalue t\ssp\,\roman d\,s:t\in I\sp\,}\,$. Then $\,
c\ssp$ is a standard differentiable curve in the space $\ssp E\sp$ with also $\,
\cder_E c\ssp = \sp\Seq{\,E\text{\,-\sp}\int_{\,0}^{\,1}
 \cder_E(\sp\Gamma\sp\fvalue s\sp)\fvalue t\ssp\,\roman d\,s:t\in I\sp\,}\,$.

  \end{proposition}

\begin{proof} Let $\sp\Cal U\sp,\sp g$ be as in
Definitions \ref{defi curve dif:ty} above, and arbitrarily fix \œ$t\ar 0 \in I\sp
$, and\linebreak put \œ$ x = E\text{\,-\sp}\int_{\,0}^{\,1}
          \cder_E(\sp\Gamma\sp\fvalue s\sp)\fvalue t\ar 0\,\roman d\,s\ssp$.
Then it suffices to show that $\sp c\sp$ is differentiable to $x$ in $E$ at $
t\ar 0\,$. Writing $\,\roman F\,s\,s\ar 1\ssp t\sp = \sp
t^{\sp\mminus 1}\sp(\sp\Gamma\sp\fvalue s\fvalue(\sp t\ar 0+s\ar 1\ssp t\sp)
                     -  \Gamma\sp\fvalue s\fvalue t\ar 0)
              - s\ar 1\ssp\cder_E(\sp\Gamma\sp\fvalue s\sp)\fvalue t\ar 0\,$,\linebreak
and arbitrarily fixing a closed convex \œ$\sp U\sn\in\ymp E\,$, using
Proposition \ref{conti is int} and Lemma \ref{int-mvT} above, one quickly
checks that it further suffices to show (e) existence of $\ssp\delta\in\Rep$
such that for all $t\in\Re$ with $t\ar 0+t\in I$ and $0<|\ssp t\ssp|<\delta$
we have $\,\roman F\,s\,1\,t\in U\spp$.

To get (e) above, we first note that using $\tau_{_{\Re\,}}$--\,compactness of
$\sp I$ and continuity \œ$\sp\Cal U\to\taurd E$ of $\sp g\,$, one deduces
existence of some \œ$\sp\delta\in\Rep$ such that for all \œ$t\in\Re$\linebreak
 with $t\ar 0+t\in I$ and $|\ssp t\ssp|<\delta$ we have\vskip.5mm\noin\ %
(x) $\mhyppy 6  \cder_E(\sp\Gamma\sp\fvalue s\sp)\fvalue(\sp t\ar 0+t\sp)
              - \cder_E(\sp\Gamma\sp\fvalue s\sp)\fvalue t\ar 0
= g\fvalue(\sp s\ssp,t\ar 0+t\sp) - g\fvalue(\sp s\ssp,t\ar 0) \in U\spp$. \vskip.5mm

\noin Now arbitrarily fixing $t$ as required in (e) we get \œ$\,
\roman F\,s\,1\,t\in U\sp$ by Proposition \ref{mvT} and (x) when considering
the curve $\,\seq{\,\sp\roman F\,s\,s\ar 1\ssp t:s\ar 1\snn\in I\sp\,}\,$.
  \end{proof}


  \binsubsubhead D{Seip's higher order differentiability classes}\label{subsec D}

In a rather non-uniform manner, Seip scattered his definitions of his first,
higher, and infinite order continuous differentiabilities in
\cite[pp.\ 60, 75, 85]{Se72}\ssp. In Definitions \ref{def Seip d:bilities}
below, we give our own reformulation of these by constructing the
differentiabi- lity classes $\ssp\CalDSeip^k$ for $\sp k\in\infty\ssp\yplus$
at a (multiple) stroke. In Theorem \ref{our equivv Seip} these turn out to
exactly capture Seip's concepts. Theorem \ref{Seip = BGN} shows that these are
obtained as particular cases of the general BGN\,--\,construction.

\begin{definitions}\label{def Seip d:bilities}

For all classes $\sp\tilde f\ssp,\sp k\ssp$ first generally letting\vskip.5mm

\noin (1) \ $\varSei^{}Ï\tilde f=\bigcap\ssp\{\ssp(E\skcap E\ssp,F\sp,g\sp):
\exi f\,\tilde f=(E\ssp,F\sp,f\sp)\in
\roman{scgVS}\ssp(\biit R\ssp)^{\sp\times2.}\ssn\times\snn\Univ$ \par \hyppy{6mm}
and $f\in F^{\sp/\sp E}$ and $\dom\sn f\in\taurd E$ and $
     g = \{\ssp(\sp x\ssp,u\ssp,y\sp):x\in\dom\sn f$ and \par $\mhyppy{6}
u\in\vecs E$ and $y\in\vecs F$ and $\all{\sp V\sn\in\Nbh(\sp y\ssp,
\taurd F\ssp)}\,\exi{\delta\in\Rep}\,\all{t\in\Re}$ \par $\null\hfill
0<|\ssp t\ssp|<\delta\imply t^{\sp\mminus 1}\sp(\sp
f\sp\fvalue(\sp x+t\,u\sp)- f\sp\fvalue x\sp)\in V\ssp\,\}\sp\}\,$,\KP{10}\vskip.3mm

\noin (2) \ $\varSei^kÏ\tilde f = \bigcap\ssp\{\ssp(\ssp
\vsprod_cgVS_\Re(\sp k\ssp\yplus\ssn\times\sn\{\sp E\ssp\}\sp)\ssp,F\sp,\sp
g\ssp):k:k\in\No$ and \par $\mhyppy{15}
E\ssp,F\in\roman{scgVS}\ssp(\biit R\ssp)$ and $\sp\exi{f\sp,\bosy f}\,
f\in F^{\sp/\sp E}$ and $\tilde f=(E\ssp,F\sp,f\sp)$ and \par $\mhyppy{11}
(\sp\emptyset\,,f\circ\seq{\,\seq{\ssp x\ssp}:x=x\,}\inve\spp)\ssp,
(\sp k\ssp,g\sp)\in\bosy f\in\Univ\,^{k\ssp+\sp 1.}$ and $\ssp\all{i\in k}$ \par $\mhyppy{7}
\bosy f\sp\fvalue i\ssp\yplus =
\{\ssp(\spp\seq{\sp x\sp}\concc\bosy u\concc\seq{\sp u\sp}\ssp,y\sp):\all h\,
h=\{\ssp(\sp x\yr 1\spp,y\yr 1):(\spp\seq{\ssp x\yr 1\sp}\concc\bosy u\ssp,
y\yr 1)\in\bosy f\sp\fvalue i\,\}$ \par $\null\hfill
\imply(\sp x\ssp,u\ssp,y\sp)\in\taurd\varSei^{}Ï(E\ssp,F\sp,h\ssp)\ssp\}\sp\}\,$,\KP{10}\vskip.3mm

\noin (3) \ $\CalCSe0 = \{\ssp(E\ssp,F\sp,f\sp):f\in F^{\ssp/\sp E}$ and $\sp
f\sp\inve\images\taurd F\inc\taurd E\sp\,\}
\,|\,\sp\roman{scgVS}\ssp(\biit R\ssp)^{\sp\times2.}\,$,\vskip.7mm

\noin any $\sp\tilde f$ with \œ$\dom\taurd\tilde f\times\vecs\ssigrd\tilde f
\inc\dom\taurd\varSe\tilde f\not=\Univ\sp$ we say to be {\it directionally
Seip \biggerlineskip5 differentiable\ssp}. Any $\sp\tilde f$ with \œ$
\tilde f\sp,\varSe\tilde f \in \CalCSe0$ and \œ$
\dom\taurd\tilde f\times\vecs\ssigrd\tilde f\inc\dom\taurd\varSe\tilde f$ we
say to be {\it Seip differentiable\ssp}. We say that $\tilde f$ is {\it order
$k$ simply Seip\,--\,differenti- \biggerlineskip3 able\ssp} if and only if we
have \œ$k\in\infty\ssp\yplus$ and for all \œ$i\in k+1\adot$ it holds that $
\varSei^iÏ\tilde f\in\CalCSe0$ and $\{\,\biit x:\biit x\fvalue\emptyset \in
\dom\taurd\sp\tilde f\ssp\,\}\cap(\sp\vecs\sp\ssigrd\sp\tilde f\,
)^{\,i\ssp+\sp 1.}\inc\dom\taurd\varSei^iÏ\tilde f\ssp$. We further let\vskip.5mm

\noin (4) \ $\DerSe\tilde f = \bigcap\ssp\{\,\tilde f\ssp':\exi{E\ssp,F}\,
\tilde f$ is Seip differentiable and\nopagebreak\par\nopagebreak$\null\hfill
(E\ssp,F\sp)=\sigrd\sp\tilde f$ and $\tilde f\ssp'=(E\ssp,\Lincg(E\ssp,F\sp)\ssp
, (\sp\taurd\varSe\tilde f\,)\snn\yvee\spp)\ssp\}\,$,\KP{10}\vskip.3mm

\noin (5) \ $\CalDSeip^k=\{\,\tilde f:k\in\infty\ssp\yplus$ and $\exi{\biit f}\,
(\sp\emptyset\,,\tilde f\,)\in\biit f\in\CalCSe0^{k\ssp+\sp 1.}$ and \par $\null\hfill
\all{i\in k}\,\bosy f\sp\fvalue i\sp$ is Seip differentiable and $
\bosy f\sp\fvalue i\ssp\yplus =
\varSe\snn(\sp\bosy f\sp\fvalue i\sp)\ssp\}\,$,\KP{10}\vskip.5mm

\noin and also \œ$\cgbarDelta\tilde f=\CalCSe0\vbarDelta_R\tilde f\sp$, and
suggest that $\varSei^{}Ï\tilde f$ and $\varSei^kÏ\tilde f$ may be referred to
by the phrases the Seip{\it\,--\,variation\ssp} and the order $\sp k\ssp$
Seip\,--\,variation of $\sp\tilde f\sp$. The class $\DerSe\tilde f\sp$ may be
called the Seip {\it derivative\ssp} of $\sp\tilde f\sp$.  \end{definitions}

Note that any directionally Seip differentiable class $\sp\tilde f$
necessarily is a real vector map since by our definitions from \œ$
\dom\taurd\varSe\tilde f\not=\Univ\sp$ it follows that there are $\sp
E\ssp,F\sp,f$ \Biggerlineskip1 with $\sp
E\ssp,F\in\roman{scgVS}\ssp(\biit R\ssp)\sp$ and $\ssp
f\in F^{\sp/\sp E}\sp$ and $\,\tilde f=(E\ssp,F\sp,f\sp)\,$.

An exercise in set theory and logic shows that we have either $\ssp
\varSei^kÏ\tilde f=\Univ\,$ or \vskip-3.5mm {\leftskip4mm\null\hskip-8mm%
(u) \ \œ$k\in\No$ and there are unique \œ$\bosy g\in\Univ\,^{k\ssp+\sp 1.}$
    and \œ$E\ssp,F\in\roman{scgVS}\ssp(\biit R\ssp)$ and \œ$\sp f \in
    F^{\sp/\sp E}$ with \œ$\tilde f=(E\ssp,F\sp,f\sp)$ and \œ$
\bosy g\fvalue\emptyset=\{\ssp(\spp\seq{\ssp x\ssp}\ssp,y\sp):(\sp x\ssp,y\sp)
\in f\sp\,\}$ and $\,\bosy g\fvalue i$ a function for\linebreak all \œ$ i \in
k\ssp\yplus$ such that in the case \œ$i\not=k$ we have \œ$
\{\,x:\seq{\ssp x\ssp}\concc\bosy x\in\dom\snn(\ssp\bosy g\fvalue i\sp)\ssp\}
\in\taurd E$\linebreak for every \œ$\bosy x\,$, and for all $\bosy z\ssp,y$ we
have \œ$(\sp\bosy z\ssp,y\sp)\in\bosy g\fvalue i\ssp\yplus$ if and only if
there are $x\ssp,u\in\vecs E$ and $\bosy x\in(\sp\vecs E\ssp)^{\,i}$ with $
\bosy z=\seq{\ssp x\ssp}\concc\bosy x\concc\seq{\ssp u\ssp}$ and \vskip0mm\centerline{$
y = \taurd F\text{\,-\,}\underset{t\sp\to\sp 0}\lim\,t^{\sp\mminus 1}\sp(\ssp
     \bosy g\fvalue i\sp\fvalue(\spp\seq{\,x+t\,u\,}\concc\spp\bosy x\sp)
    - \bosy g\fvalue i\sp\fvalue(\spp\seq{\ssp x\ssp}\concc\bosy x\sp))
       \not= \Univ\,$.}\vskip.3mm
                                   \par} 

\noin Here we understand that a limit equals $\Univ$ if it \q{does not exist}.
In the case (u) we also have $ \varSei^kÏ\tilde f =
(E\expnota^k\ssp+\sp 1.]_{kvs}\spp,\sp F\sp,\ssp\bosy g\fvalue k\ssp)\,$,
\,for $\,E\expnota^I]_{kvs} = \prod\sn\RHB{.5}{_{_{\roman{cgVS\sp\,
\sixroman{dom}^2\sp\sixmath\tauu\sixmath\sigmaa_{\sn r\sn d}}\sp
      \sixmath E\,}}}(\sp I\timesn\{\sp E\ssp\}\sp)\,$.

\begin{proposition}\label{Seip0 is BGN class}

\hfill $\CalCSe0\ssp$ is a \,\erm{BGN}\,--\,class on $\,\scgVS(\bosy R\ssp)$
over $\ssp\bosy R\,$. \hfill\null                          \end{proposition}

\begin{proof} First, to see that $\CalCSe0$ is a productive class on $\sp
\scgVS(\bosy R\ssp)$ over $\sp\bosy R\,$, we must verify that for all \œ$
E\ssp,F\in\scgVS(\bosy R\ssp)$ there is $G$ with $\sp
\CalCSe0\text{\ssp-\,prod}\ssp\subtext{mcl}\sp(E\ssp,F\sp,G\sp)$ as specified
in Definitions \ref{BGN axioms} above. To verify this for \œ$G=E\skcap F\sp$,
first note that by Corollary \ref{Seip vs products} we have \œ$ G \in
\scgVS(\bosy R\ssp)\,$, and then that for \œ$(\spp H\sp,E\ssp,f\sp)\ssp,
(\spp H\sp,F\sp,g\sp)\in\CalCSe0$ we have $
(\spp H\sp,G\sp,[\,\sp f\sp,\spp g\sp\,]\ssp)\in\CalCSe0$ by
Lemmas \ref{open and clos comp gen} and \ref{basic comp gen property}\ssp.

For the verification of the postulates (1)\œ$\ssp,\ldots\,$(6) in
Definitions \ref{BGN axioms}\ssp, first note that by
Lemma \ref{for product family} we get (6) directly from the definitions since
\œ$\ssp \CalCSe0\vProd\KN1\fvalue(\bosy R\,,E\sp) = \bosy R\skcap E $ and
trivially \œ$\ssp \tau_{_{I\!\!R}}\snn\ttimes\taurd E \inc
\tau_{_{I\!\!R}}\snn\ktimes\taurd E=\taurd(\bosy R\skcap E\ssp)$ when \œ$
E\in\scgVS(\bosy R\ssp)\,$. We have (1)\œ$\ssp,\ldots\,$(5) as trivialities or
well-known facts given by straighforward verifications, noting e.g.\ that
(\ref{locality}) just means that a function is continuous with open domain if
it can be expressed as a union of such, and that for (\ref{compos rule}) under
the premise we have $\,(\ssp g\sp\circ f\ssp)\inve\images\taurd G \inc
f\ssp\inve\images(\ssp g\sp\inve\images\taurd G\ssp) \inc
f\ssp\inve\images\taurd F\inc\taurd E\,$.                        \end{proof}

In \cite[Theorem 3.8, pp.\ 82\,--\,83]{Se79}\ssp, Seip incidentally
                                          established the following

\begin{lemma}\label{T3.8 of Se79}

Let \œ$\ssp\tilde f=(E\ssp,F\sp,f\sp)\,$. If $\ssp\tilde f$ is directionally
Seip\,--\,differentiable and $\,\varSe\tilde f$ is continuous, there is $\ssp
g$ with the property that \œ$\ssp(E\spp\skcap E\spp\skcap\biit R\,,F\sp,g\sp)
\in\CalCSe0\,,$ and also such that for all $\ssp x\ssp,u\ssp,t\ssp,y\ar 1$
with $\ssp y\ar 1\in\vecs F\ssp$ it holds that\vskip.3mm\centerline{$
  f\sp\fvalue(\sp x+t\,u\sp)=f\sp\fvalue x+t\,y\ar 1\not=\sp\Univ\,\equivv\,
   \exi y\,(\sp x\ssp,u\ssp,t\ssp,y\sp)\in g$ and $\,
   [\ t\not=0\imply y=y\ar 1\sp\,]\,$.}                          \end{lemma}

Putting \ $\Cal D\yr 1\KN{1.8}\subtext{Seip\,72} = \{\ssp(E\ssp,F\sp,f\sp):
              f$ is {\it stetig differenzierbar\ssp}

$\null\hfill E\supset\dom\sn f\to F$ in the sense of
               \cite[Definition 4.2, p.\ 60]{Se72}$\ssp\,\}\,$, \,we have

\begin{corollary}\label{D^1 iff BGN^1 iff Seip^1(1972)}

For all $\,\tilde f$ and $\,\tilde g\sp$ it holds that\vskip.7mm

{\rm(a)} \ \ $\CalDSeip^{1.} = \sp
              \CalDBGN^{1.}(\sp\CalDSeip^{0.}\ssp,\spp\biit R\ssp) = \sp
              \Cal D\yr 1\KN{1.8}\subtext{Seip\,72}\,,$\vskip.4mm

{\rm(b)} \ \ $\tilde f$ is Seip\,--\,differentiable $\,
               \imply\,\DerSe\tilde f\,$ is a vector map$\,,$\vskip.4mm

{\rm(c)} \ \ $\tilde f$ is Seip\,--\,differentiable $\,\equivv\,
              \tilde f \in \CalDSeip^{1.}\,\equivv$\par$\null\hfill
              \tilde f$ is directionally Seip\,--\,differentiable and
                                  $\,\varSe\tilde f\in\CalCSe0\,,$ \KP{15} \vskip.4mm

{\rm(d)} \ \ $\tilde f$ and $\,\tilde g$ are Seip\,--\,differentiable and $\,
\tsigrd\sp\tilde f=\ssigrd\sp\tilde g$\par$\null\hfill
\imply\,\taurd\varSe\sn(\ssp\tilde g\Circ\snn\tilde f\,) =
\taurd\varSe\tilde g\circ[\,\sp\taurd\tilde f\circ\spp\roman{pr}\ar 1\ssp,\sp
                   \taurd\varSe\tilde f\ssp\,]\,$. \KP{15}   \end{corollary}

\begin{proof} The equality \œ$\ssp \CalDSeip^{1.} = \sp
\CalDBGN^{1.}(\sp\CalDSeip^{0.}\ssp,\spp\biit R\ssp)$ follows from
Lemma \ref{T3.8 of Se79} above by a nice exercise in logic. For (b)
assuming that \œ$\sp\tilde f=(E\ssp,F\sp,f\sp)\in\CalDSeip^{1.}\,$, it
suffices for \linebreak each fixed \œ$x\in\dom\sn f$ to establish linearity \œ$
\sigrd E\to\sigrd F$ of $\sp\taurd\varSe\tilde f\ssp(\sp x\ssp,\cdot\sp)\,$.
This follows from \cite[Proposition 2.2, p.\ 223]{BGN}\ssp.

To get \œ$\ssp\CalDBGN^{1.}(\sp\CalDSeip^{0.}\ssp,\spp\biit R\ssp) = \sp
\Cal D\yr 1\KN{1.8}\subtext{Seip\,72}\,$, we first note that for \œ$\tilde f =
(E\ssp,F\sp,f\sp)$ we have \œ$\tilde f\in\Cal D\yr 1\KN{1.8}\subtext{Seip\,72}
$ if and only if \œ$E\ssp,F\in\scgVS(\biit R\ssp)$ and for \œ$L=\Lincg(E\ssp,
F\ssp)$ and \œ$\ssp U=\dom\sn f$ and \œ$U\aar 3=\{\ssp(\sp x\ssp,u\ssp,t\sp):
x\ssp,x+t\,u\in U\sp\,\}\,$, we have \œ$U\in\taurd E$ and \œ$f \in
(\sp\vecs F\ssp)\,^U\sp$, and there are \œ$f\aar 1\in(\sp\vecs L\sp)\,^U\sp$
and \œ$h\in(\sp\vecs F\ssp)\,^{U_3}$ such that $f\aar 1$ is continuous \œ$
\taurd E\to\taurd L$ and $h$ is continuous \œ$\taurd(E\spp\skcap E\spp\skcap
\biit R\ssp)\to\taurd F$ with \œ$U\sn\times\sn(\sp\vecs E\ssp)\sn\times\sn
\{0\}\sn\times\sn\{\ssp\bnull F\}\inc h\,$,\linebreak and such that for $
(\sp x\ssp,u\ssp,t\sp)\in U\aar 3$ it holds that\vskip.3mm\centerline{$
f\sp\fvalue(\sp x+t\,u\sp) = f\sp\fvalue x + f\aar 1\KN{.7}\fvalue x\fvalue
 (\sp t\,u\sp) + t\,h\fvalue(\sp x\ssp,u\ssp,t\sp)\,$.}\vskip.3mm

Noting that by Proposition \ref{lin expo} continuity of $(E\ssp,L\ssp,f\aar 1)
$ is equivalent to that of \œ$\ssp(E\spp\skcap E\ssp,F\sp,f\aar 1\KN1\ywed\sp)\,
$, and recalling the above established property that \œ$\sp\tilde f \in
\CalDSeip^{1.}$ implies linearity of $(\sp\sigrd E\ssp,\sigrd F\sp,\taurd
\varSe\tilde f\ssp(\sp x\ssp,\cdot\sp))$ for all \œ$x\in U\spp$, it is a
simple exercise in logic (to the reader) to verify that $\ssp
\CalDBGN^{1.}(\sp\CalDSeip^{0.}\ssp,\spp\biit R\ssp) = \sp
\Cal D\yr 1\KN{1.8}\subtext{Seip\,72}\,$.

For (c) assuming that $\tilde f$ is directionally Seip\,--\,differentiable
with \œ$\sp\varSe\tilde f\in\CalCSe0\,$, we only have to prove continuity of $\sp
\tilde f\sp$. Taking \œ$t=1$ in Lemma \ref{T3.8 of Se79} above, we have \œ$\sp
f\sp\fvalue(\sp x+u\sp)=f\sp\fvalue x+g\fvalue(\sp x\ssp,u\ssp,1\sp)$ for $u$
close to $\bnull E\,$, whence the assertion. For (d) we refer the reader to
see \cite[Proposition 3.1, p.\ 225]{BGN}\ssp.                    \end{proof}

Property (d) in Corollary \ref{D^1 iff BGN^1 iff Seip^1(1972)} may be called
the {\it first order chain rule\ssp}. The differentiability classes
$\sp\CalDSeip^k$ for different $k$ are related to one another via the
recursion rules given in the following

\begin{proposition}\label{D^k recursion}

For all $\,k\ssp,\tilde f$ it holds that \œ$\,\tilde f \in
\CalDSeip^{k\ssp+\sp 1.}\sp\equivv\sp\text{(a)}\sp\equivv\sp\text{(b)}\,,$
when the con- ditions are as given below. It also holds that $\,
\CalDSeip^\infty=\{\,\tilde f:\all{k\in\No}\,\tilde f\in\CalDSeip^k\ssp\}\,$.\vskip.7mm

{\rm(a)} \ \ $\tilde f\in\CalDSeip^{1.}$ and \œ$\ssp
               \varSei^{}Ï\tilde f \in\CalDSeip^k\,,$\vskip.3mm

{\rm(b)} \ \ $\tilde f$ is directionally Seip\,--\,differentiable and $\ssp
               \varSei^{}Ï\tilde f\in\CalDSeip^k\,$.       \end{proposition}

\begin{proof} One obtains \œ$\sp\tilde f \in\CalDSeip^{k\ssp+\sp 1.}\equivv\sp{}
$(a) similarly as in the proof of Proposition \ref{BGN recurs} above. For
(a)\œ${}\equivv{}$(b)\ssp, it suffices to observe that \œ$\sp\tilde f \in
\CalDSeip^{1.}$ follows from directional differentiability and continuity of
the variation by Corollary \ref{D^1 iff BGN^1 iff Seip^1(1972)}\ssp(c) above.

For \œ$\sp\Cal D\ar 0=\{\,\tilde f:\all{k\in\No}\,\tilde f\in\CalDSeip^k\ssp\}\,
$, to establish \œ$\sp\CalDSeip^\infty=\Cal D\ar 0\,$, note that we trivially
have \œ$\CalDSeip^\infty\inc\Cal D\ar 0\,$. To get the converse, arbitrarily
fixing \œ$\sp\tilde f\in\Cal D\ar 0\,$, and considering the class\vskip.4mm

$\Gamma=\{\,\bosy f:\exi{k\in\infty\ssp\yplus}\,\all i\,(\sp\emptyset\,,
\tilde f\sp)\in\bosy f\in\CalCSe0^{k\ssp+\sp 1.}$ and

$\null\hfill[\ i\in k\imply\bosy f\ssp\fvalue i$ is Seip differentiable and $
\bosy f\ssp\fvalue(\sp i\ssp\yplus) =
\varSe(\sp\bosy f\ssp\fvalue i\sp)\ ]\sp\,\}\,$,\KP6\vskip.4mm

\noin we first see for $\bosy f\ssp,\ssp\bosy g\in\Gamma$ that either \œ$\ssp
\bosy f\inc\bosy g\,$ or \œ$\,\bosy g\inc\bosy f\ssp$, \,and consequently that
$\ssp\bigcup\,\Gamma$ is a function. In view of \œ$\sp\tilde f\in\Cal D\ar 0\,
$, this further gives \œ$\ssp\bigcup\,\Gamma\in\Gamma$ with $\ssp
\dom\sn\bigcup\,\Gamma=\No\ssp$, whence $\tilde f\in\CalDSeip^\infty$ follows
directly by definition.                                          \end{proof}

Some further basic properties of the classes $\sp\CalDSeip^k$ are given in the
following

\begin{proposition}\label{basic D^k properties}

For all $\,E\ssp,F\sp,G\sp,H\sp,P\sp,f\sp,g\ssp,k\,$ it holds that\vskip.8mm

\begin{enumerate}

\item\label{pr1a} \ \ $\{\ssp(E\ssp,F\sp,\sp U\snn\times\sn\{\ssp y\sp\}\sp) :
  U\snn\in\taurd E\text{ and }\sp y\in\vecs F\sp\,\}\,|\,\sp
  \roman{scgVS}\ssp(\biit R\ssp)^{\sp\times 2.}\inc\CalDSeip^\infty\,,$\vskip.4mm

\item\label{pr1b} \ \ $\{\ssp(E\ssp,F\sp,\sp\ell\,) :
  \ell\in\vecs\Lincg(E\ssp,F\ssp)\ssp\}\,|\,\sp
  \roman{scgVS}\ssp(\biit R\ssp)^{\sp\times 2.}\inc\CalDSeip^\infty\,,$\vskip.4mm

\item\label{pr1x} \ \ $\RHB{1.5}{\fiveroman{bi}\sn}\Cal L\sp\subtext{scgVS}\spp
                            (\biit R\ssp)\inc\CalDSeip^\infty\,,$\vskip.4mm

\item\label{pr1c} \ \ $\,(E\ssp,F\sp,f\sp)\ssp,(E\ssp,G\sp,g\sp)\in\CalDSeip^k
  \imply(E\ssp,F\skcap G\sp,\spp[\,\sp f\sp,\spp g\sp\,]\ssp)\in\CalDSeip^k\,,$\vskip.4mm

\item\label{pr1d} \ \ $\,(E\ssp,F\sp,f\sp)\ssp,(\sp F\spp,G\sp,g\sp) \in
  \CalDSeip^k\imply(E\ssp,G\sp,\spp g\spp\circ f\sp)\in\CalDSeip^k\,,$ \hfill
                                               {\rm(chain rule)}\vskip.4mm

\item\label{pr1e} \ \ $\,(E\ssp,F\sp,f\sp)\ssp,(H\ssp,G\sp,g\sp)\in\CalDSeip^k
  \imply(E\spp\skcap H\sp,F\skcap G\sp,f\risti2 g\sp)\in\CalDSeip^k\,,$\vskip.4mm

\item\label{pr1l} \ \ $\,P=(E\ssp,F\sp)\in\scgVS(\bosy R\ssp)^{\sp\times 2.}$
  and $\sp f$ is a function and \par\hyppy{-1.1mm}
  $[\sp\ \all{\smb Z}\,\exi h\,\smb Z\in f\imply\smb Z\in h\inc f$ and $\,
  (\spp P\spp,h\sp)\in\CalDSeip^k\ssp\,] \imply
  (\spp P\spp,f\sp)\in\CalDSeip^k\,$.     \end{enumerate}  \end{proposition}

\begin{proof} We get (\ref{pr1a}) by a trivial induction using
Proposition \ref{D^k recursion} above when for fixed \œ$\sp\tilde f =
(E\ssp,F\sp,\sp U\snn\times\sn\{\ssp y\sp\}\sp)$ we note that \œ$
\varSe\tilde f=(E\spp\skcap E\ssp,F\sp,\sp U\snn\times\sn(\sp\vecs E\ssp)\sn
\times\sn\{\ssp\bnull F\}\sp)\,$. We get (\ref{pr1b}) similarly, noting that $
\varSe\sn(E\ssp,F\sp,\sp\ell\,)=(E\spp\skcap E\ssp,F\sp,\sp\ell\spp\circ
\roman{pr}\ar 2\ssp|\,(\sp\vecs E\ssp)^{\sp\times2.})\,$. In or- der to
establish (\ref{pr1x})\,, note that for $\tilde b=(E\ar 2\spp,F\sp,b\ssp)\in
\RHB{1.5}{\fiveroman{bi}\sn}\Cal L\sp\subtext{scgVS}\spp(\biit R\ssp)$ we have\vskip.3mm\centerline{$
\taurd\varSe\tilde b=\seq{\,b\fvalue(\sp x\ssp,v\sp)+b\fvalue(\sp u\ssp,y\sp):
\smb W=(\sp x\ssp,y\,;u\ssp,v\sp)\in(\sp\vecs E\ar 2)^{\sp\times2.}\ssp}\,$,}\vskip.3mm

\noin and hence that also $\varSe\tilde b\in
\RHB{1.5}{\fiveroman{bi}\sn}\Cal L\sp\subtext{scgVS}\spp(\biit R\ssp)\,$.

For (\ref{pr1c}) the case \œ$k=\emptyset$ follows from
Proposition \ref{Seip0 is BGN class} above. For the inductive step, letting \œ$
f\aar 1 = \taurd\varSe\sn(E\ssp,F\sp,f\sp)$ and \œ$g\ar 1 =
\taurd\varSe\sn(E\ssp,G\sp,g\sp)\,$, note that we have\linebreak $
\varSe\sn(E\ssp,F\skcap G\sp,\spp[\,\sp f\sp,\spp g\sp\,]\ssp) = $ $
(E\spp\skcap E\ssp,F\skcap G\sp,\spp[\,\sp f\aar 1\ssp,\spp g\ar 1\,]\ssp)\,$.

For (\ref{pr1d}) the case \œ$k=\emptyset$ is got from
Proposition \ref{Seip0 is BGN class}\ssp. For \œ$\ssp k\in\N\sp$, one proceeds
by induction, noting that with \œ$f\aar 1 = \taurd\varSe\sn(E\ssp,F\sp,f\sp)$
and \œ$\ssp g\ar 1 = \taurd\varSe\sn(\sp F\spp,G\sp,g\sp)\,$, we have \œ$
\varSe\sn(E\ssp,G\sp,\spp g\sp\circ f\sp)$ \œ$=(E\spp\skcap E\,,\sp G\ssp,\sp
g\ar 1\snn\circ[\,\sp f\circ\spp\roman{pr}\ar 1\ssp,\sp f\aar 1\,]\ssp)\sp$ by
the first order chain rule. For the inductive step, assuming that we have the
result with \œ$\sp k\in\No$ and that \œ$
(E\ssp,F\sp,f\sp)\ssp,(\sp F\spp,G\sp,g\sp)\in\CalDSeip^{k\ssp+\sp 1.}\sp$,
first note that \œ$(E\skcap E\,,F\skcap F\sp,
[\,\sp f\circ\spp\roman{pr}\ar 1\ssp,\sp f\aar 1\,]\ssp)\in\CalDSeip^k$ by
(\ref{pr1b}) and (\ref{pr1c})\ssp, and then use the inductive assumption again.

Noting that $f\risti2 g = [\,\sp f\snn\circ\sp\roman{pr}\ar 1\ssp,\sp
g\circ\roman{pr}\ar 2\,]\,$, we get (\ref{pr1e}) from (\ref{pr1b}) and
(\ref{pr1c}) and (\ref{pr1d})\ssp.

Finally, for the inductive proof of (\ref{pr1l}) the case \œ$k=\emptyset$ has
been \q{explained\sp} in the proof of Proposition \ref{Seip0 is BGN class}
above. Now assuming that we have (\ref{pr1l}) for a fixed \œ$k\in\No$ and for
all $P\sp,E\ssp,F\sp,f\sp$, and that also \œ$ P = (E\ssp,F\ssp) \in
\scgVS(\bosy R\ssp)^{\sp\times 2.}$ and $\sp f$ is a function such that
\,(\œ$*)\ \ [\sp\ \all{\smb Z}\,\exi h\,\smb Z\in f\imply\smb Z\in h\inc f\text{
and }\sp(\spp P\spp,h\sp)\in\CalDSeip^{k\ssp+\sp 1.}\ ]\,$ holds, to establish
\œ$(\spp P\spp,f\sp)\in\CalDSeip^{k\ssp+\sp 1.}\sp$, by
Proposition \ref{D^k recursion}\ssp(b) it suffices that $(\spp P\spp,f\sp)$ is
directionally Seip\,--\,differentiable and that $\ssp
             \varSei^{}Ï(\spp P\spp,f\sp)\in\CalDSeip^k\,$.

Having $(\spp P\spp,h\sp)$ directionally Seip\,--\,differentiable when \œ$
(\spp P\spp,h\sp)\in\CalDSeip^{k\ssp+\sp 1.}\inc\CalDSeip^{1.}\,$, inspection
of the definitions shows that from ($*$) we get that $(\spp P\spp,f\sp)$ is
directionally Seip\,--\,differentiable. To get \œ$\ssp
\varSei^{}Ï(\spp P\spp,f\sp)\in\CalDSeip^k\,$, putting \œ$G=E\skcap E$ and \œ$
Q=(\spp G\sp,F\ssp)$ and \œ$f\aar 1=\taurd\varSei^{}Ï(\spp P\spp,f\sp)\,$, we
apply the inductive assumption with $Q\ssp,G\sp,F\sp,f\aar 1$ in place of $
P\sp,E\ssp,F\sp,f\sp$. Indeed, if \œ$\smb W\in f\aar 1\ssp$, there are $x\ssp,
v$ with \œ$\smb W=(\sp x\ssp,v\ssp,f\aar 1\KN1\fvalue(\sp x\ssp,v\sp))\,$. For
\œ$\smb Z=(\sp x\ssp,f\sp\fvalue x\sp)$ then \œ$\smb Z\in f$ holds whence by
($*$) there is $h$ with \œ$\smb Z\in h\inc f$ and \œ$(\spp P\spp,h\sp) \in
\CalDSeip^{k\ssp+\sp 1.}\sp$. For \œ$h\ar 1=\taurd\varSei^{}Ï(\spp P\spp,h\sp)
$ it follows (by an exercise to the reader) that \œ$\smb W\in h\ar 1 \inc
f\aar 1$ with \œ$(\sp Q\ssp,h\ar 1)\in\CalDSeip^k\,$. Consequently, the
premise in the inductive assumption holds, hence also the conclusion, that is
$\ssp\varSei^{}Ï(\spp P\spp,f\sp)=(\sp Q\ssp,f\aar 1)\in\CalDSeip^k\,$.
  \end{proof}

\begin{proposition}\label{Seip iff simply Seip}

For all $\ssp k\sp$ it holds that\vskip.3mm\centerline{$
\CalDSeip^k=\{\,\tilde f:\tilde f\text{ is order $k$ simply
                Seip\,--\,differentiable}\ssp\,\}\,$.}     \end{proposition}

\begin{proof} Let \œ$\sp\Cal D\sbi k=\{\,\tilde f:\tilde f\text{ is order $k$
simply Seip\,--\,differentiable}\ssp\,\}\,$, and let $(k)\subtext A$ mean that
\œ$\sp\CalDSeip^k=\Cal D\sbi k$ holds. From
Definitions \ref{def Seip d:bilities} above, one first checks that we have \œ$\sp
\Cal D\sbi k=\emptyset$ for $\sp k\not\in\infty\ssp\yplus\sp$, and that \œ$\sp
\Cal D\sbi\infty=\{\,\tilde f:\all{k\in\No}\,\tilde f\in\Cal D\sbi k\ssp\}\,$.
In view of Proposition \ref{D^k recursion} above, it hence suffices to
establish $\ssp\all{k\in\No}\,(k)\subtext A\,$.

Given \œ$ \tilde f = (E\ssp,F\sp,f\sp) \in
\roman{scgVS}\ssp(\biit R\ssp)^{\sp\times2.}\ssn\times\snn\Univ$ with \œ$
f\in F^{\sp/\sp E}\sp$, and considering the linear homeomorphisms \œ$
\Iota\ar 1\sn:E\to E\expnota^1.]_{kvs}$ and \œ$\Iota\ar 2 \sn : E\skcap E \to
E\expnota^2.]_{kvs}$ given by \œ$x\mapsto\seq{\ssp x\ssp}$ and \œ$
(\sp x\ssp,y\sp)\mapsto\seq{\,x\ssp,y\,}\,$, respectively, from the result
around (u) after Definitions \ref{def Seip d:bilities}\linebreak we see that \œ$\,
f=\taurd\varSei^0.Ï\tilde f\circ\Iota\ar 1$ and \œ$\,[\
\taurd\varSei^{}Ï\tilde f=\taurd\varSei^1.Ï\tilde f\circ\Iota\ar 2$ or \œ$
\ssp[\ \dom\sn f\not\in\taurd E$ and\linebreak \œ$\varSei^{}Ï\tilde f = 
\varSei^1.Ï\tilde f=\Univ\ ]\ ]\,$. Using this, by
Definitions \ref{def Seip d:bilities} we get $\ssp
               \all{k\in 2\adot}\,(k)\subtext A\,$.

We next assume that \œ$k\in\N$ with $(k)\subtext A\,$, and proceed to get $
(\sp k\ssp\yplus)\subtext A$ as follows. First letting \œ$\tilde f =
(E\ssp,F\sp,f\sp)\in\CalDSeip^{k\ssp+\sp 1.}\,$, then \œ$\tilde f \in
\CalDSeip^k\,$, and by Proposition \ref{D^k recursion}\ssp(a) also \œ$\tilde f
\in\CalDSeip^{1.}$ and \œ$\ssp\varSei^{}Ï\tilde f \in\CalDSeip^k\,$. Using $
(k)\subtext A\,$, we obtain \œ$\tilde f\ssp,\sp\varSei^{}Ï\tilde f \in
\Cal D\sbi k\,$. A \q{straightfor- ward} induction on $l\in k\ssp\yplus$ shows
that we have\vskip.3mm\centerline{$
\taurd\varSei^\sp l\ssp+\sp 1.Ï\tilde f\sp\fvalue\bosy x=
\taurd\varSei^\sp lÏ\varSei^{}Ï\tilde f\sp\fvalue
\big(\sp\seq{\ssp(\sp\bosy x\fvalue\emptyset\,,
                             \bosy x\fvalue 1\adot)\ssp}\concc
\seq{\,(\sp\bosy x\fvalue i\ssp\yplus\sn\yplus\spp,\bnull E):i\in l\ssp\,}
\sp\big)\,$,}\vskip.3mm

\noin for \œ$\bosy x\in(\sp\vecs E\ssp)\,^{l\ssp+\ssp 2.}$ with \œ$
\bosy x\fvalue\emptyset\in\dom\sn f\sp$. Taking here \œ$\ssp l=k\ssp$, we get
\œ$\varSei^k\ssp+\sp 1.Ï\tilde f\in\CalCSe0\,$, and see that \œ$
\{\,\biit x:\biit x\fvalue\emptyset\in\dom\sn f\ssp\,\}\cap(\sp\vecs
E\ssp)^{\,k\ssp+\ssp 2.}\inc\dom\taurd\varSei^k\ssp+\sp 1.Ï\tilde f\ssp$.
Consequently, noting $\tilde f\in\Cal D\sbi k\,$, we get
              $\sp\tilde f\in\Cal D\sbi{k\ssp+\sp 1.}\,$.

For the converse proof we first introduce the notation\vskip.3mm\centerline{$
\roman R\,x\,\bosy u\,i\,\bosy x=\seq{\,x\ssp,\bosy u\fvalue i\,}\concc(\sp
\bosy x\,|\,i\sp)\concc\big(\sp\bosy x\circ\seq{\,i+j+1\adot\sn:j\in\No\ssp}\big)\,$.}\vskip.3mm

\noin Next, using the linearity from the proof of
(b) in Corollary \ref{D^1 iff BGN^1 iff Seip^1(1972)} above, by induction on
\œ$l\in\No$ one (= the reader\ssp!\sp)\ proves ($*$) that for an order $\ssp
l\ssp\yplus$ simply Seip\,--\,differentiable map \œ$\tilde f=(E\ssp,F\sp,f\sp)$
and for any \œ$\bosy x\in(\sp\vecs E\ssp)\,^l$ and \œ$x\in\dom\sn f$ we have a
linear map $\seq{\,\taurd\varSei^\sp l\ssp+\sp 1.Ï\tilde f\sp\fvalue
(\spp\seq{\,x\ssp,u\,}\concc\spp\bosy x\ssp):u\in\vecs E\sp\,}
               : \sigrd E\to\sigrd F\sp$.\vskip.3mm

Now, assuming that \œ$\tilde f=(E\ssp,F\sp,f\sp)\in\Cal D\sbi{k\ssp+\sp 1.}\,
$, then \œ$\tilde f\in\Cal D\sbi{1.}$ whence by $(\sp 1\adot)\subtext A$ we
have \œ$\tilde f\in\CalDSeip^{1.}\,$. Consequently, by
Proposition \ref{D^k recursion}\ssp(a) and $(k)\subtext A$ to get \œ$\tilde f
\in\CalDSeip^{k\ssp+\sp 1.}$ it suffices to establish \œ$\varSei^{}Ï\tilde f
\in\Cal D\sbi k\,$. To get this, by induction on \œ$l\in k\ssp\yplus$ one
proves (x) that for $x\in\dom\sn f$ and $u\in\vecs E$ and $\bosy x\ssp,\bosy u
\in(\sp\vecs E\ssp)\,^l$ we have \vskip.6mm

$\taurd\varSei^\sp lÏ\varSei^{}Ï\tilde f\sp\fvalue\big(\spp
\seq{\,(\sp x\ssp,u\sp)\,}\concc\ssp [\,\bosy x\ssp,\bosy u\,]\ssp\big) = $ \vskip.3mm $\null\hfill
\taurd\varSei^\sp l\ssp+\sp 1.Ï\tilde f\sp\fvalue\big(
                               \seq{\,x\ssp,u\,}\concc\bosy x\sp\big) +
\sum_{\,i\ssp\in\ssp l}\taurd\varSei^\sp lÏ\tilde f\sp\fvalue\sp
            \roman R\,x\,\bosy u\,i\,\bosy x\,$.\KP{15}\vskip.5mm

\noin For the inductive step here, noting ($*$) above, and also letting \œ$\sp
y\ssp,v\in\vecs E\,$, one checks the computation\vskip.7mm

$\taurd\varSei^\sp l\ssp+\sp 1.Ï\varSei^{}Ï\tilde f\sp\fvalue
  \big(\seq{\,(\sp x\ssp,u\sp)\,}\concc\sp [\,\bosy x\ssp,\bosy u\,]\sp
   \concc\seq{\,(\sp y\ssp,v\sp)\,}\sp\big)$ \vskip.5mm $\mhyppy{11.7}
= \taurd F\text{\,-\,}\underset{t\to 0}\lim\ t^{\sp\mminus 1}\sp
   \big(\sp\taurd\varSei^\sp lÏ\varSei^{}Ï\tilde f\sp\fvalue\big(\spp
    \seq{\,(\sp x+t\,y\ssp,u+t\,v\sp)\,}\concc\ssp [\,\bosy x\ssp,
                               \bosy u\,]\ssp\big)$\vskip-1mm $\mhyppy{51}
 - \taurd\varSei^\sp lÏ\varSei^{}Ï\tilde f\sp\fvalue\big(\spp\seq{\,
   (\sp x\ssp,u\sp)\,}\concc\ssp [\,\bosy x\ssp,\bosy u\,]\ssp\big)\big)$ \vskip.5mm $
= \taurd F\text{\,-\,}\underset{t\to 0}\lim\ t^{\sp\mminus 1}\sp\big(\sp
   \taurd\varSei^\sp l\ssp+\sp 1.Ï\tilde f\sp\fvalue\big(\seq{\,x+t\,y\ssp,
                   u+t\,v\,}\concc\bosy x\sp\big)
 - \taurd\varSei^\sp l\ssp+\sp 1.Ï\tilde f\sp\fvalue\big(
               \seq{\,x\ssp,u\,}\concc\bosy x\sp\big) $ \vskip-.5mm $\null\hfill
 + \sum_{\,i\ssp\in\ssp l}\sp\big(\ssp
      \taurd\varSei^\sp lÏ\tilde f\sp\fvalue\sp
           \roman R\,(\sp x+t\,y\sp)\,\bosy u\,i\,\bosy x
    - \taurd\varSei^\sp lÏ\tilde f\sp\fvalue\sp
           \roman R\,x\,\bosy u\,i\,\bosy x\ssp\big)\big)$\KP{4.5}\vskip.5mm$
= \taurd\varSei^\sp l\ssp+\ssp 2.Ï\tilde f\sp\fvalue\big(
     \seq{\,x\ssp,u\,}\concc\bosy x\concc\seq{\ssp y\ssp}\sp\big)
 + \taurd\varSei^\sp l\ssp+\sp 1.Ï\tilde f\sp\fvalue\big(
                           \seq{\,x\ssp,v\,}\concc\bosy x\sp\big)$ \par $\mhyppy{46.5}
 + \sum_{\,i\ssp\in\ssp l}
    \taurd\varSei^\sp l\ssp+\sp 1.Ï\tilde f\sp\fvalue\sp
    \roman R\,x\,\bosy u\,i\,(\sp\bosy x\concc\seq{\ssp y\ssp}\spp)$\vskip.5mm $
= \taurd\varSei^\sp l\ssp+\ssp 2.Ï\tilde f\sp\fvalue\big(
   \seq{\,x\ssp,u\,}\concc\bosy x\concc\seq{\ssp y\ssp}\sp\big)
 +  \sum_{\,i\ssp\in\ssp l\ssp+\sp 1.}
         \taurd\varSei^\sp l\ssp+\sp 1.Ï\tilde f\sp\fvalue\sp
         \roman R\,x\,(\sp\bosy u\concc\seq{\ssp v\ssp}\spp)\,i\,
                      (\sp\bosy x\concc\seq{\ssp y\ssp}\spp)\,$.\vskip.7mm

Finally, from (x) we \q{directly} see that $\varSei^lÏ\varSei^{}Ï\tilde f \in
\CalCSe0$ with also\vskip.2mm\centerline{$
\{\,\biit z:\biit z\fvalue\emptyset\in\dom\taurd\varSei^{}Ï\tilde f\ssp\,\}
\cap(\sp\vecs(E\skcap E\ssp))\,^{l\ssp+\sp 1.} \inc
             \dom\taurd\varSei^lÏ\varSei^{}Ï\tilde f$}\vskip.3mm

\noin for all $\sp l\in k\ssp\yplus\sp$. Consequently $\sp\varSei^{}Ï\tilde f
\in\Cal D\sbi k\,$, \,as was to be established.                  \end{proof}

A not quite short analytic\ssp-\ssp combinatorial proof shows that for an
order $l$ simply Seip\,--\,differentiable map \œ$\tilde f=(E\ssp,F\sp,f\sp)$
and for \œ$x\in\dom\sn f$ and \œ$\bosy x\in(\sp\vecs E\ssp)\,^l$ and\linebreak
any bijection \œ$\sigma:l\to l$ we have $
\taurd\varSei^\sp lÏ\tilde f\sp\fvalue
                       \big(\sp\seq{\ssp x\ssp}\concc\spp\bosy x\ssp\big) =
\taurd\varSei^\sp lÏ\tilde f\sp\fvalue
   \big(\sp\seq{\ssp x\ssp}\concc\spp(\sp\bosy x\circ\sigma\sp)\big)\,$. If we
had first established this result, in the above proof we could have written
also\vskip.3mm\centerline{$
\taurd\varSei^\sp lÏ\tilde f\sp\fvalue\sp\roman R\,x\,\bosy u\,i\,\bosy x =
\taurd\varSei^\sp lÏ\tilde f\sp\fvalue\big(\sp\seq{\ssp x\ssp}\concc
             (\sp\bosy x\,|\,(\sp\Univ\setminus\sn\{\sp i\sp\})\cup
              \{\ssp(\sp i\ssp,\bosy u\fvalue i\sp)\ssp\}\spp)\big)\,$.}

\begin{construction}

With any classes $\,\tilde g\,,\sp k\,$ we associate\vskip.4mm

$I^{\,k}\sp\tilde g\ssp=\ssp\bigcap\sp\big\{\ssp(E\skcap E\,,F\sp,\spp g\ar 1)
:k\in\No$ and $E\ssp,F$ are real topologized \par\nopagebreak\hyppy{7mm}
vector spaces and $\sp g\in F^{\sp/\sp E}$ and $\sp\tilde g =
(E\ssp,F\sp,\spp g\sp)$ and $g\ar 1 = \big\{\ssp(\sp x\ssp,u\ssp,y\sp):$ \par\nopagebreak\centerline{$
[\ \all s\,0\le s\le 1\imply x+s\,u\in\dom g$ ] and $y =
      F\text{\,-\sp}\int_{\,0}^{\,1}s^{\,k}\sp
      g\fvalue(\sp x+s\,u\sp)\ssp\,\roman d\,s\,\big\}\sp\big\}\sp\,$.}

  \end{construction}

By elementary set theoretic manipulations, we have either \œ$
I^{\,k}\sp\tilde g=\bigcap\ssp\emptyset=\Univ\,$, or \œ$k\in\No$ and there are
real topologized vector spaces $E\ssp,F$ and \œ$\sp g\in F^{\sp/\sp E}$ with \œ$\sp
\tilde g=$ $(E\ssp,F\sp,\spp g\sp)\,$. In the latter case \œ$
I^{\,k}\sp\tilde g=(E\skcap E\,,F\sp,\spp g\ar 1)$ where $g\ar 1$ is the
function defined on the set of the $(\sp x\ssp,u\sp)$ with \œ$x\in\dom g$ and
\œ$u\in\vecs E$ such that \œ$x+s\,u\in\dom g$ for \œ$0\le s\le 1\ssp$, and
such that also the function \œ$s\mapsto s^{\,k}\sp g\fvalue(\sp x+s\,u\sp)$ on
$[\,0\,,1\,]$ in $F$ is Riemann integrable to some $y\,$. Deeper properties
important for our purposes of the integration families $
 J:\CalCSe0\owns\tilde g\mapsto I^{\,k}\sp\tilde g$ are given by the following

\pagebreak

\begin{lemma}\label{g in D^i imply Ig in D^i when i=0 or i=1}

Let \œ$\ssp k\in\No$ and \œ$\ssp J=\seq{\,I^{\,k}\sp\tilde g:\tilde g \in
\CalCSe0\ssp}\,$. Then \œ$\ssp\rng J\inc\CalCSe0=\dom J$ and \œ$\ssp
J\ssp\image\sp\CalDSeip^{1.}\inc\CalDSeip^{1.}$ hold. Also if \œ$\,\tilde g =
(E\ssp,F\sp,g\sp)\in\CalDSeip^{1.}\,,$ for all $\ssp x\ssp,y\ssp,u\ssp,v$ such
that \œ$\ssp(\sp x\ssp,u\sp)\in\dom\taurd\spp I^{\,k}\sp\tilde g$ and $\ssp
y\ssp,v\in\vecs E\,$ it holds that \ $
\taurd\varSe I^{\,k}\sp\tilde g\fvalue(\sp x\ssp,u\,;y\ssp,v\sp) = $\vskip.5mm\centerline{$
F\text{\,-\sp}\int_{\,0}^{\,1}(\sp s^{\,k}\sp\taurd\varSe\tilde g\fvalue(\sp
x+s\,u\ssp,y\sp)+s^{\,k\ssp+\sp 1.}\sp\taurd\varSe\tilde g\fvalue(\sp x +
s\,u\ssp,v\sp))\ssp\,\roman d\,s\,$.}                            \end{lemma}

\begin{proof} Write \œ$I=[\,0\,,1\,]\,$. Arbitrarily fixing \œ$\tilde g =
(E\ssp,F\sp,g\sp)\in\CalCSe0\,$, we first show that \œ$I^{\,k}\sp\tilde g
\in\CalCSe0\,$. Putting \œ$G=E\skcap E$ and \œ$F\aar 1 =
\taukv\KN1\fvalue\snn F$ and \œ$g\ar 1=\taurd\spp I^{\,k}\sp\tilde g\,$, note
that since $\taurd G$ is a compactly generated topology, and since we have \œ$
F=\kappatv\KN1\fvalue\snn F\aar 1\ssp$, it suffices that $
(\ssp G\sp,F\aar 1\sp,\sp g\ar 1)$ has open domain and is continuous. Hence,
arbitrarily fixing \œ$z=(\sp x\ssp,u\sp)\in\dom g\ar 1$ and a closed convex \œ$
V\in\ymp F\aar 1\ssp$, it suffices to show (e) existence of \œ$
W\in\Nbh(\sp z\ssp,\taurd G\sp)$ such that for all \œ$w\in W$ we have \œ$
g\ar 1\KN1\fvalue w - g\ar 1\KN1\fvalue z\in V\spp$. Note that this implies \œ$
W\inc\dom g\ar 1$ since by \cite[Theorem 69, p.\ 261]{Ky} we have $
g\ar 1\KN1\fvalue w - g\ar 1\KN1\fvalue z = \Univ \not\in V$ if $w\not\in
\dom g\ar 1\ssp$.

To establish (e) above, we consider the set\vskip.3mm

$\mhyppy3  \Cal N\aar 0 = \{\ssp(\ssp W_{\spp},N\sp):W\in\Nbh(\sp z\ssp,
\taurd G\sp)$ and $\sp\exi{t\in I}\,N\in\Nbh(\sp t\ssp,\tau_{_{I\!\!R}})$ and \par\nopagebreak
$\mhyppy{38} \all{x\yr 1\spp,u\yr 1\spp,s}\,(\sp x+x\yr 1\spp,u+u\yr 1)\in W$
              and $s\in I\cap N$ \par\nopagebreak
$\mhyppy{38.5}\imply g\fvalue(\sp x+x\yr 1\snn+s\,(\sp u+u\yr 1))
                   - g\fvalue(\sp x+s\,u\sp)\in V\sp\,\}\,$.\vskip.3mm

\noin Putting \œ$m = \seq{\,x+s\,u:\smb X=(\sp x\ssp,u\ssp,s\sp)\in\vecs(\spp
G\sp\sqcap\biit R\ssp)\,}\,$, note that $m$ is a continuous second order
polynomial function $G\sp\sqcap\biit R\to E$ since in its decomposition\vskip.3mm

$\mhyppy{13.5}\smb X=(\sp x\ssp,u\ssp,s\sp)\mapsto(\sp x\,;s\ssp,u\sp) =
(\sp x\ssp,\smb Z\sp)\mapsto(\sp x\ssp,\tsigrd E\sp\fvalue\smb Z\sp)$

$\mhyppy{52} = (\sp x\ssp,v\sp)=\smb Y\mapsto\ssigrd E\sp\fvalue\sp\smb Y =
               m\fvalue\smb X$\vskip.3mm

\noin besides the continuous linear projections as factors there are the
continuous linear and bilinear \œ$(\spp G\sp,E\ssp,\ssigrd E\sp)$ and \œ$
(\biit R\sp\sqcap E\ssp,E\ssp,\tsigrd E\sp)\,$. Consequently, since $
(E\ssp,F\aar 1\sp,g\sp)$ has open domain and is continuous, the same holds
also for $(\spp G\sp\sqcap\biit R\,,F\aar 1\sp,\spp g\circ m\sp)\,$.

Using this, we see that \œ$I\inc\bigcup\sp\rng\Cal N\aar 0\ssp$, whence $
\tau_{_{I\!\!R}}\ssp$--\,compactness of $\sp I$ gives existence of a finite \œ$
\Cal N\inc\Cal N\aar 0$ with \œ$I\inc\bigcup\sp\rng\Cal N\sp$. Putting \œ$W =
\bigcap\sp\dom\Cal N\sp$, we have \œ$W\in\Nbh(\sp z\ssp,\taurd G\sp)\,$, and
for (e) it now remains to verify that for arbitrarily fixed $w =
z+(\sp x\yr 1\spp,u\yr 1)\in W$ we have $
g\ar 1\KN1\fvalue w-(\sp g\ar 1\KN1\fvalue z\sp)\in V\spp$. Letting\vskip.3mm\centerline{$
 \gamma = \seq{\,s^{\,k}\sp(\ssp g\fvalue(\sp x+x\yr 1\snn+s\,(\sp u+u\yr 1))
                 - g\fvalue(\sp x+s\,u\sp)):s\in I\sp\,}\,$,}\vskip.3mm

\noin we have \œ$g\ar 1\KN1\fvalue w-(\sp g\ar 1\KN1\fvalue z\sp) =
F\text{\,-}\int_{\,0}^{\,1}\gamma = F\aar 1\text{\ssp-}\int_{\,0}^{\,1}\gamma\,
$, whence by our Lemma \ref{int-mvT} it suffices\Biggerlineskip1 that $\rng
\gamma\inc V\spp$, which in turn directly follows from our arrangement above.

Next assuming that \œ$\sp\tilde g\in\CalDSeip^{1.}\,$, we show that for \œ$\sp
z=(\sp x\ssp,u\sp)\in\dom g\ar 1$ and \œ$\ssp w=$ $(\sp y\ssp,v\sp)$ with \œ$
y\ssp,v\in\vecs E\,$ we have the asserted variation formula $\sp
\taurd\varSe I^{\,k}\sp\tilde g\fvalue(\sp z\ssp,w\sp)$ \œ$ = \smb Y\ssp$,
putting \œ$\ssp\smb Y\sp=
 F\text{\,-\sp}\int_{\,0}^{\,1}(\sp s^{\,k}\sp\taurd\varSe\tilde g\fvalue(\sp
 x+s\,u\ssp,y\sp) + s^{\,k\ssp+\sp 1.}\sp\taurd\varSe\tilde g\fvalue(\sp x +
 s\,u\ssp,v\sp))\ssp\,\roman d\,s\,$.\linebreak Noting that the topologies $
\taurd F$ and $\taurd F\aar 1$ have the same convergent sequences, one quickly
deduces that for arbitrarily given closed convex \œ$V\in\ymp F\aar 1\ssp$, it
suffices to show (x) existence of some \œ$\ssp\delta\in\Rep$ such that for all
\œ$t\in\Re$ with \œ$0<|\ssp t\ssp|<\delta$ we have $\,
    t^{\sp\mminus 1}\sp(\ssp g\ar 1\KN1\fvalue(\sp z+t\,w\sp) -
                             g\ar 1\KN1\fvalue z\sp) - \smb Y\in V\spp$. \vskip.3mm

To establish (x) above, letting\vskip.3mm\centerline{$
\Gamma\,v\ar 1\ssp s\,t =
  \taurd\varSe\tilde g\fvalue(\sp x+s\,u+t\,(\sp y+s\,v\sp)\ssp,v\ar 1)
    - \taurd\varSe\tilde g\fvalue(\sp x+s\,u\ssp,v\ar 1)\,$,}\vskip.3mm

\noin  for each fixed \œ$v\ar 1\in\vecs E$ the function \œ$\sp\Gamma\,v\ar 1\sn
:\smb S=(\sp s\ssp,t\sp)\mapsto\Gamma\,v\ar 1\ssp s\,t$ defined exactly for
the \œ$\smb S\in\Re\,^{\times2.}$ having \œ$x+s\,u+t\,(\sp y+s\,v\sp)\in\dom g
$ is continuous \œ$\ssp\Cal T=\tau_{_{\Re}}\sn\ttimes\tau_{_{\Re}}$\linebreak
\œ$\to\taurd F\aar 1$ and satisfies \œ$I\timesn\{0\} \inc
(\sp\Gamma\,v\ar 1)\inve\image\snn\{\ssp\bnull F\}$ and \œ$
\dom\snn(\sp\Gamma\,v\ar 1)\in\Cal T\spp$. Using this \linebreak in
conjunction with compactness of $\sp I\sp$, one deduces existence of \œ$\ssp
\delta\in\Rep$ such that we have $\,\Gamma\,v\,s\,t\ssp,\Gamma\,y\,s\,t \in
\frac 12\,V$ for all $s\in I$ and $t\in\Re$ with $|\ssp t\ssp|<\delta\ssp$.

We now arbitrarily fix \œ$t\in\Re$ with \œ$0<|\ssp t\ssp|<\delta\ssp$, and
consider in the space $F\aar 1$ the curve \œ$\, c = \seq{\,t^{\sp\mminus 1}\sp
(\ssp g\ar 1\KN1\fvalue(\sp z+s\ar 1\ssp t\,w\sp)
-  g\ar 1\KN1\fvalue z\sp) - s\ar 1\ssp\smb Y : s\ar 1\in I\ssp\,}\,$. Note
that for every \œ$s\ar 1\in I$ we have \œ$x+s\,u+s\ar 1\ssp t\,(\sp y+s\,v\sp)
\in\dom g$ for all \œ$s\in I\sp$, and also that the function \œ$s \mapsto
g\fvalue(\sp x+s\,u+s\ar 1\ssp t\,(\sp y+s\,v\sp))$ is continuous. Recalling
Proposition \ref{conti is int} above, we hence have \œ$s\ar 1\in\dom c\,$, and
thus indeed \œ$\dom c=I\sp$. Noting that we have got (x) once \œ$\,c\fvalue 1
\in V$ is established, by the
{\it mean value theorem\ssp} (Proposition \ref{mvT} above) it suffices to show
that $c\sp$ is differentiable with $\rng\cder_{F_1}\sn c\inc V\spp$.

For this arbitrarily fixing \œ$s\ar 1\in I\sp$, using
Proposition \ref{diff under int} and noting the linearity from the proof of
(b) in Corollary \ref{D^1 iff BGN^1 iff Seip^1(1972)} above, a direct
computation gives\vskip.3mm\centerline{$
\cder_{F_1}\sn c\fvalue s\ar 1 =  F\aar 1\text{\ssp-}\int_{\,0}^{\,1}(\sp
    s^{\,k\,}\Gamma\,y\,s\,(\sp s\ar 1\ssp t\sp)
 +  s^{\,k\ssp+\sp 1.}\ssp\Gamma\,v\,s\,(\sp s\ar 1\ssp t\sp)
                            )\ssp\,\roman d\,s\,$,}\vskip.3mm

\noin whence by Lemma \ref{int-mvT} we get $\,
      \cder_{F_1}\sn c\fvalue s\ar 1\in\frac 12\,V+\frac 12\,V\inc V\spp$. \vskip.3mm

We have now (x) and by the variation formula we also have the decomposition\vskip.3mm\centerline{$
\taurd\varSe I^{\,k}\sp\tilde g = \ssigrd F\circ(
(\sp\taurd\sp I^{\,k}\sp\varSe\tilde g\ssp)\risti2
(\sp\taurd\sp I^{\,k\ssp+\sp 1.}\sp\varSe\tilde g\ssp))\circ\sp\ell\,$,}\vskip.3mm

\noin where the continuous linear map \œ$\ssp
       \ell : (E\skcap E\ssp)\skcap(E\skcap E\ssp) = E\ar 4 \to
                                    E\ar 4\snn\skcap\spp E\ar 4$ is given by \œ$
 (\sp x\ssp,u\,;y\ssp,v\sp) \mapsto
 (\sp x\ssp,y\,;u\ssp,\bnull E\,;(\sp x\ssp,v\,;u\ssp,\bnull E))\,$. Having
already established the inclusion $\rng J\inc\CalCSe0\,$, this gives $
\varSe I^{\,k}\sp\tilde g\in\CalCSe0$ when applied also to \œ$k+1\adot$ in
place\linebreak of $\sp k\ssp$, and consequently by
Corollary \ref{D^1 iff BGN^1 iff Seip^1(1972)}\ssp(c) also $\ssp
J\ssp\image\sp\CalDSeip^{1.}\inc\CalDSeip^{1.}$ is verified.     \end{proof}

In proving that $\CalDSeip^k=\CalDBGN^k(\sp\CalDSeip^{0.}\ssp,\spp\biit R\ssp)
$ we shall utilize the following

\begin{corollary}\label{g in D^k imply Ig in D^k}

For any $\sp k$ and $\ssp l\in\No$ and $\ssp\tilde g\in\CalDSeip^k\,,$ it
holds that $\sp I^{\,l}\sp\tilde g\in\CalDSeip^k\ssp$.       \end{corollary}

\begin{proof} Since \œ$\sp\CalDSeip^k=\emptyset$ if \œ$k \not\in
\infty\ssp\yplus\sp$, it suffices to prove by induction that for \œ$k\in\No$
the assertion holds for all $\sp l$ and $\sp\tilde g\,$. The case \œ$k =
\emptyset$ or \œ$k=1\adot$ having been settled by
Lemma \ref{g in D^i imply Ig in D^i when i=0 or i=1} above, assuming that the
assertion holds for a fixed \œ$k\in\N\sp$, and that \œ$\ssp\tilde g =
(E\ssp,F\sp,g\sp)\in\CalDSeip^{k\ssp+\sp 1.}\,$, by
Proposition \ref{D^k recursion} we have \œ$\sp\tilde g\in\CalDSeip^{1.}$ and\linebreak
\œ$\varSe\sp\tilde g\in\CalDSeip^k\,$. By the inductive assumption then \œ$\sp
I^{\,\sp l}\sp\varSe\tilde g\,,\sp I^{\,\sp l\ssp+\sp 1.}\sp\varSe\tilde g\in
\CalDSeip^k\,$, and by Lemma \ref{g in D^i imply Ig in D^i when i=0 or i=1}
also \œ$\sp I^{\,\sp l}\sp\tilde g\in\CalDSeip^{1.}$ holds together with the
decomposition\vskip.3mm\centerline{$
 \taurd\varSe I^{\,\sp l}\sp\tilde g =
 \ssigrd F\circ((\sp\taurd\sp I^{\,\sp l}\sp\varSe\tilde g\ssp)\risti2
 (\sp\taurd\sp I^{\,\sp l\ssp+\sp 1.}\sp\varSe\tilde g\ssp))\circ\sp\ell\,$,}\vskip.3mm

\noin where the continuous linear map \œ$
(E\ar 4\ssp,E\ar 4\snn\skcap\spp E\ar 4\ssp,\sp\ell\,)$ is as in the proof of
Lemma \ref{g in D^i imply Ig in D^i when i=0 or i=1} above. By
Proposition \ref{D^k recursion} and items (\ref{pr1b})\ssp, (\ref{pr1e}) and
(\ref{pr1d}) of Proposition \ref{basic D^k properties} we hence obtain the
membership $\ssp I^{\,\sp l}\sp\tilde g\in\CalDSeip^{k\ssp+\sp 1.}$ as
required.                                                        \end{proof}

\begin{theorem}\label{Seip = BGN}

For all $\ssp k\in\infty\ssp\yplus$ it holds that $\ssp\CalDSeip^k =
\CalDBGN^k(\sp\CalDSeip^{0.}\ssp,\spp\biit R\ssp)\,$.          \end{theorem}

\begin{proof} It clearly suffices to treat the cases \œ$\ssp k\in\infty =
\No\ssp$. Proceeding by induction, we first note that the case \œ$k=\emptyset$
is trivial, and that the case \œ$k=1\adot$ has been settled in (a) of
Corollary \ref{D^1 iff BGN^1 iff Seip^1(1972)} above.

Let now \œ$k\in\N$ with \œ$\CalDSeip^k =
\CalDBGN^k(\sp\CalDSeip^{0.}\ssp,\spp\biit R\ssp)\,$. If \œ$\tilde f \in
\CalDBGN^{k\ssp+\sp 1.}(\sp\CalDSeip^{0.}\ssp,\spp\biit R\ssp)\,$, then \œ$
\tilde f\in\CalDBGN^{1.}(\sp\CalDSeip^{0.}\ssp,\spp\biit R\ssp) \inc
\CalDSeip^{1.}$ and \œ$\cgbarDelta\tilde f \in
\CalDBGN^k(\sp\CalDSeip^{0.}\ssp,\spp\biit R\ssp)\inc\CalDSeip^k\,$. Putting \œ$
\tilde\ell=$\linebreak \œ$(E\skcap E\,,E\skcap E\skcap\biit R\,,\ssp\ell\,)$
where \œ$\ell=\seq{\,(\sp x\ssp,u\ssp,0\sp):z=(\sp x\ssp,u\sp)\in(\sp\vecs
E\sp)^{\sp\times 2.}\ssp}\,$, we have $\ssp\tilde\ell$\Biggerlineskip2 a
continuous linear map, hence \œ$\ssp\tilde\ell\in\CalDSeip^k\,$, and also \œ$
\varSe\tilde f=\cgbarDelta\tilde f\snn\Circ\sp\tilde\ell\,$, whence the chain
rule gives $\varSe\tilde f\in\CalDSeip^k\,$. Using $\tilde f\in\CalDSeip^{1.}$
we get $\tilde f\in\CalDSeip^{k\ssp+\sp 1.}\,$.

To establish the converse, assuming that \œ$\tilde f=(E\ssp,F\sp,f\sp)\in
\CalDSeip^{k\ssp+\sp 1.}\ssp$, we have \œ$\tilde f\in$ $
\CalDBGN^{1.}(\sp\CalDSeip^{0.}\ssp,\spp\biit R\ssp)\,$, whence to get \œ$
\tilde f\in\CalDBGN^{k\ssp+\sp 1.}(\sp\CalDSeip^{0.}\ssp,\spp\biit R\ssp)\,$,
it suffices to obtain \œ$\cgbarDelta\tilde f\in$ \œ$\CalDSeip^k\inc
\CalDBGN^k(\sp\CalDSeip^{0.}\ssp,\spp\biit R\ssp)\,$. By
Proposition \ref{basic D^k properties}\ssp(\ref{pr1l})\ssp, considering \œ$
\smb W=(\sp x\ar 0\ssp,u\ar 0\ssp,t\ar 0\ssp,y\sp) \in$ \œ$ f\yr 1 =
\taurd\cgbarDelta\tilde f\sp$, with \œ$G=E\skcap E\skcap\biit R\,$, it
suffices to establish $h$ with \œ$\smb W\in h\inc f\yr 1$ and $
(\spp G\sp,F\sp,h\sp)\in\CalDSeip^k\ssp$. If $t\ar 0\not=0\,$, we may take $
 h = f\yr 1\ssp|\,
   (\ssp\roman{pr}\ar 2\KN1\inve\ssp[\,\sp\Univ\setminus\sn\{0\}\,]\ssp)\,$.

To handle the case \œ$t\ar 0=0\,$, letting \œ$\tilde m =
(\spp G\sp,E\skcap E\skcap(E\skcap E\sp)\ssp,m\sp)\,$, where $m = $ \œ$
\seq{\,(\sp x\ssp,u\,;t\,u\ssp,\bnull E):\smb X=(\sp x\ssp,u\ssp,t\sp)\in\vecs
G\sp\,}\,$, we have $\tilde m$ a continuous second order polynomial map, and
hence \œ$\tilde m\in\CalDSeip^k\ssp$. By \œ$\tilde f \in
\CalDSeip^{k\ssp+\sp 1.}$ we have \œ$\varSe\tilde f\in\CalDSeip^k\,$, whence
by Corollary \ref{g in D^k imply Ig in D^k} above we get \œ$
I^{\,\sp 0.}\varSe\tilde f\in\CalDSeip^k\,$, and consequently for \œ$\tilde h=
(\spp G\sp,F\sp,h\sp)$ $=I^{\,\sp 0.}\varSe\tilde f\snn\Circ\tilde m$ the
chain rule gives $\tilde h\in\CalDSeip^k\ssp$. For $\smb X =
(\sp x\ssp,u\ssp,t\sp)\in\dom h$ having\vskip.5mm

$\mhyppy{10.63} h\fvalue\smb X
= \sp\taurd I^{\,\sp 0.}\varSe\tilde f\sp\fvalue(\sp m\fvalue\smb X\sp)
= \sp\taurd I^{\,\sp 0.}\varSe\tilde f\sp\fvalue
                                     (\sp x\ssp,u\,;t\,u\ssp,\bnull E)$\vskip.6mm

\centerline{$\phantom{h\fvalue\smb X} =
 F\text{\,-}\int_{\,0}^{\,1}\taurd\varSe\tilde f\sp\fvalue(\sp x+s\,t\,u\ssp,
                                                   u\sp)\ssp\,\roman d\,s =
 f\yr 1\fvalue(\sp x\ssp,u\ssp,t\sp) = f\yr 1\fvalue\smb X\,$,}\vskip.6mm

\noin we see that $\smb W\in h\inc f\yr 1$, completing the inductive step, and
the whole proof.                                                 \end{proof}

\begin{lemma}\label{f(,) D^k imply f^v D^k}

Let \œ$\ssp E\ssp,F\sp,G\in\scgVS(\biit R\ssp)$ and \œ$\ssp L =
\Lincg(\spp F\sp,G\ssp)$ and \œ$\,U\sn\in\taurd E$ and also \œ$f\aar 1 \in
\big((\sp\vecs G\ssp)^{\,\upsilon_s\sp F\,}\big){^{\,U}}$ and \œ$\ssp f =
f\aar 1\KN1\ywed\sp$. For the condition {\ssp\rm(n)} below, for all \œ$\ssp k
\in\No$ it then holds that \ $[\ \text{(n) and }\ssp
(E\spp\skcap F\sp,G\sp,f\sp)\in\CalDSeip^k\ ]\sp\equivv\sp
(E\ssp,L\ssp,f\aar 1)\in\CalDSeip^k\,$.\vskip.5mm

{\rm(n)} \ \ $f\aar 1\KN1\fvalue x$ is linear $\sigrd F\to\sigrd G$
              for every $x\in U\spp$.                            \end{lemma}

\begin{proof} If \œ$\ssp k=\emptyset\,$, the assertion follows from
Proposition \ref{lin expo} above. We next prove the assertion for \œ$\sp k =
1\adot\ssp$. Let \œ$\sp\tilde f=(E\spp\skcap F\sp,G\sp,f\sp)$ and \œ$\ssp
\tilde f\aar 1=(E\ssp,L\ssp,f\aar 1)\,$. First assum- ing (n) and that \œ$\sp
\tilde f\in\CalDSeip^{1.}\,$, to prove that \œ$\sp\tilde f\aar 1 \in
\CalDSeip^{1.}\,$, by Proposition \ref{D^k recursion} it suffices that $
\tilde f\aar 1$ is directionally Seip\,--\,differentiable and that $
\varSe\tilde f\aar 1$ is continuous.\vskip.3mm

Putting \œ$g\ar 1=\seq{\,\seq{\,\taurd\varSe\tilde f\sp\fvalue(\sp x\ssp,v\,;
u\ssp,\bnull F):v\in\vecs F\sp\,}\subtext{old}\sn:\smb Z=(\sp x\ssp,u\sp)\in
O\sp\,}\,$, where\linebreak \œ$O=U\snn\times\snn(\sp\vecs E\ssp)\,$, to
establish directional differentiability with \œ$\taurd\varSe\tilde f\aar 1 =
g\ar 1\ssp$, we\linebreak arbitrarily fix \œ$\smb Z=(\sp x\ssp,u\sp)\in O\ssp
$, and take any \œ$\sp G\ar 0 \in
\kappatv\KN1\inve\image\sn\{\ssp G\ssp\}\cap\LCS(\biit R\ssp)\,$. In view of
Lemma \ref{seq conv in L_co(F,G)} it then suffices for arbitrarily given $
\taurd F\,$--\,compact $K$ and closed convex \œ$\sp V\sn\in\ymp G\ar 0$ to
show (e) existence of \œ$\ssp\delta\in\Rep$ such that for all \œ$t\in\Re$ and\linebreak
$v\in K$ with $0<|\ssp t\ssp|<\delta\sp$ we have $\ssp
\varDelta\ssp(\sp t\ssp,v\sp)\in V\spp$, upon putting\vskip.3mm\centerline{$
\varDelta\ssp(\sp t\ssp,v\sp) = t^{\sp\mminus 1}\sp(\sp
 f\sp\fvalue(\sp x+t\,u\ssp,v\sp) - f\sp\fvalue(\sp x\ssp,v\sp)) -
 \taurd\varSe\tilde f\sp\fvalue(\sp x\ssp,v\,;u\ssp,\bnull F)\,$.}\vskip.3mm

To establish (e)\ssp, using continuity of $\sp\varSe\tilde f$ and of the
algebraic operations of $E$ and compactness of $\sp K\sp$, we first see that
there is $\delta\in\Rep$ such that\vskip.3mm

\noin$(*)$ \hfill $
  \taurd\varSe\tilde f\sp\fvalue(\sp x+t\,u\ssp,v\,;u\ssp,\bnull F)
- \taurd\varSe\tilde f\sp\fvalue(\sp x\ssp,v\,;u\ssp,\bnull F)\in V$ \hfill \phantom{$(*)$}\vskip.3mm

\noin whenever \œ$t\in\Re$ and \œ$v\in K$ with \œ$|\ssp t\ssp|<\delta\ssp$.
To get \œ$\varDelta\ssp(\sp t\ssp,v\sp)\in V$ for arbitrarily fixed such $
t\ssp,v$ with \œ$t\not=0\,$, by the mean value theorem considering the
differentiable curve $c\,$ on $\ssp[\,0\,,1\,]\sp$ with values in $\sp G\ar 0$
and defined by\vskip.3mm\centerline{$
s \mapsto t^{\sp\mminus 1}\sp(\sp
  f\sp\fvalue(\sp x+s\,t\,u\ssp,v\sp) - f\sp\fvalue(\sp x\ssp,v\sp)) -
  s\,\sp\taurd\varSe\tilde f\sp\fvalue(\sp x\ssp,v\,;u\ssp,\bnull F)\,$,}\vskip.3mm

\noin it suffices that \œ$\ssp\rng\cder_{G_0}\sn c\inc V\spp$. By $(*)$ this
indeed is the case. Having now \œ$\,\varSe\tilde f\aar 1 = $ \œ$
(E\spp\skcap E\ssp,L\ssp,g\ar 1)\,$, to prove that $\varSe\tilde f\aar 1$ is
continuous, we just note that this directly follows from
Proposition \ref{lin expo} by continuity of $\sp\varSe\tilde f\sp$.

Next assuming that \œ$\sp\tilde f\aar 1\in\CalDSeip^{1.}\,$, trivially (n)
holds, and to prove that \œ$\sp\tilde f\in\CalDSeip^{1.}\,$, we again
establish directional differentiability and continuity of the variation. With
$\,O\ar 2=U\snn\times\snn(\sp\vecs F\ssp)\snn\times\snn(\sp\vecs E\snn\times\snn
(\sp\vecs F\ssp))\,$, \,letting\vskip.3mm\centerline{$
g=\seq{\,\taurd\varSe\tilde f\aar 1\KN1\fvalue(\sp x\ssp,u\sp)\fvalue v +
         f\sp\fvalue(\sp x\ssp,v\aar 1) : \smb W =
         (\sp x\ssp,v\,;u\ssp,v\aar 1)\in O\ar 2\ssp}\,$,}\vskip.3mm

\noin to prove directional differentiability together with the equality \œ$\sp
\taurd\varSe\tilde f=g\,$, for arbitrarily fixed $\,\smb W=(\sp x\ssp,v\,;
u\ssp,v\aar 1)\in O\ar 2\ssp$, for $0\not=t\in\Re$ with $|\ssp t\ssp|$ small,
having\vskip.5mm\centerline{$
t^{\sp\mminus 1}\sp(\sp f\sp\fvalue(\sp x+t\,u\ssp,v+t\,v\aar 1) -
                        f\sp\fvalue(\sp x\ssp,v\sp)) -
\taurd\varSe\tilde f\aar 1\KN1\fvalue(\sp x\ssp,u\sp)\fvalue v -
         f\sp\fvalue(\sp x\ssp,v\aar 1)$}\vskip.5mm

$\mhyppy{2.2} = t^{\sp\mminus 1}\sp(\sp f\sp\fvalue(\sp x+t\,u\ssp,v\sp) -
                                        f\sp\fvalue(\sp x\ssp,v\sp)) -
     \taurd\varSe\tilde f\aar 1\KN1\fvalue(\sp x\ssp,u\sp)\fvalue v$

$\mhyppy{66.7} + f\sp\fvalue(\sp x+t\,u\ssp,v\aar 1) -
                 f\sp\fvalue(\sp x\ssp,v\aar 1)\,$,\vskip.4mm

\noin we see that directional differentiability is immediate by continuity of
$\sp\tilde f$ following from that of $\sp\tilde f\aar 1$ by
Proposition \ref{lin expo} above. Continuity of $\sp\varSe\tilde f$ follows
similarly.

Finally, for fixed $F\sp,G$ letting $\ssp\roman P\ssp(k\spp)$ mean that the
assertion of the Lemma to be established holds for $k$ and for all the
appropriate $E\ssp,f\aar 1\ssp$, we prove \œ$\ssp\all{k\in\No}\,
\roman P\ssp(k\spp)$ by induction as follows. The case \œ$\ssp k=\emptyset$ or
\œ$\ssp k=1\adot$ already having been established, we assume $\ssp
\roman P\ssp(k\spp)\,$, and let $\sp g$ and $\sp g\ar 1$ be as above. Using
Proposition \ref{D^k recursion} and assuming that \œ$\tilde f \in
\CalDSeip^{k\ssp+\sp 1.}\sp$, we have \œ$\tilde f\in\CalDSeip^{1.}$ and \œ$
\varSe\tilde f\in\CalDSeip^k\,$. From \œ$\varSe\tilde f\in\CalDSeip^k$ we get
\œ$(E\spp\skcap E\spp\skcap F\sp,G\sp,g\ar 1\KN1\ywed\sp)\in\CalDSeip^k\,$,
whence by $\ssp\roman P\ssp(k\spp)$ it follows that \œ$\varSe\tilde f\aar 1 =
(\sp\sigrd\varSe\tilde f\aar 1\sp,g\ar 1)$ \œ$\in\CalDSeip^k\,$. Since by \œ$
\tilde f\in\CalDSeip^{1.}$ and $\ssp\roman P\ssp(\sp 1\adot)$ we have \œ$
\tilde f\aar 1\in\CalDSeip^{1.}\,$, then \œ$\tilde f\aar 1 \in
\CalDSeip^{k\ssp+\sp 1.}$ follows. Conversely, letting \œ$\tilde f\aar 1 \in
\CalDSeip^{k\ssp+\sp 1.}\sp$, we have \œ$\tilde f\aar 1\in\CalDSeip^{1.}$ and
(x) that \œ$\tilde f\aar 1\ssp,\sp\varSe\tilde f\aar 1\in\CalDSeip^k\,$. From
(x) by $\ssp\roman P\ssp(k\spp)$ it follows that \œ$\tilde f\ssp,\sp
(E\spp\skcap E\spp\skcap F\sp,G\sp,g\ar 1\KN1\ywed\sp)\in\CalDSeip^k\,$. Since
with \œ$H=$ \œ$E\spp\skcap F\skcap(E\spp\skcap F\ssp)$ and \œ$\ell\ar 0 =
\seq{\,(\sp x\ssp,v\ar 1):\smb W=(\sp x\ssp,v\,;u\ssp,v\aar 1)\in\vecs H\sp\,}
$ and \œ$\ell\ar 1=\langle\,(\sp x\ssp,u\ssp,$\linebreak \œ$v\sp):\smb W =
(\sp x\ssp,v\,;u\ssp,v\aar 1)\in\vecs H\sp\,\rangle$ we have the continuous
linear maps \œ$(\spp H\sp,E\spp\skcap F\sp,\ell\ar 0)$ and \œ$(\spp H\sp,
E\spp\skcap E\spp\skcap F\sp,\ell\ar 1)$ with the above established \œ$g =
\ssigrd G\circ[\,\sp g\ar 1\KN1\ywed\sn\circ\spp\ell\ar 1\ssp,\sp
f\circ\sp\ell\ar 0\,]\,$,\linebreak by (\ref{pr1b})\ssp, (\ref{pr1c}) and
(\ref{pr1d}) of Proposition \ref{basic D^k properties} it follows that \œ$
\varSe\tilde f = (\sp\sigrd\varSe\tilde f\sp,g\sp)\in\CalDSeip^k\,$. Since by
$\tilde f\aar 1\in\CalDSeip^{1.}$ and $\ssp\roman P\ssp(\sp 1\adot)$ we have $
\tilde f \in \CalDSeip^{1.}\,$, then $\tilde f\in\CalDSeip^{k\ssp+\sp 1.}$
follows.                                                         \end{proof}

\begin{corollary}\label{varSe D^k iff derSe D^k}

Let $\ssp\tilde f\in\CalDSeip^{1.}\,$. Then $\ssp\varSe\tilde f\in\CalDSeip^k\sp
\equivv\sp\DerSe\tilde f\in\CalDSeip^k\ssp$ for all $\ssp k\in\No\ssp$.\KN{10}

  \end{corollary}

\begin{proof} Applying Lemma \ref{f(,) D^k imply f^v D^k} above with $E$ in
place of $F$ and $F$ in place of $G$ and $\DerSe\tilde f$ in place of $
\tilde f\aar 1$ and $\varSe\tilde f$ in place of $\tilde f\sp$, we see that
under the assumption that $\tilde f\in\CalDSeip^{1.}$ we have $\varSe\tilde f
\in\CalDSeip^k$ if and only if $\DerSe\tilde f\in\CalDSeip^k\,$. \end{proof}

\begin{theorem}\label{our equivv Seip}

Let $\tilde f=(E\ssp,F\sp,f\sp)\,$. Then $\tilde f \in
\Cal D_{\fiveroman{Seip}}^{\,k}\equivv$ {\ssp\rm(a)} or {\ssp\rm(b)} or
{\ssp\rm(c)} or {\ssp\rm(d)} for\vskip1mm

{\rm(a)} \ \ $k=\emptyset$ and $E\ssp,F\in\scgVS(\bosy R\sp)$ and \newline $\null$ \hfill
$f\in F^{\sp/\sp E}$ and $\dom\sn f\in\taurd E$ and
                       $f$ continuous $\taurd E\to\taurd F\sp,$\KP8\vskip.5mm

{\rm(b)} \ \ $k=1\adot$ and $f:E\supset\dom\sn f\to F$ is {\ssp\rm stetig
differenzierbar} \newline $\null$ \hfill
in the sense of {\,\rm\cite[Definition 4.2, p.\ 60]{Se72}\ssp,}\KP8\vskip.5mm

{\rm(c)} \ \ $k\in\N\snn\setminus\sn 2\adot$ and $f:E\supset\dom\sn f\to F$ is
$k\,${\rm--\,mal stetig differenzierbar} \newline $\null$ \hfill
 in the sense of {\,\rm\cite[Definition 5.1, p.\ 75]{Se72}\ssp,}\KP8\vskip.5mm

{\rm(d)} \ \ $k=\infty$ and $f:E\supset\dom\sn f\to F$ is {\ssp\rm unendlig
 oft stetig differenzierbar} \newline $\null$ \hfill in the sense of
 {\,\rm\cite[Definition 6.1, p.\ 85]{Se72}\ssp}.\KP8

  \end{theorem}

\begin{proof} First note that \œ$\sp\CalDSeip^k=\emptyset\sp$ if \œ$\sp k
\not\in\infty\ssp\yplus\sp$. Then we see that the assertion concerning (a) is
trivial, and that for (b) it is given in
Corollary \ref{D^1 iff BGN^1 iff Seip^1(1972)}\ssp(a) above. For a moment
supposing we already also know (c) we obtain (d) by \œ$\sp\CalDSeip^{1.} \inc
\CalDSeip^{0.}$ and $\sp\CalDSeip^\infty=\{\,\tilde f:\all{k\in\No}\,\tilde f
\in\CalDSeip^k\ssp\}$ obtained in Proposition \ref{D^k recursion} above.\vskip.3mm

Now, for (c) considering the classes \vskip.3mm\centerline{$
\Cal D\sbi k = \{\,\tilde f:\exi{\bosy f}\,(\sp\emptyset\,,\tilde f\,) \in
\bosy f\in\Univ\,^{k\ssp+\sp 1.}$ and $\ssp\all{i\in k}\,
\bosy f\sp\fvalue i\ssp\yplus = \DerSe(\sp\bosy f\sp\fvalue i\sp)\ssp\}\,$,}\vskip.3mm

\noin since by the above we have \œ$\sp\Cal D\yr 1\KN{1.8}\subtext{Seip\,72}
\overset{_*}= \{\,\tilde f:\tilde f\text{ is Seip\,--\,differentiable}\sp\,\}\,
$, in view of our definition of the Seip derivative, it suffices to prove that
\œ$\sp\CalDSeip^k=\Cal D\sbi k$ for all $k\in\N\sp$. By $\overset{_*}=$
already having the case \œ$k=1\adot\ssp$, assuming that with \œ$\sp k\in\N\sp$
we have \œ$\sp\CalDSeip^k=\Cal D\sbi k\,$, we prove that also \œ$\sp
\CalDSeip^{k\ssp+\sp 1.}=\Cal D\sbi{k\ssp+\sp 1.}$ as follows. By
Proposition \ref{D^k recursion}\ssp(a) and
Corollary \ref{varSe D^k iff derSe D^k} above, for all $\tilde f$ we have\vskip.3mm

\centerline{$\tilde f\in\CalDSeip^{k\ssp+\sp 1.}\sp\equivv\sp
\tilde f\in\CalDSeip^{1.}$ and $\varSe\tilde f\in\CalDSeip^k\sp\equivv\sp
\tilde f\in\CalDSeip^{1.}$ and $\DerSe\tilde f\in\CalDSeip^k$}\vskip.3mm

$\mhyppy{23}\equivv\sp\tilde f\in\Cal D\sbi{1.}$ and $\DerSe\tilde f \in
\Cal D\sbi k\,$, \hfill and as in the proof of Proposi-\vskip.5mm

\noin tion \ref{BGN recurs} it is seen that \ \ $\tilde f\in\Cal D\sbi{1.}$
and $\DerSe\tilde f\in\Cal D\sbi k\sp\equivv\sp
\tilde f\in\Cal D\sbi{k\ssp+\sp 1.}\,$.                          \end{proof}

As a consequence of Seip's non-uniform manner of putting his definitions, the
proof of our Theorem \ref{our equivv Seip} above became a bit clumsy. See in
particular \cite[Definition 5.1, p.\ 75]{Se72} whose clarity surely is not the
best possible. \hfill We interpreted the vague presentation there so that in
Theorem \ref{our equivv Seip} we have the equivalence\vskip.3mm\centerline{%
(c)\ $\equivv\,k\in\N\snn\setminus\sn 2\adot$
                   and $\ssp(E\ssp,F\sp,f\ssp)\in\Cal D\sbi k\,$.}

\begin{example}\label{Seip d:fm group}

We prove that \œ$\tilde f\in\CalDSeip^\infty$ for the map \œ$\tilde f =
(E\skcap E\,,\kappatv\KN1\fvalue E\,,\sp f\ssp)\,$, where with \œ$E =
\Cal D\ssp(\Re\sp)$ and \œ$F=\Cinfty(\Re\sp)$ and $\Iota=\roman{id\,\ssp}\Re\sp\,
$, we have \œ$f=f\aar 1\ssp|\,(\sp\vecs E\ssp)^{\sp\times 2.}$ for the
function $f\aar 1=\seq{\,x\circ(\sp\Iota+y\sp):z=(\sp x\ssp,y\sp)\in(\sp\vecs
F\ssp)^{\sp\times 2.}\ssp}\,$.

We prove first that \œ$(\spp F\sp\sqcap F\sp,F\sp,f\aar 1)\in\CalDSeip^\infty\,
$. Note that since $F$ is a Fr\'echet space, it is Seip\,--\,convenient and \œ$
F\sp\sqcap F=F\sp\skcap F$ holds. Letting $G^{\,H}$ denote Seip's canonical
function space $\Cal D_{\fiveroman{Seip}}^{\,\sp\infty}(\sp\vecs H\sbi H\ssp,
G\ssp)\,$, which in \cite[Kapitel 6]{Se72} would be written \œ$\ssp
\roman D^{\,\infty}\spp(\spp H\supset\vecs H\sp,G\ssp)\,$, we hence wish to
have \œ$f\aar 1\in\vecs\big(\big(\sp
\biit R\,\sp^{\bmii7R}\ssp\big)\sp^{F\,\sqcap\ssp F}\sp\big)\,$. By
\cite[Satz 6.18, pp.\ 91\,--\,92]{Se72} this holds if{}f \œ$\ywed\sn f\aar 1
\in\vecs\big(\sp\bosy R\,\sp^{\bosy R\sp\,\sqcap\ssp(F\,\sqcap\ssp F\sp)}\sp
\big)\,$. For \œ$t\in\Re$ and \œ$x\ssp,y\in\vecs F$ we have \œ$
\ywed\sn f\aar 1\KN{.7}\fvalue(\sp t\,;x\ssp,y\sp) =
f\aar 1\KN{.7}\fvalue(\sp x\ssp,y\sp)\fvalue t =
x\fvalue(\sp t+(\sp y\fvalue t\sp))\,$, whence we see that $\ywed\sn f\aar 1$
has\linebreak the decomposition\vskip.3mm\centerline{$
(\sp t\,;x\ssp,y\sp)\mapsto(\sp t\ssp,y\ssp,x\sp)\mapsto
(\sp t\ssp,y\fvalue t\ssp,x\sp)=(\sp t\ssp,s\ssp,x\sp)\mapsto
(\sp t+s\ssp,x\sp)=(\sp r\sp,x\sp)\mapsto x\fvalue r\,$,}\vskip.3mm

\noin and consequently that for \œ$\ywed\sn f\aar 1 \in
\vecs\big(\sp\bosy R\,\sp^{\bosy R\sp\,\sqcap\ssp(F\,\sqcap\ssp F\sp)}\sp\big)
\,$, letting \œ$\Eps =\roman{ev}\,|\,(\Re\snn\times\sn(\sp\vecs F\ssp))\,$, it
suffices that \œ$\Eps\in
\vecs\big(\sp\bosy R\,\sp^{\bosy R\sp\,\sqcap\ssp F}\sp\big)\,$. This in turn
follows if we establish \œ$\yvee\Eps \in
\vecs\big(\big(\sp\biit R\,\sp^{\bmii7R}\ssp\big)\sp^F\sp\big)$ \œ$ =
\vecs\big(\spp F\,\sp^F\sp\big)\,$. For \œ$t\in\Re$ and \œ$x\in\vecs F$ having
\œ$\yvee\Eps\fvalue x\fvalue t = \Eps\fvalue(\sp t\ssp,x\sp) = x\fvalue t\ssp
$, we see that $\yvee\Eps=\idv F\sp$, hence that $\yvee\Eps \in
\vecs\big(\spp F\,\sp^F\sp\big)$ is trivial, in view of
Proposition \ref{basic D^k properties}\ssp(\ref{pr1b})\ssp.

To obtain \œ$\tilde f\in\CalDSeip^\infty\,$, we deduce as follows. Letting \œ$
\roman E\sp\,l = \Cal D\ssp\sbi{[\mminus l\sp,\ssp l\,]\ssp}(\Re\sp)$ and \œ$
\roman G\sp\,l = $ $\roman E\sp\,l\sp\sqcap\roman E\sp\,l$ and \œ$
\tilde{\roman f}\sp\,l = (\sp\roman G\sp\,l\ssp,\spp\roman E\sp\,l\ssp,\sp
              f\aar 1\ssp|\,\vecs\roman G\sp\,l\ssp)\,$, note that for every \œ$
l\in\Zepp$ we have $\roman E\sp\,l$ a closed topological linear subspace both
in $E$ and in $F\sp$, and that for \œ$x\ssp,y\in\vecs\roman E\sp\,l$ and \œ$
t\in\Re\spp\setminus[\minus l\ssp,\sp l\ssp\,]$ we have \œ$
f\aar 1\KN{.7}\fvalue(\sp x\ssp,y\sp)\fvalue t =
x\fvalue(\sp t+(\ssp y\fvalue t\sp)) = x\fvalue t = 0\,$, hence \œ$
f\aar 1\KN{.6}\fvalue(\sp x\ssp,y\sp)\in$ $\vecs\roman E\sp\,l\ssp$, and
further \œ$\rng\taurd\tilde{\roman f}\sp\,l\inc\vecs\roman E\sp\,l\ssp$. From
\œ$(\spp F\sp\sqcap F\sp,F\sp,f\aar 1)\in\CalDSeip^\infty\,$ by
Proposition \ref{Seip iff simply Seip} it follows by induction on \œ$k\in\No$
that $\varSei^kÏ\tilde{\roman f}\sp\,l$ is global with \œ$
\varSei^kÏ\tilde{\roman f}\sp\,l\in\CalCSe0$ for all \œ$l\in\Zepp$. Using this
and noting that every $\taurd E\,$--\,compact $K$ is $\taurd\roman E\sp\,l\,
$--\,compact for some \œ$l\in\Zepp$, it follows that $\varSei^kÏ\tilde f$ is
global with \œ$\varSei^kÏ\tilde f\in\CalCSe0$ for \œ$k\in\No\ssp$, whence by
Proposition \ref{Seip iff simply Seip} we get $\tilde f\in\CalDSeip^\infty\,$.

  \end{example}

Seip systematically used the concept of being {\it scharf differenzierbar\ssp}
in \cite{Se72} when stating the premises in his various inverse and implicit
function theorems. The next example shows that this property is so strong that
these theorems are practically useless, at least when considering problems
where maps of the form \œ$x\mapsto\varphi\circ x$ between spaces of smooth
functions are involved. Indeed, the only maps below of this form which are
everywhere scharf differenzierbar are those given by some affine $\,\varphi :
t \mapsto\alpha\sp\,t+\beta\,$ with fixed $\,\alpha\,,\sp\beta\in\Re\sp\,$.

\begin{example}\label{scharf d:bar}

We fix a smooth nonaffine bijection \œ$\varphi:\Re\to\Re$ with \œ$0\not\in\rng
\varphi\ssp'\sp$, and with \œ$I=[\,0\,,1\,]$ and \œ$G=\Cinfty(I\sp)$ and \œ$f=
\seq{\,\varphi\circ x:x\in\vecs G\sp\,}\,$, consider the map \œ$\tilde f =
(\spp G\sp,G\sp,f\sp)\,$. Since $G$ is Fr\'echet, it is Seip\,--\,convenient
and \œ$\CinftyPi(\biit R\ssp)\,|\,\{\ssp(\spp G\sp,G\sp)\ssp\}\inc$ $
\CalDSeip^\infty\,$, whence by \cite[Remarks 3.7\ssp(a)]{Hic} it follows that
$\tilde f\ssp,\sp(\spp G\sp,G\sp,f\sp\inve\spp)\in\CalDSeip^\infty\,$.

Note that we cannot obtain Seip\,--\,smoothness of $\tilde f$ and its inverse
by (directly) applying any exponential law in Seip's theory since the closed
interval $I$ is not an admissible domain there. However, one can prove that
there is a continuous linear map, extension operator \œ$\Eps:G\to E =
\Cinfty(\Re\sp)$ with \œ$x\inc\Eps\sp\fvalue x$ for \œ$x\in\vecs G\sp$. With
the aid of Seip's exponential law, we obtain smoothness \œ$E\to E$ of the map
\œ$f\aar 1\sn:x\mapsto\varphi\circ x\ssp$. Since also \œ$\rho:E\to G$ given by
\œ$y\mapsto y\,|\,I$ is a continuous linear map, in view of $\,f =
\rho\sp\circ f\aar 1\snn\circ\spp\Eps\,$ we can get $\,\tilde f \in
\CalDSeip^\infty\,$, and similarly for the inverse.

We next show indirectly that for any fixed \œ$\xi\in\Re$ such that \œ$
\varphi\ssp''(\xi)\not=0\,$, when we take \œ$x=I\snn\times\sn\{\sp\xi\sp\}\,$,
the above map $\tilde f$ is not scharf differenzierbar at the point $x\ssp$.
Indeed, suppose that $\tilde f$ were scharf differenzierbar at $x\ssp$. This
means that for\vskip.3mm\centerline{$
\varrho = \seq{\,\varphi\circ(\sp x+u\sp)-\varphi\circ x -
           \varphi\ssp'\circ x\snn\cdot\snn u:u\in\vecs G\sp\,}$}\vskip.3mm

\noin it holds (s) that for every \œ$\ssp\ell\in\Cal L\,(\spp G\sp,G\sp)$
there are open zero neighborhoods $U\spp,V$ in $G$ such that for \œ$h =
\seq{\,u+\ell\ssp\fvalue v:u\in U\text{ and }(\sp u\ssp,v\sp)\in\varrho\sp\,}$
we have a homeomorphism \œ$h:\taurd G\lei U\to\taurd G\lei V\spp$, and in
particular a bijection \œ$h:U\to V\spp$. Note that for all $u\in\vecs G$ we
have \ \ $\varrho\fvalue u = \varphi\circ(\sp x+u\sp)-\varphi\circ x -
\varphi\ssp'\circ x\snn\cdot\snn u = $\vskip.5mm\centerline{$
\int_{\,0}^{\,1}(\sp\varphi\ssp'\circ(\sp x+s\,u\sp) -
 \varphi\ssp'\circ x\sp)\snn\cdot\snn u\ssp\,\roman d\,s =
\int_{\,0}^{\,1}\int_{\,0}^{\,1}s\,\sp\varphi\ssp''\circ(\sp x + s\ar 1\ssp
s\,u\sp)\snn\cdot\snn u^{\,2}\,\roman d\,s\ssp\,\roman d\,s\ar 1\,$.}\vskip.5mm

We now consider \œ$\ssp\ell=\seq{\,u\ssp':u\in\vecs G\sp\,}$ in condition (s)
and proceed to show that $h$ is not even injective. Fixing \œ$\eps\in\Rep$
sufficiently small, for \œ$u\ar 0=I\snn\times\sn\{\sp\eps\sp\}$ we have $
u\ar 0\in U\spp$, and by $\varphi\ssp''(\xi)\not=0$ we may also arrange that\vskip.5mm\noin$
(*)$ \hfill $\int_{\,0}^{\,1}\int_{\,0}^{\,1}(\sp
  s\ar 1\ssp s^{\,2}\sp\varphi\ssp'''(\sp\xi+s\ar 1\ssp s\,\eps\sp)\,\eps
+ 2\,s\,\sp\varphi\ssp''(\sp\xi+s\ar 1\ssp s\,\eps\sp))
   \,\roman d\,s\ssp\,\roman d\,s\ar 1\not=0\,$. \hfill \phantom{$(*)$}\vskip.5mm

\noin We have $u\ar 0+(\sp\varrho\fvalue u\ar 0)\ssp'=u\ar 0\,$, and defining
the smooth $\chi:I\snn\times\sn(\Re\times\sn\Re\ssp)\to\Re\,$ by\vskip.5mm$\mhyppy{12.1}
(\sp t\,;\eta\ssp,\eta\sp\yplk) \mapsto
 \eta - \eps + \int_{\,0}^{\,1}\int_{\,0}^{\,1}(\sp
  s\ar 1\ssp s^{\,2}\sp\varphi\ssp'''(\sp\xi+s\ar 1\ssp s\,\eta\sp)\,
   \eta\sp\yplk\ssp\eta^{\,2}$\nopagebreak\par\nopagebreak$\mhyppy{53}
+ 2\,s\,\sp\varphi\ssp''(\sp\xi+s\ar 1\ssp s\,\eta\sp)\,
   \eta\,\eta\sp\yplk\spp)\,\roman d\,s\ssp\,\roman d\,s\ar 1\,$,\vskip.5mm

\noin it now holds that \œ$u\ar 0\in V_{\spp}$, and $u\ar 0$ satisfies the
differential equation \œ$\chi\circ[\,\sp\roman{id}\,;u\ssp,u\ssp'\sp\,]=$ $
I\snn\times\sn\{0\}$ with initial condition $u\fvalue 0=\eps\,$. More
generally, we have\vskip.3mm\centerline{$
h\inve\image\snn\{\ssp u\ar 0\} = U\cap\{\,u:\chi\circ[\,\sp\roman{id}\,;
u\ssp,u\ssp'\sp\,] = I\snn\times\sn\{0\}\sp\}\,$.}\vskip.3mm

\noin Consequently, we are done once we show that there is at least one \œ$
u\in U\setminus\sn\{\ssp u\ar 0\}$ satisfying the differential equation $
\chi\circ[\,\sp\roman{id}\,;u\ssp,u\ssp'\sp\,]=I\snn\times\sn\{0\}\,$.

To establish this, we utilize our previous result \cite[Theorem 5.2]{Hic} from
which it follows that if \œ$0\not\in\{\,\partial\ar 3\ssp\chi\fvalue(\sp t\,;
\eps\ssp,0\sp):t\in I\sp\,\}\,$, there is $\gamma$ with \œ$(\biit R\,,G\sp,
\gamma\sp)\in\CinftyPi(\biit R\ssp)$ and \œ$(\sp\eps\ssp,u\ar 0)\in\gamma\,$,
and such that we have \œ$\chi\circ[\,\sp\roman{id}\,;u\ssp,u\ssp'\sp\,]=
I\snn\times\sn\{0\}$ and \œ$u\fvalue 0=\eta$ whenever \œ$(\sp\eta\ssp,u\sp)\in
\gamma\,$. Then $\gamma$ is in particular continuous \œ$\tau_{_{I\!\!R}}\to
\taurd G\sp$, whence for all $\eta$ sufficiently close to $\eps$ we have \œ$
\gamma\fvalue\eta=u\in U\spp$, and here \œ$u\not=u\ar 0$ if we take \œ$\eta
\not=\eps$ since $u\ar 0\KN1\fvalue 0=\eps\not=\eta=u\fvalue 0\,$. By $(*)$ we
indeed have $0\not\in\{\,\partial\ar 3\ssp\chi\fvalue(\sp t\,;\eps\ssp,0\sp):
t\in I\sp\,\}\,$.

Note also that essentially by the same method as above, we could have proved
more generally that if \œ$x\in\vecs G$ only is such that \œ$0\not\in
\varphi\ssp''\ssp[\,\rng x\,]\,$, then $\tilde f$ is not scharf
differenzierbar at $x\ssp$. In this case the function \œ$p\ar 3 =
\seq{\,\partial\ar 3\ssp\chi\fvalue(\sp t\,;\eps\ssp,0\sp):t\in I\sp\,}$ is
not necessarily any more constant but we may still arrange \œ$0\not\in\rng
p\ar 3$ by choosing \œ$\eps>0$ sufficiently small. For example, if we take \œ$
\varphi=\seq{\,t+\roman e^{\,t}:t\in\Re\ssp\,}\,$, then $\tilde f$ is not
scharf differenzierbar at any $x\ssp$.                         \end{example}

In the preceding example, we did not make explicit the space in which the
function space integrals are taken in the expression for $\varrho\fvalue u\ssp
$. Most conveniently, they are understood to be in $\sp\biit R\expnota^I]_{tvs}\,
$, but they could have been considered also in $G$ or $\sp
\biit R\expnota^I]_{tvs}{\sn_{/\sp\upsilon_s G}}\ssp$, \,of which the latter
is not sequentially complete.

\begin{example}\label{"norm"}

Let \œ$E=(X\sp,\Cal T\,)\in\LCS(\biit R\ssp)\,$, and let $\Cal C$ be the set
of bounded abso- lutely convex closed sets in $E\,$. If one wants to properly
present the exact content of \cite[Definition 8.6, p.\ 108]{Se72}\ssp, tracing
the matters back to
\cite[Definitions 0.10, 8.3, pp.\ 4\,--\,5, 107\,--\,108]{Se72} one arrives at
the conclusion that a {\it full Seip norm\ssp} in $E$ is
(cf.\ Remark \ref{Weihe pairs} above) for example some \œ$
\bosy\nu=\bar n\weco\bar\varphi\ar 1\weco\bar\varphi\ar 2$ where $\bar n$ is a
functor between certain small \q{additive} categories (it not being worth
while specifying them here) and the $\bar\varphi\sbi\iota$ are certain natural
transformations between certain associated functors. However, it turns out
that the $\bar\varphi\sbi\iota$ once existing are unique, and that $\bar n$ is
uniquely determined by its object component $n$ for which we have (n) that \œ$
n\in\Cal C\,^{\upsilon_s\sp E}$ with \œ$
z\in n\fvalue z\inc\ClT(\sp n\fvalue x+n\fvalue(\sp z-x\sp))$ for all \œ$
x\ssp,z\in\vecs E\,$. Conversely, every such $n$ determines a unique full Seip
norm $\bosy\nu\ssp$. Hence, if one wants to keep matters as simple as
possible, it is advisable to abandon the redundant \cite[Definition 8.6]{Se72}
and in- stead define a {\it Seip norm\ssp} in $E$ to be any $n$ satisfying (n)
  above.

To give a nontrivial example of a Seip norm, taking \œ$
E = \biit R\expnota^\Re]_{tvs}$ and \œ$ n = $\linebreak \œ$
\seq{\,\Re\sp\,^{\Re}\cap\{\,z:\all{s\in\Re}\,|\,z\fvalue s\,| \le
|\,x\fvalue s\,|\,\}:x\in\vecs E\sp\,}\,$, then $n$ is a Seip norm in $E\,$.
Note that from the triangle inequality \œ$
|\,z\fvalue s\,|\le|\,x\fvalue s\,|+|\,z\fvalue s-x\fvalue s\,|$ it easily
fol- lows that we even have $z\in n\fvalue z\inc
        n\fvalue x+n\fvalue(\sp z-x\sp)$ for $x\ssp,z\in\vecs E\,$.

Supposing that also \œ$\varphi\in\Re\sp\,^{\Re\ssp\times\ssp\Re}$ with \œ$
|\,\varphi\fvalue(\sp s\ssp,t\ar 1) - \varphi\fvalue(\sp s\ssp,t\ar 2)\,| \le
\frac12\,|\,t\ar 1\snn - t\ar 2\ssp|$\linebreak for all \œ$
s\ssp,t\ar 1\ssp,t\ar 2\in\Re\sp\,$, and considering \œ$
f=\seq{\,\varphi\circ[\,\roman{id\,},x\,]:x\in\vecs E\sp\,}\,$, then we have \œ$
n\fvalue(\sp f\sp\fvalue x-f\sp\fvalue y\sp)\inc\frac 12\,n\fvalue(\sp x-y\sp)\,$
for all \œ$\ssp x\ssp,y\in\vecs E\,$. Consequently, by
\cite[Satz 8.8, p.\ 109]{Se72} there is $x$ with \œ$
\dom\snn(\sp f\sp\cap\id\sn)=\{\sp x\sp\}\,$. However, the same result can be
obtained directly from Banach's fixed point theorem by considering for each
fixed $s\in\Re\,$ the $\,
\{\ssp(\sp s\ssp,t\ssp,|\,s-t\,|\ssp):s\ssp,t\in\Re\sp\,\}\,$--\sp\,contractor
$\,\varphi\,(\sp s\ssp,\sp\cdot\sp)\,$.                        \end{example}


\end{document}